%% file: 00_main.tex
\documentclass[pdflatex,sn-mathphys-num,iicol]{sn-jnl}

\usepackage{graphicx}%
\usepackage{multirow}%
\usepackage{amsmath,amssymb,amsfonts}%
\usepackage{amsthm}%
\usepackage{mathrsfs}%
\usepackage[title]{appendix}%
\usepackage{xcolor}%
\usepackage{textcomp}%
\usepackage{manyfoot}%
\usepackage{booktabs}%
\usepackage{algorithm}%
\usepackage{algorithmicx}%
\usepackage{algpseudocode}%
\usepackage{listings}%
\usepackage{matlab-prettifier}
\usepackage{cleveref}
\usepackage{anyfontsize}
\usepackage{tikz}
\usepackage{pgfplots}
\pgfplotsset{compat=newest}
\pgfplotsset{grid style={dotted,gray}}
\usetikzlibrary{plotmarks,pgfplots.fillbetween,shapes.geometric, arrows.meta,patterns}

\newcommand{\R}{\mathbb{R}}
\newcommand{\N}{\mathbb{N}}
\newcommand{\mx}{\mathbf{x}}
\newcommand{\mK}{\mathbf{K}}
\newcommand{\mF}{\mathbf{F}}
\newcommand{\mU}{\mathbf{U}}
\newcommand{\dd}{\,\mathrm{d}}
\DeclareMathOperator{\dist}{dist}
\DeclareMathOperator*{\argmin}{arg\,min}

\theoremstyle{thmstyleone}%
\theoremstyle{thmstyletwo}%
\newtheorem{remark}{Remark}%

\theoremstyle{thmstylethree}%

\begin{document}

\title[A 140 line MATLAB code for topology optimization problems with probabilistic parameters]{A 140 line MATLAB code for topology optimization problems with probabilistic parameters}


\author*[1]{\fnm{Andrian} \sur{Uihlein}}\email{andrian.uihlein@fau.de}
\author[2]{\fnm{Ole} \sur{Sigmund}}\email{olsi@dtu.dk}
\author[1]{\fnm{Michael} \sur{Stingl}}\email{michael.stingl@fau.de}

\affil[1]{\orgdiv{Department of Mathematics}, \orgname{Friedrich-Alexander-Universität Erlangen-Nürnberg}, \orgaddress{\street{Cauerstraße}, \city{Erlangen}, \postcode{91058}, \country{Germany}}}

\affil[2]{\orgdiv{Department of Civil and Mechanical Engineering}, \orgname{Technical University of Denmark}, \orgaddress{\street{Nils Koppels Allé}, \city{Kongens Lyngby}, \postcode{2800}, \country{Denmark}}}



\abstract{We present an efficient 140 line MATLAB code for topology optimization problems that include probabilistic parameters. It is built from the \texttt{top99neo} code by Ferrari and Sigmund and incorporates a stochastic sample-based approach. Old gradient samples are adaptively recombined during the optimization process to obtain a gradient approximation with vanishing approximation error. The method's performance is thoroughly analyzed for several numerical examples. While we focus on applications in which stochastic parameters describe local material failure, we also present extensions of the code to other settings, such as uncertain load positions or dynamic forces of unknown frequency. The complete code is included in the Appendix and can be downloaded from \url{www.topopt.dtu.dk.}}




\maketitle

\section{Introduction}
Topology optimization is a highly active field of research in optimization and engineering alike. Thus, it is hardly surprising that an overwhelming amount of free as well as commercial software for efficient topology optimization is available, e.g.,~\cite{Wang2003,Andreassen2014,Talischi2012,Suresh2010,Challis2009,Xia2015,Liu2005,Sanders2018,Allaire2006,Sok2010,Amir2013,Liu2014}. The field was given an extensive review in~\cite{Wang2021}. The most popular implementations for linear elasticity in MATLAB are the famous \texttt{top99}~\cite{99lineCode} and its descendants \texttt{top88}~\cite{88lineCode} and \texttt{top99neo}~\cite{99neoCode}, which we build upon in this contribution.

With such a repertoire of strong numerical tools, the complexity of optimization tasks has grown rapidly over time. In recent years, more and more focus has been placed on topology optimization under uncertainty, i.e., optimization problems in which some of the involved parameters are unknown and need to be modeled probabilistically. Such problems arise naturally in many applications, for example when the forces applied to the structure are unknown~\cite{Liu2018,Nishioka2023,Maute2014}. Likewise, stochastic models can be used to describe local material impurities~\cite{De2022,Long2018,Greifenstein2020,Jansen2013,Zhou2016} or geometry imperfections~\cite{De2020,Krger2024,Komini2023,Schevenels2011} of the structure, as they may appear in the context of additive manufacturing.

In the literature, many different schemes have been proposed to deal with model uncertainties of different kinds. For example, a natural approach is to perform a worst-case optimization. Here, the main challenge is to design schemes that efficiently find (or approximate) this worst-case, e.g.,~\cite{Brittain2011,Gao2021,Kharmanda2004,Lombardi1998,Zhou2016}. Alternatively, it is common to consider a combination of average and standard deviation over all possibilities, since there exist powerful tools from statistics for this setting~\cite{Gao2022,Krger2024,Tootkaboni2012,Schevenels2011}. Lastly, there are several fully stochastic optimization approaches, in which stochasticity is dealt with by different sampling techniques~\cite{De2020,De2021,De2022,Campi2010}. Note that it is also possible to use stochastic optimization schemes to solve deterministic problems, as it is often observed that noisy gradient steps can help the optimizer escape poor local minima in the objective function~\cite{pmlr-v80-kleinberg18a}. In this work, however, we focus on optimization problems in which stochasticity is part of the model itself and not (artificially) added by the optimization algorithm.

The method proposed in this contribution follows this sample-based philosophy and is closely related to the famous stochastic gradient scheme~\cite{Monro1951}. To be precise, we adapt the techniques proposed in~\cite{CSGoriginal}. It is important to note that this method stores sampled gradient values to adaptively recombine old and new information for the current gradient approximation. It can be shown that the associated approximation error vanishes over the course of iterations~\cite{CSG1}, resulting in an optimization scheme that aims to minimize the number of times we need to solve the system state equation during the optimization process.

Our goal is to provide an easy-to-use educational tool for students and newcomers to the field. For maximum accessibility, the implementation is designed and presented as an add-on to the original 88/99 line codes. Therefore, large parts of the code are identical to \texttt{top99neo} and we use similar notation as in~\cite{88lineCode,99lineCode,99neoCode}. Ideally, users which are familiar with these codes will be able to use \texttt{topS140} without difficulty.

The paper is structured as follows: We briefly introduce the probabilistic topology optimization problem in~\Cref{sec:ProblemFormulation}. In~\Cref{sec:Method}, we present the theoretical background of the optimization approach. Here, we place particular emphasis on the introduction and motivation of the stochastic sample-based approximation techniques of~\cite{CSG1}. The code itself is explained in~\Cref{sec:Implementation}, where we use an exemplary optimization problem to motivate individual steps. Afterwards, we analyze the code's performance in several numerical examples in~\Cref{sec:NumExamples}. Lastly, we discuss possible variations of the code in~\Cref{sec:Generalizations}. The full code is printed in~\Cref{sec:AppendixCode,sec:AppendixCode_load} and can be dowlnoaded from~\cite{Zenodo} as well as \url{www.topopt.dtu.dk}.
\section{Problem formulation}\label{sec:ProblemFormulation}
We consider a two-dimensional rectangular design domain $\Omega$ and associated discretization $\Omega_h$. Specifically, $\Omega_h$ consists of $m$ square finite elements $\Omega_e$, ${e=1,\ldots,m}$. We denote by ${\hat{\mx}_e\in[0,1]}$ the physical material density on $\Omega_e$ and collect all information in the physical design vector ${\hat{\mx}\in[0,1]^m}$. To be precise, $\hat{\mx}$ is obtained from the pseudo-density vector ${\mx\in[0,1]^m}$ (design vector) by a two-step procedure: First, given a filter radius ${r_\text{min}>0}$, we apply the linear filter
\begin{equation}\label{eq:Filtering}
    \tilde{x}_e := \frac{\sum_{i=1}^m x_i H_{e,i}}{\sum_{i=1}^m H_{e,i}},
\end{equation}
where
\begin{equation*}
    H_{e,i} := \max\big\{ 0,r_\text{min}-\dist(\Omega_i,\Omega_e)\big\}.
\end{equation*}
Afterwards, we project each element via a relaxed Heavyside function (see~\cite{Wang2010})
\begin{equation}\label{eq:Heavyside}
    \mathcal{H}(\tilde{x}_e,\eta,\beta) := \frac{\tanh(\beta\eta)+\tanh\big(\beta(\tilde{x}_e-\eta)\big)}{\tanh(\beta\eta)+\tanh\big(\beta(1-\eta)\big)}.
\end{equation}
Lastly, we define ${\mathcal{A},\mathcal{P}_0,\mathcal{P}_1\subset\{1,\ldots,m\}}$ as disjoint sets of active design variables, passive void elements ($x_e=0$) and passive solid elements ($x_e=1$), respectively.

Based on the entries of $\hat{\mx}$, the stiffness matrix $\mK$ is constructed using a SIMP interpolation~\cite{Bendse1999}
\begin{equation}\label{eq:SIMPinterpol}
    E(\hat{x}_e) := E_\text{min} + \hat{x}_e^p(E_0-E_\text{min}),
\end{equation}
where $E_\text{min}$ and $E_0$ are the Young's moduli of void and solid. Furthermore, we allow the stiffness matrix to depend on a random parameter ${\xi\in\Xi\subset\R^{d_\Xi}}$, which is distributed according to a probability measure $\mu_\xi$. For example, in~\Cref{sec:ClampNonsym}, we model local material failure by punching a hole in the design, i.e., locally reducing the value of $\hat{\mx}$ at a small subdomain. In this setup, we may use ${\xi\in\R^2}$ to describe the (random) position of this damage region and $\mu_\xi$ to model the probability of the damage being located at a specific point, e.g., a uniform distribution ${\mu_\xi\equiv\vert\Xi\vert^{-1}}$ if each possible position has the same probability. 

Likewise, we allow the load vector $\mF$ to depend on a second random parameter ${\psi\in\Psi\subset\R^{d_\Psi}}$ with probability distribution $\mu_\psi$, which may be used to model unknown positions, directions or magnitudes of the applied forces.

The corresponding displacement $\mU_{\xi,\psi}$ is obtained by solving the state equation
\begin{equation}\label{eq:StateEquation}
    \mK(\hat{\mx},\xi)\mU_{\xi,\psi}=\mF(\psi).
\end{equation}
Lastly, we define the physical volume
\begin{equation*}
    V(\hat{\mx}) := \sum_{e=1}^m \vert\Omega_e\vert\hat{x}_e,
\end{equation*}
compliance
\begin{equation*}
    c(\hat{\mx},\xi,\psi) := \mF(\psi)^\top\mU_{\xi,\psi},
\end{equation*}
and expected compliance
\begin{equation*}
    \mathbb{E}_{\xi,\psi}\left[ c(\hat{\mx},\xi,\psi) \right] := \int_{\Xi\times\Psi} \! c(\hat{\mx},\xi,\psi) \dd(\mu_\xi\times\mu_\psi)(\xi,\psi).
\end{equation*}
Here, we use ${\dd(\mu_\xi\times\mu_\psi)}$ to denote the product measure of the individual probability distributions. For example, if $\mu_\xi$ and $\mu_\psi$ are independent and have probability densities $\rho_\xi$ and $\rho_\psi$, the expected compliance can be calculated as
\begin{equation*}
\int_\Xi\int_\Psi c(\hat{\mx},\xi,\psi) \rho_\xi(\xi)\rho_\psi(\psi)\dd\psi\dd\xi.
\end{equation*}
Now, given a volume fraction ${f\in(0,1)}$, we consider the following optimization problem:
\begin{align}\label{eq:ProblemFormulation}
    \begin{split}
        \min_{\mx_\mathcal{A}\in[0,1]^{\vert\mathcal{A}\vert}} \quad & \mathbb{E}_{\xi,\psi}\left[ c(\hat{\mx},\xi,\psi) \right], \\
        \text{s.t.} \quad & V(\hat{\mx}) \le f\vert\Omega_h\vert.
    \end{split}
\end{align}
\section{Method}\label{sec:Method}
\subsection{Optimality criterion method}
It is straightforward to calculate that the gradients of compliance and volume constraint are given by
\begin{align*}
    \nabla_{\hat{\mx}}c(\hat{\mx},\xi,\psi) &= -\big(\mU_{\xi,\psi}^\top\nabla_{\hat{\mx}}\mK(\hat{\mx},\xi,\psi)\mU_{\xi,\psi}\big)\odot \chi_\mathcal{A}, \\
    \nabla_{\hat{\mx}} V(\hat{\mx}) &= \frac{1}{m}\chi_{\mathcal{A}},
\end{align*}
where $\odot$ indicates element-wise multiplication and ${\chi_\mathcal{A}\in\{0,1\}^m}$ is the indicator vector
\begin{equation*}
    (\chi_\mathcal{A})_e=\begin{cases} 1, & e\in\mathcal{A}, \\ 0, & \text{otherwise.}\end{cases}
\end{equation*}
Thus, we have
\begin{equation*}
    \nabla_{\tilde{\mx}}\mathcal{H}(\tilde{\mx}) = \beta\frac{1-\tanh\big(\beta(\tilde{\mx}-\eta)\big)^2}{\tanh(\beta\eta)+\tanh\big(\beta(1-\eta)\big)}
\end{equation*}
and
\begin{align*}
    \nabla_\mx c(\hat{\mx},\xi,\psi) &= \nabla_{\tilde{\mx}}\mathcal{H}(\tilde{\mx})\odot \big( \nabla_\mx\tilde{\mx}^\top\nabla_{\hat{\mx}}c(\hat{\mx},\xi,\psi) \big), \\
    \nabla_\mx V(\hat{\mx}) &= \nabla_{\tilde{\mx}}\mathcal{H}(\tilde{\mx})\odot\big( \nabla_\mx\tilde{\mx}^\top\nabla_{\hat{\mx}}V(\hat{\mx}) \big).
\end{align*}
Due to the regularity of all involved functions, the full objective function gradient can be calculated as
\begin{equation}\label{eq:FullGradient}
    \nabla_\mx \mathbb{E}_{\xi,\psi}\left[ c(\hat{\mx},\xi,\psi) \right] = \mathbb{E}_{\xi,\psi}\left[ \nabla_\mx c(\hat{\mx},\xi,\psi) \right].
\end{equation}
Assuming that this gradient information is available in each iteration, problem~\eqref{eq:ProblemFormulation} can be solved using the optimality criterion method (OCM,~\cite{Bendse1995,99lineCode}). For this purpose, in iteration ${k\in\N}$, define
\begin{equation}\label{eq:OCMrule}
    \mathcal{F}_{k,e} := x_{k,e}\sqrt{-\frac{\partial_e \mathbb{E}_{\xi,\psi}\left[ c(\hat{\mx}_k,\xi,\psi) \right]}{\tilde{\lambda}_k\partial_e V(\hat{\mx}_k)}}.
\end{equation}
Here, $\tilde{\lambda}_k$ is an approximation to the true Lagrange multiplier $\lambda_k^\ast$ of the volume constraint and is obtained by a bisection subroutine for the equation
\begin{equation*}
    V\big( \hat{\mx}_{k+1}(\tilde{\lambda}) \big) -f\vert\Omega_h\vert = 0,
\end{equation*}
where ${\hat{\mx}_{k+1}(\tilde{\lambda})}$ is calculated via the element-wise design update given by
\begin{equation}\label{eq:UpdateNoML}
    x_{k+1,e} = \min\big\{ 1, \max\{ 0, \mathcal{F}_{k,e} \} \big\}.
\end{equation}
For improved stability of the method, it is common to define a move limit $\Delta_\text{move}\in(0,1)$ and replace~\eqref{eq:UpdateNoML} with
\begin{equation*}
    x_{k+1,e} = \min\big\{ \delta_+, \max\{ \delta_-, \mathcal{F}_{k,e} \} \big\},
\end{equation*}
where
\begin{align*}
    \delta_+ & := \min\{1,x_{k,e}+\Delta_\text{move}\}, \\
    \delta_- & := \max\{0,x_{k,e}-\Delta_\text{move}\}.
\end{align*}
\subsection{Stochastic approximations}
Note that the OCM update $\mathcal{F}_{k}$ in~\eqref{eq:OCMrule} requires calculation of ${\nabla_\mx \mathbb{E}_{\xi,\psi}\left[ c(\hat{\mx},\xi,\psi) \right]}$. In practice, this step is exceptionally time consuming, as evaluating the integral over ${\Xi\times\Psi}$ requires us to solve the state equation~\eqref{eq:StateEquation} for a large number of random parameters ${(\xi,\psi)}$.

In order to reduce this computational effort, we approximate the true gradient using the stochastic sample-based approach proposed in~\cite{CSGoriginal,CSG1}. For this purpose, let ${(\xi_k,\psi_k)_{k\in\N}\subset\Xi\times\Psi}$ be random samples, which we assume to be independent and identically distributed according to ${\mu_\xi\times\mu_\psi}$. Then, in each iteration, we only calculate one gradient sample
\begin{equation*}
    g_k := \nabla_\mx c(\hat{\mx}_k,\xi_k,\psi_k).
\end{equation*}
In contrast to other stochastic approximation techniques, these samples are not discarded after the iteration. Instead, all available gradient samples are combined to form an approximation $G_k$ to the current gradient
\begin{equation*}
    \nabla_\mx \mathbb{E}_{\xi,\psi}\left[ c(\hat{\mx},\xi,\psi) \right] \approx G_k := \sum_{i=1}^k \alpha_i g_i.
\end{equation*}
The coefficients ${\alpha_i\in\R_{\ge0}}$, called \emph{integration weights}, can be efficiently calculated on-the-fly, as they are based on a constant nearest neighbor surrogate model. Due to the continuity of ${\nabla_\mx c(\hat{\mx},\xi,\psi)}$, we know that the difference in gradient values
\begin{equation*}
    \big\Vert \nabla_\mx c(\hat{\mx}_1,\xi_1,\psi_1) - \nabla_\mx c(\hat{\mx}_2,\xi_2,\psi_2)\big\Vert
\end{equation*}
is small, as long as the distance between $(\hat{\mx}_1,\xi_1,\psi_1)$ and $(\hat{\mx}_2,\xi_2,\psi_2)$ is small. In order to quantify the distance between two sample points, we choose a norm ${\Vert\cdot\Vert_\ast}$ on the product space ${[0,1]^m\times\Xi\times\Psi}$. Now, unknown gradient values at the current iterate are simply approximated by the closest sample in memory: 
\begin{equation*}
    \nabla_\mx c(\hat{\mx}_k,\xi,\psi) \approx g_{j_{k,\xi,\psi}},
\end{equation*}
where
\begin{equation}\label{eq:NNIndex}
    j_{k,\xi,\psi} \in \argmin_{j=1,\ldots,k} \big\Vert (\hat{\mx}_k,\xi,\psi) - (\hat{\mx}_j,\xi_j,\psi_j)\big\Vert_{\ast}
\end{equation}
denotes the index of the closest sample point to ${(\hat{\mx}_k,\xi,\psi)}$ with respect to ${\Vert\cdot\Vert_{\ast}}$. It is important to note that, since all appearing spaces are of finite dimension, all possible choices of ${\Vert\cdot\Vert_\ast}$ are equivalent. However, the practical performance of the method does of course vary. Furthermore, some specific choices simplify the weight calculation process, see, e.g.,~\eqref{eq:NormDecomposition}.

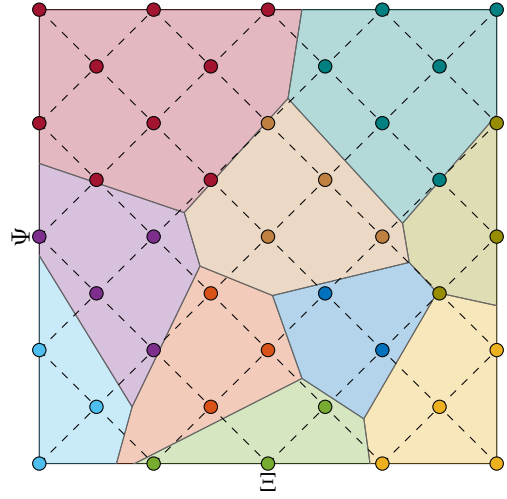
\begin{figure}
    \centering
    \input{Graphics/voronoi}
    \caption{For fixed $\hat{\mx}_k$, the integrand ${\nabla_\mx c(\hat{\mx}_k,\cdot,\cdot)}$ needs to be integrated over the space ${\Xi\times\Psi}$ (whole square). We do so by partitioning ${\Xi\times\Psi}$ into the sets $\mathcal{V}_{i,k}$ (colored polygons), on which ${\nabla_\mx c(\hat{\mx}_k,\cdot,\cdot)}$ is approximated by the piecewise constant values $g_i$. Then, the quadrature rule (dashed lines) evaluates this model on all quadrature points ${(\xi_t,\psi_t)_{t=1,\ldots,T}}$ (colored dots), to obtain the integral approximation $G_k$.}
    \label{fig:voronoi}
\end{figure}

Since the constructed nearest neighbor model is piecewise constant, the integration yields
\begin{align*}
    \nabla_\mx \mathbb{E}_{\xi,\psi}\big[ &c(\hat{\mx}_k,\xi,\psi) \big] \\ 
    & = \int_{\Xi\times\Psi} \! \nabla_\mx c(\hat{\mx}_k,\xi,\psi) \dd(\mu_\xi\times\mu_\psi)(\xi,\psi) \\
    & \approx \int_{\Xi\times\Psi} \! g_{j_{k,\xi,\psi}} \dd(\mu_\xi\times\mu_\psi)(\xi,\psi) \\
    & = \sum_{i=1}^k g_i \cdot(\mu_\xi\times\mu_\psi)\big[\mathcal{V}_{i,k}\big],
\end{align*}
with
\begin{equation*}
    \mathcal{V}_{i,k} := \left\{ (\xi,\psi)\in\Xi\times\Psi\,:\, j_{k,\xi,\psi} = i \right\}
\end{equation*}
denoting the sets on which the gradient in iteration $k$ is approximated by sample $i$. Note that, by construction the sets $(\mathcal{V}_{i,k})_{i=1,\ldots,k}$ form a partition of ${\Xi\times\Psi}$ for all ${k\in\N}$. An illustration is provided in~\Cref{fig:NN_model} and~\Cref{fig:voronoi}. Thus, the integration weights can be calculated on-the-fly and are given by the measure of the sets $\mathcal{V}_{i,k}$, i.e.,
\begin{equation*}
    \alpha_i = (\mu_\xi\times\mu_\psi)\big[\mathcal{V}_{i,k}\big].
\end{equation*}
Again, assuming that $\mu_\xi$ and $\mu_\psi$ are independent and have probability densities $\rho_\xi$ and $\rho_\psi$, this simplifies to
\begin{equation*}
\alpha_i = \iint_{\mathcal{V}_{i,k}}  \rho_\xi(\xi)\rho_\psi(\psi)\dd\psi\dd\xi.
\end{equation*}
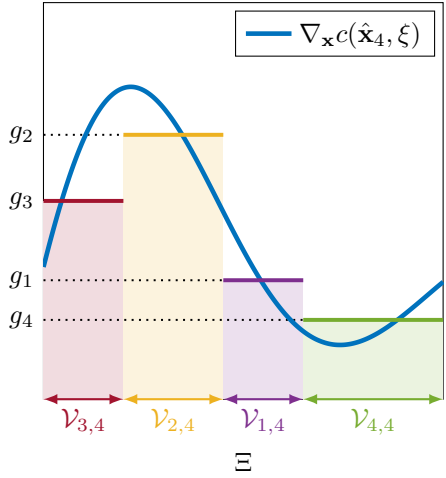
\begin{figure}
    \centering
    \input{Graphics/NN_exmpl}
    \caption{For $\Xi\subset\R$ ($x$-axis) and $k=4$, the true gradient values at the current design $\nabla_\mx c(\hat{\mx}_4,\xi)$ (blue) are approximated by the constant nearest neighbor model (colored horizontal lines). Integrating the model yields a weighted sum of the gradient samples $g_i$, where each sample is weighted by the measure of the set $\mathcal{V}_{i,4}$.}
    \label{fig:NN_model}
\end{figure}

While computing the integration weights requires additional numerical effort, a key advantage of this approach against other stochastic methods lies in the approximation property
\begin{equation*}
    \left\Vert \nabla_\mx \mathbb{E}_{\xi,\psi}\big[ c(\hat{\mx}_k,\xi,\psi) \big] - G_k\right\Vert \xrightarrow{\text{a.s.}}0\quad(k\to\infty),
\end{equation*}
see~\cite{CSG1}, which guarantees that the associated approximation error almost surely converges to zero over the course of iterations.

In practice, the approximations can efficiently be obtained by fixing an appropriate quadrature rule and evaluating the constant nearest neighbor model instead of the true gradient: 
\begin{align*}
    \int_{\Xi\times\Psi} \! \nabla_\mx &c(\hat{\mx}_k,\xi,\psi) \dd(\mu_\xi\times\mu_\psi)(\xi,\psi) \\*
    & \approx \sum_{t=1}^T w_t \nabla_\mx c(\hat{\mx}_k,\xi_t,\psi_t)\tag{Quadrature}\label{eq:Quadrature} \\*
    & \approx \sum_{t=1}^T w_t g_{j_{k,\xi_t,\psi_t}}\tag{NN model} \\*
    & = \sum_{i=1}^k \alpha_i g_i,\tag{change order}
\end{align*}
where
\begin{equation}
    \alpha_i = \sum_{t=1}^T w_t\cdot\delta_{i\vert j_{k,\xi_t,\psi_t}}\label{eq:weights_pseudo}
\end{equation}
with $\delta_{i\vert j}$ denoting the Kronecker delta
\begin{equation*}
    \delta_{i\vert j} = \begin{cases} 1, & i=j, \\ 0 & i\neq j,\end{cases}
\end{equation*}
and $w_t$ corresponding to the quadrature weight for integration point ${t\in\{1,\ldots,T\}}$. An illustration is given in~\Cref{fig:voronoi}. Since finding the nearest sample is computationally much cheaper than solving the state equation~\eqref{eq:StateEquation}, the number of quadrature points $T$ can typically be chosen rather large.
\subsubsection{Sample management}\label{sec:LimitedWeights}
As the method aggregates more and more gradient samples during the optimization process, both the computational cost of calculating the integration weights as well as the memory required to store the samples increases. To counteract this problem, we may limit the maximum number of stored gradient samples beforehand. If this limit is reached and a new sample is collected, we simply remove all information related to the sample point associated with the smallest integration weight, as it can be considered to be the least valuable information.
\section{Implementation}\label{sec:Implementation}
The \texttt{topS140} code provided in~\Cref{sec:AppendixCode} is implemented in a way that minimizes the changes to \texttt{top99neo}~\cite{99neoCode}. It is called by 
\begin{lstlisting}[frame=single,style=Matlab-editor]
topS140(nelx,nely,volfrac,penal,rmin,ft,ftBC,eta,beta,move,pnorm,maxit)
\end{lstlisting}
Here, \texttt{nelx} and \texttt{nely} are the number of finite elements in the horizontal/vertical direction. The parameters \texttt{volfrac}, \texttt{penal} and \texttt{rmin} correspond to the volume fraction $f$ in~\eqref{eq:ProblemFormulation}, $p$ in the SIMP interpolation~\eqref{eq:SIMPinterpol} and the filter radius $r_\text{min}$ in~\eqref{eq:Filtering}, respectively. The filtering technique is controlled by \texttt{ft}, with the following options:
\begin{enumerate}
    \item \texttt{ft=1}: Density filtering only.
    \item \texttt{ft=2}: Density filter and Heavyside projection~\eqref{eq:Heavyside} with parameters \texttt{beta} and \texttt{eta}. This option is not volume-preserving.
    \item \texttt{ft=3}: Density filter, Heavyside projection and adaptive modification of \texttt{eta} such that volume is preserved. This option changes the value of \texttt{eta} in each iteration. 
\end{enumerate}
The filter boundary conditions are set by \texttt{ftBC}, with \texttt{ftBC='D'} and \texttt{ftBC='N'} corresponding to zero-Dirichlet and zero-Neumann boundary conditions, respectively. Move limits and the maximum number of iterations are specified through \texttt{move} and \texttt{maxit}. Lastly, the objective in problem~\eqref{eq:ProblemFormulation} is slightly generalized by a $P$-norm approach
\begin{align*}
    \begin{split}
        \min_{\mx_\mathcal{A}\in[0,1]^{\vert\mathcal{A}\vert}} \quad & \frac{1}{P}\mathbb{E}_{\xi,\psi}\left[ c(\hat{\mx},\xi,\psi)^P \right], \\
        \text{s.t.} \quad & V(\hat{\mx}) \le f\vert\Omega_h\vert.
    \end{split}
\end{align*}
The value of $P$ is set by the parameter \texttt{pnorm}. 
\subsection{Exemplary setup}\label{sec:ExmplSetup}
\begin{figure}
    \centering
    \input{Graphics/Clamp/clamp_setup}
    \caption{Design domain with support at left and right boundary. A uniform downward facing force (red) is applied at the top. A possible damage region is indicated by the blue square.}
    \label{fig:clamp_setup}
\end{figure}
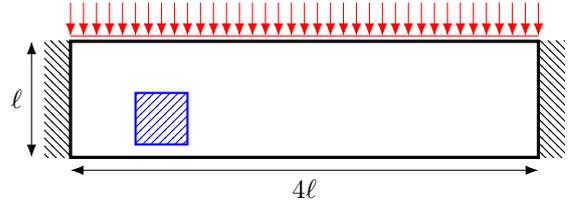

To better understand the individual steps in the algorithm, we define a generic optimization problem. Thus, consider a rectangular design domain of width $4\ell$ and height $\ell$, which is fully supported at both the left and right boundary. A uniform load is applied at the upper boundary. Moreover, following the ideas of~\cite{Jansen2013}, we model local material failure by introducing a square hole within the design region. The position of this damage region is assumed to be uniformly random distributed. An illustration of the setup can be found in~\Cref{fig:clamp_setup}.

Initialization of the code is done in 8 major steps
\begin{lstlisting}[frame=single,style=Matlab-editor,basicstyle={\fontsize{7}{8}\ttfamily}]
PRE. 1) MATERIAL AND CONTINUATION PARAMETERS
PRE. 2) DISCRETIZATION FEATURES
PRE. 3) LOADS, SUPPORTS AND PASSIVE DOMAINS
PRE. 4) DAMAGE
PRE. 5) STOCHASTIC MODEL
PRE. 6) DEFINE IMPLICIT FUNCTIONS
PRE. 7) PREPARE FILTER
PRE. 8) ALLOCATE AND INITIALIZE OTHER PARAMETERS
\end{lstlisting}
and the optimization loop consists of 7 groups of operations
\begin{lstlisting}[frame=single,style=Matlab-editor,basicstyle={\fontsize{7}{8}\ttfamily}]
RL. 1) COMPUTE PHYSICAL DENSITY FIELD
RL. 2) SETUP AND SOLVE EQUILIBRIUM EQUATIONS
RL. 3) COMPUTE SENSITIVITIES
RL. 4) SAMPLE MANAGEMENT AND INTEGRATION WEIGHTS
RL. 5) NEAREST NEIGHBOR APPROXIMATIONS
RL. 6) UPDATE DESIGN VARIABLES AND APPLY CONTINUATION
RL. 7) PRINT CURRENT RESULTS AND PLOT DESIGN
\end{lstlisting}
Due to the large overlap with \texttt{top99neo} and \texttt{top88}, we only cover the additional steps of \texttt{top99neoS} here and refer to~\cite{88lineCode,99neoCode} for details.

\subsection{Parameters for the damage region}
In this simple setting, block \texttt{PRE. 4)} is carried out in line 41:
\begin{lstlisting}[frame=single,style=Matlab-editor,basicstyle={\fontsize{7}{8}\ttfamily}]
[L, nonD, dmg_fac] = deal(20, 5, 1);
\end{lstlisting}
Since the damage is assumed to be square-shaped, it is modeled by the single parameter \texttt{L}, corresponding to the side length (in number of elements). As applying the damage at the very top of our design domain would result in arbitrary large compliance values (due to unsupported loads), we restrict the possible damage positions to exclude the top \texttt{nonD} rows of design elements. Lastly, the parameter \texttt{dmg\_fac} controls how strongly the structure is damaged in the design region, where \texttt{dmg\_fac=1} corresponds to completely removing all material, while \texttt{dmg\_fac=0} will not remove any material at all.
\subsection{Setting up the stochastic model}\label{sec:StoModSetup}
Initialization of variables associated with the stochastic nearest neighbor model (\texttt{PRE. 5}) is implemented in lines 43-54:
\begin{lstlisting}[frame=single,style=Matlab-editor,basicstyle={\fontsize{7}{8}\ttfamily}]
rng('default');
com0 = 100;
[y1, y2] = meshgrid( 1:(nelx-L+1) , 
    1:(nely-L+1-nonD)  ); 
y = [y1(:),y2(:)]';
n_disc = size(y,2);
X = [randi(nelx-L+1,1,maxit) ; 
    randi(nely-L+1-nonD,1,maxit)];
maxsmpl = 2000;
[x_birth,x_ind,leavers] = deal(1:maxsmpl,
    1:maxsmpl,1);
y_weight = volfrac*sqrt(nEl);
y_diff = pdist2(y',X');
y_diff = y_diff/max(max(y_diff,1e-10),[],'all');
[Gra,DesH,ComH] = deal( zeros(nEl,maxsmpl), zeros(nEl,maxsmpl), zeros(1,maxsmpl) );
\end{lstlisting}
At first, the random number generator is fixed, to ensure reproducibility. The parameter \texttt{com0} in line 44 serves as an initial guess to the structure's expected compliance and will later be used to rescale the objective function. As quadrature points \texttt{y}, on which the nearest neighbor model will be evaluated (see~\eqref{eq:Quadrature}), we choose all possible damage cases, which are constructed in lines 45\&46. The sample sequence \texttt{X} is drawn in line 48 and \texttt{maxsmpl} in line 49 limits the amount of samples stored during the optimization process, see~\Cref{sec:LimitedWeights}. Here, we also initialize the auxiliary variables \texttt{x\_birth}, \texttt{x\_ind} and \texttt{leavers}, which later contribute to the sample management in \texttt{RL. 4)}.

For a more efficient integration weight calculation scheme, it is beneficial to decompose the product norm $\Vert\cdot\Vert_\ast$ in~\eqref{eq:NNIndex} into a sum of norms over each individual space
\begin{equation}
    \Vert (\hat{\mx},\xi,\psi)\Vert_\ast = \Vert\hat{\mx}\Vert+c_1\Vert\xi\Vert+c_2\Vert\psi\Vert.\label{eq:NormDecomposition}
\end{equation}
Abusing the independence of $\hat{\mx}$ and $\xi$ in this norm, we can construct and normalize the full random space distance matrix
\begin{align*}
    &\texttt{y\_diff}\in[0,1]^{\texttt{n\_disc}\times\texttt{maxit}}, \\
    &\big(\texttt{y\_diff}\big)_{i,j}:= \frac{\Vert \texttt{y(:,i)}-\texttt{X(:,j)}\Vert_2}{\max_{i,j}\Vert \texttt{y(:,i)}-\texttt{X(:,j)}\Vert_2}
\end{align*}
before the optimization procedure (lines 52\&53). As a result, only the design differences ${\Vert\hat{\mx}_k-\hat{\mx}_{1,\ldots,k}\Vert}$ need to be calculated in each iteration. The parameters $c_1,c_2\in\R_{\ge0}$ in~\eqref{eq:NormDecomposition} allow for fine-tuning the importance of individual components. In our setting, we want to consider both spaces to be (roughly) of equal importance. However, due to the normalization of \texttt{y\_diff}, the maximum possible distance between a pair ${(\xi_i,\xi_j)}$ is normalized to be 1. On the other hand, the Euclidean distance between the empty design region and a full design region scales like $\sqrt{\texttt{nelx}\cdot\texttt{nely}}$, practically rendering the difference in $\xi$ irrelevant. Thus, in line 51, we adjust the parameter \texttt{y\_weight}, corresponding to $c_1$ in~\eqref{eq:NormDecomposition}, accordingly. All in all, we have
\begin{align*}
    \big\Vert (\hat{\mx}_i,&\xi_i)-(\hat{\mx}_j,\xi_j)\big\Vert_\ast \\
    & = \Vert\hat{\mx}_i-\hat{\mx}_j\Vert_2 + \texttt{y\_weight}\cdot\texttt{y\_diff(i,j)}.
\end{align*}
Also, it should be noted that choosing the Euclidean norm instead of, e.g., $\Vert\cdot\Vert_1$, is arbitrary.

Lastly, \texttt{Gra}, \texttt{DesH} and \texttt{ComH} initialize the vectors in which gradient samples and associated designs as well as samples of the objective function are stored, respectively.
\subsection{Sampling and integration weight calculation}
To set up the state equation, we first damage the structure according to the current random sample (lines 89-91):
\begin{lstlisting}[frame=single,style=Matlab-editor,basicstyle={\fontsize{7}{8}\ttfamily}]
D = zeros(nely,nelx);
D( nely+1-(X(2,loop):X(2,loop)+L-1), 
    X(1,loop):X(1,loop)+L-1 ) = 1;
x_dmg = max(0, min(1, xPhys-dmg_fac*D(:)));
\end{lstlisting}

After solving the state equation, the current gradient sample \texttt{dc}, compliance sample \texttt{F'*U} and physical design \texttt{xPhys} are stored (line 103).

To calculate the current integration weights, for each integration point \texttt{y(:,i)}, we have to find the index of the closest sample point in memory, see~\eqref{eq:NNIndex}. By our specific choice of norm and precalculation of the distance matrix \texttt{y\_diff} (\Cref{sec:StoModSetup}), this can be done by finding the position of the minimum entry in the vector
\begin{align*}
    \Vert \hat{\mx}_\text{loop} - \texttt{DesH}&\texttt{(:,1:loop)}\Vert_2 \\
    &+ \texttt{y\_weight}\cdot\texttt{y\_diff(i,1:loop)}.
\end{align*}
Afterwards, for each sample, we add up the associated weights $w_t$ in the quadrature rule used for the construction of \texttt{y}, see~\eqref{eq:weights_pseudo}. Thus, a straight-foward implementation of this process could be
\begin{lstlisting}[frame=single,style=Matlab-editor,basicstyle={\fontsize{7}{8}\ttfamily}]
for t = 1:T
    [~,i_t] = min(
        vecnorm(xPhys-DesH(:,1:ulim),2,1) 
         + y_weight*y_diff(t,x_ind(1:ulim)));
    weights(i_t) = weights(i_t) + w_t;
end
\end{lstlisting}
In the case of ${w_t=\tfrac{1}{T}}$ for all ${t=1,\ldots,T}$, this can be vectorized as implemented in lines 104\&105:
\begin{lstlisting}[frame=single,style=Matlab-editor,basicstyle={\fontsize{7}{8}\ttfamily}]
[~,csw] = min( vecnorm(xPhys-DesH(:,1:ulim),2,1)
    + y_weight*y_diff(:,x_ind(1:ulim)), [], 2);
weights = sum(csw==1:ulim)/n_disc;
\end{lstlisting}
Adjusting this expression to more general cases for $w_t$ is done in~\Cref{sec:load}.
If ${\texttt{maxsmpl}<\texttt{maxit}}$, i.e., we will be required to remove old samples from memory at some point (\Cref{sec:LimitedWeights}), lines 106-113 
\begin{lstlisting}[frame=single,style=Matlab-editor,basicstyle={\fontsize{7}{8}\ttfamily}]
ind_can = find(weights-min(weights)<1e-8);
[~,iind] = min(x_birth(ind_can));
leavers = ind_can(iind);
[x_ind(leavers), x_birth(leavers)] 
    = deal(loop+1,loop+1);
\end{lstlisting}
first determine samples with minimal associated integration weights. If there is more than one sample with minimum integration weight, we pick the oldest sample to be removed. For this, we use the auxiliary variable \texttt{x\_birth} to store the iteration in which a sample was drawn. Since this procedure leads to an unordered set of samples, we keep track of the ordering in \texttt{x\_ind}.
\subsection{Nearest neighbor approximations}
With the integration weights \texttt{weights} calculated as proposed in the previous section, integrating the nearest neighbor model~\eqref{eq:weights_pseudo} now corresponds to a weighted summation of our samples, as performed in line 115
\begin{lstlisting}[frame=single,style=Matlab-editor,basicstyle={\fontsize{7}{8}\ttfamily}]
Compl = sum(weights.*ComH(:,1:ulim),2)
\end{lstlisting}
To ensure a better scaling of the objective function, 
\begin{equation*}
    \frac{1}{P}\mathbb{E}_{\xi}\left[ c(\hat{\mx},\xi)^P \right]
\end{equation*}
is replaced by
\begin{equation*}
    \frac{1}{P}\mathbb{E}_{\xi}\left[ \left(\frac{c(\hat{\mx},\xi)}{\texttt{com0}}\right)^P \right],
\end{equation*}
where \texttt{com0} is an approximation to ${\mathbb{E}_{\xi}\left[ c(\hat{\mx},\xi,\psi) \right]}$, which gets updated every 25 iterations. The corresponding gradient \texttt{dc} is calculated in line 119.
\section{Numerical examples}\label{sec:NumExamples}
\subsection{Reference problem}\label{sec:ClampNonsym}
To solve the reference problem introduced in~\Cref{sec:ExmplSetup}, we can call
\begin{lstlisting}[frame=single,style=Matlab-editor]
topS140(180,45,0.4,3,3.2,2,'N',
    0.5,2,1e-2,1,1500);
\end{lstlisting}
to obtain the design shown in~\Cref{fig:Design_reference_problem}. The approximated objective function values $J_k$ can be found in~\Cref{fig:Objective_reference_problem}. Therein, we also included snapshots of the true objective ${\mathbb{E}_\xi\big[c(\hat{\mx}_k,\xi)\big]}$, evaluated by solving the state equation for all $3,381$ possible damage cases. Although this means that a single exact evaluation requires more than twice the amount of system solves we allowed for the full optimization process, the errors of the stochastic approximations are rather small. To be precise, we have ${J_\text{final}=9.49}$ and ${\mathbb{E}_\xi\big[c(\hat{\mx}_\text{final},\xi)]=9.34}$. 

In~\Cref{fig:Objective_reference_problem}, we also see that true objective function value is much more well-behaved than the stochastic approximation. This is somewhat expected, as the amount of available samples is rather small, meaning that the nearest neighbor model is not fully converged by the end of the optimization process. On the bright side, this means that convergence to the optimal design is typically observed faster than indicated by $J_k$. That being sad, it also means that using $J_k$ as a measure of convergence will result in a larger number of steps than necessary. Thus, for simplicity, we use the total number of iterations as a hard stopping criterion in our numerical experiments.
\begin{figure}
    \centering
    \includegraphics[width=\linewidth]{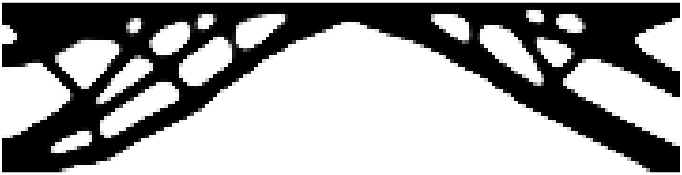}
    \caption{Final design for the reference problem (obj.: 9.34).}
    \label{fig:Design_reference_problem}
\end{figure}
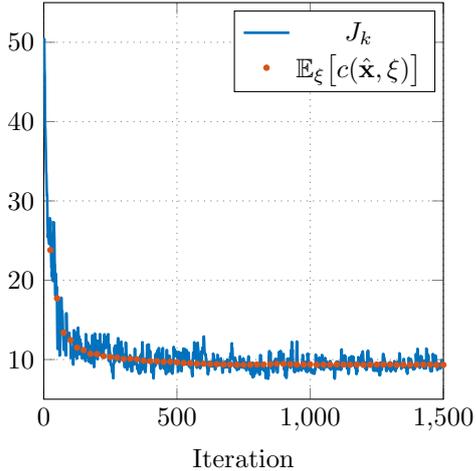
\begin{figure}
    \centering
    \input{Graphics/Clamp/Reference_obj}
    \caption{True objective function values (orange dots) and stochastic approximations of the nearest neighbor model (blue) for the reference problem over the course of iterations.}
    \label{fig:Objective_reference_problem}
\end{figure}
\subsection{Deterministic problem}
\begin{figure}
    \centering
    \includegraphics[width=\linewidth]{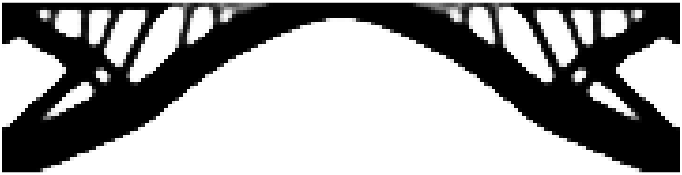}
    \caption{Final design for the damage-free deterministic reference problem (obj.: 19.38).}
    \label{fig:Design_deterministic_problem}
\end{figure}

With minor modifications, we can also use the code to solve the damage-free deterministic problem. For that, we simply set \texttt{dmg\_fac=0} in line 41 and \texttt{maxsmpl=1} in line 49. Thus, no damage is applied to the structure and the nearest neighbor approximation is essentially bypassed. The resulting final design when using the same parameters as in~\Cref{sec:ClampNonsym} is depicted in~\Cref{fig:Design_deterministic_problem}. As expected, neglecting the damage in the optimization leads to a suboptimal design, with the expected compliance of $19.38$ being roughly twice as large as for the optimal design found in~\Cref{sec:ClampNonsym}. Note that, for simplicity, we did not fix a solid layer of material at the top of the structure. Thus, it not surprising that the final design includes elements with grey material ($0<\hat{x}_e<1)$.
\subsection{Impact of random sequence}
Since the nearest neighbor model depends on the random sample sequence, we expect the routine to produce different optimal designs for different sample sequences. While the results presented in~\Cref{sec:ClampNonsym} are representative for the overall performance of the method, it is still important to analyze the differences between individual runs. Therefore, we remove line 43 from the code and perform 100 optimization runs with the same parameters as given in~\Cref{sec:ClampNonsym}. To determine the deviations between different optimization instances, we track the true objective function values ${\mathbb{E}_\xi\big[c(\hat{\mx}_k,\xi)\big]}$. Due to the high computational cost, we only calculate this value every 10 iterations.

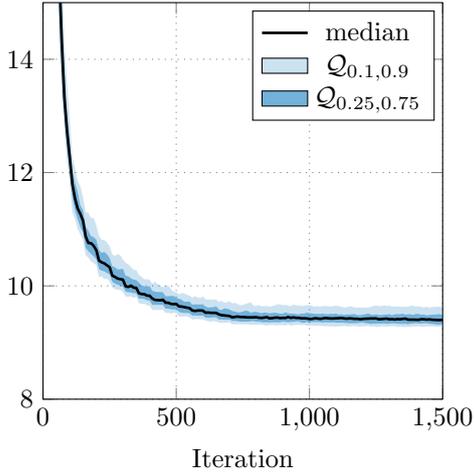
\begin{figure}
    \centering
    \input{Graphics/Clamp/Clamp_quants}
    \caption{Median of ${\mathbb{E}_\xi\big[c(\hat{\mx}_k,\xi)\big]}$ over all 100 individual optimization runs (solid black line), evaluated every 10th iteration. The shaded areas indicate $\mathcal{Q}_{0.1,0.9}(k)$ (light) and $\mathcal{Q}_{0.25,0.75}(k)$ (dark).}
    \label{fig:Quantiles_clamp}
\end{figure}

Given ${q_1\in[0,1]}$, we define the empirically observed quantiles
\begin{equation*}
    \mathcal{Q}_{q_1}(k) := \inf\big\{ t\in\R: \mathbb{P}\big[\mathbb{E}_\xi\big[c(\hat{\mx}_k,\xi)\big]\le t\big]\ge q_1 \big\},
\end{equation*}
where the probability is taken over all 100 optimization runs. For example, $\mathcal{Q}_{0.37}(251)$ is the smallest real number that is larger than the objective function value in iteration 251 for at least 37\% of all individual runs. Thus, $\mathcal{Q}_{q_1}$ allows us to ``sort out" extremely poor and extremely good individual runs, in order to focus on the ``average" performance of the method. Now, for ${0\le q_1<q_2\le 1}$, we set
\begin{equation*}
    \mathcal{Q}_{q_1,q_2}(k) := \big[ \mathcal{Q}_{q_1}(k), \mathcal{Q}_{q_2}(k)\big].
\end{equation*}
For example, $\mathcal{Q}_{0.1,0.8}(k)$ contains all values of ${\mathbb{E}_\xi\big[c(\hat{\mx}_k,\xi)\big]}$ after excluding the smallest $10\%$ and largest $20\%$ of observed values. Plots for $\mathcal{Q}_{0.1,0.9}(k)$, $\mathcal{Q}_{0.25,0.75}(k)$ as well as the median over all runs can be found in~\Cref{fig:Quantiles_clamp}. Therein, we see that the random sequence has only a minor impact on the objective function evolution.
\subsection{Incorporating symmetry}\label{sec:Symmetry}
Due to the probabilistic nature of the algorithm, the underlying symmetry (at $\texttt{nelx}/2$) in the model is lost. While simply symmetrizing the gradient through the additional line of code (to be inserted between lines 99\&100 in the original code)
\begin{lstlisting}[frame=single,style=Matlab-editor,basicstyle={\fontsize{7}{8}\ttfamily}]
dc = (dc + dc(:,nelx:-1:1));
\end{lstlisting}
is enough to reinstate this symmetry, we expect a better performance if the symmetry is additionally directly incorporated in the stochastic nearest neighbor model.

\begin{figure}
    \centering
    \includegraphics[width=\linewidth]{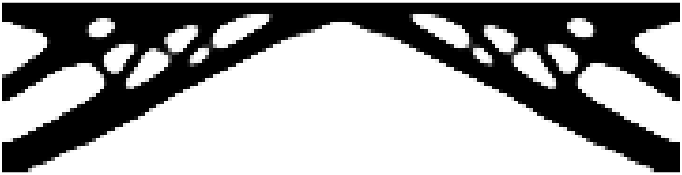}
    \includegraphics[width=\linewidth]{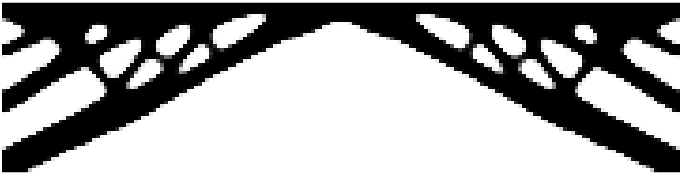}
    \caption{Final designs for enforced symmetry. Top: unmodified nearest neighbor model (obj.: 9.39). Bottom: modified nearest neighbor model (obj.: 9.41).}
    \label{fig:Design_clamp_sym}
\end{figure}

Note that, for a symmetric design, damaging the structure in the right half will yield the same compliance (and gradient) as the corresponding damage case in the left half. Thus, we may define the equivalence relation $\sim_R$ on ${\Xi\times\Xi}$ by
\begin{equation*}
    \xi\sim_R\phi\,:\,\Longleftrightarrow\vert\xi_1-2\ell\vert = \vert\phi_1-2\ell\vert\,\land\, \xi_2=\phi_2
\end{equation*}
and consider the quotient space ${\Xi_R:=\Xi\big/\sim_R}$ instead of $\Xi$. In other words: Damage cases in the right half will be reflected on the corresponding samples in the left half before calculating the norm differences. This can be implemented by replacing the calculation of \texttt{y\_diff} by the following block of code:
\begin{lstlisting}[frame=single,style=Matlab-editor,basicstyle={\fontsize{7}{8}\ttfamily}]
y_diff_l = pdist2(y',X');
y_diff_r = pdist2(y',[nelx-L+2-X(1,:);X(2,:)]');
y_diff = min(y_diff_l,y_diff_r);
y_diff = y_diff/max(max(y_diff,
    1e-10),[],'all');
\end{lstlisting}
The resulting final designs using the unmodified as well as the modified stochastic models are shown in~\Cref{fig:Design_clamp_sym}. Note that the symmetrization of \texttt{dc} (see above) is required in both cases.

\begin{figure}
    \centering
    \input{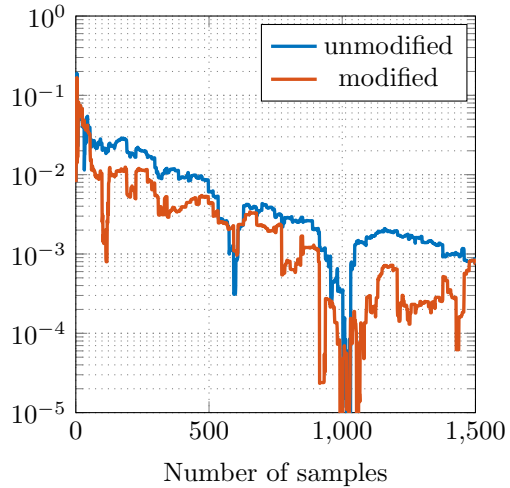}
    \caption{Relative approximation error ${\tfrac{\vert J_n-c^\ast\vert}{c^\ast}}$ for pure nearest neighbor integration using the unmodified norm (blue) and modified norm (orange).}
    \label{fig:Error_clamp_sym}
\end{figure}

Although, in this particular example, both objective function values are practically identical ($9.39$ for unmodified and $9.41$ for modified), adjusting the nearest neighbor norm according to underlying symmetry indeed leads to a better approximation of the expected value. To see this, we use the method as a plain integration tool by fixing the optimal design, \texttt{beta=16} and \texttt{move=0}. Now, we obtain iterative approximations to the true objective function value $c^\ast$, based on pure nearest neighbor integration for both norms. The resulting relative errors ${\tfrac{\vert J_n-c^\ast\vert}{c^\ast}}$ are shown in~\Cref{fig:Error_clamp_sym}.
\section{Variations}\label{sec:Generalizations}
Implementing adjusted boundary conditions, multiple load cases, etc. is straightforward and can be done in the same fashion as for \texttt{top99neo}. Similarly, an extension to three-dimensional problems can be derived from \texttt{top3D125}~\cite{99neoCode}.

Thus, we present three exemplary variations of the method. First, we analyze strategies to adjust the sampling sequence or to reduce the number of possible damage cases. Afterwards, we consider probabilistic force vectors in a data-driven context. The corresponding \texttt{topS140\_load} code is given in~\Cref{sec:AppendixCode_load}. Lastly, we extend the scheme to stochastic mini-batches for dynamic load cases.
\subsection{Damage case reduction}
\begin{figure}
    \centering
    \input{Graphics/Beam/beam_setup}
    \caption{Design domain with support at left boundary. A downward force (red) is applied at the midpoint of the right boundary. The blue square indicates a possible damage case. The location of the damaged region is assumed to be uniformly random distributed to the left of the dashed green line.}
    \label{fig:beam_setup}
\end{figure}
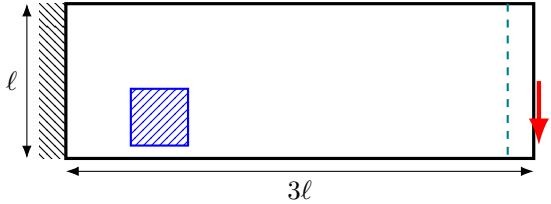

We investigate the cantilever beam optimization considered in~\cite{Jansen2013}, where a design region of width $3\ell$ and height $\ell$ is fixed at the left boundary and a downward force is applied at the middle of the right boundary. Again, a randomly placed square-shaped hole is used to model local material failure. Since the optimizer has no chance of compensating for material failure around the nodes at which the force is applied, we do not consider damage cases too close to the right boundary. An illustration of the setup can be found in~\Cref{fig:beam_setup}. Although the authors of~\cite{Jansen2013} were concerned with finding a fail-safe design, i.e., minimizing the compliance value of the worst possible damage case, we instead analyze the problem with respect to the expected compliance of the structure.

\begin{figure}
    \centering
    \includegraphics[width=\linewidth]{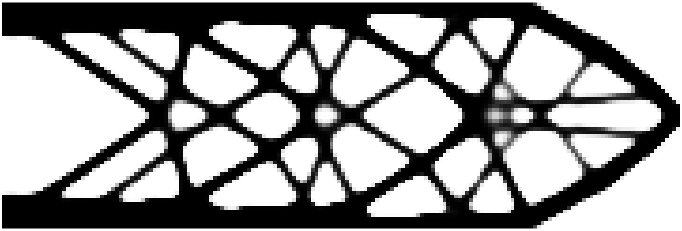}
    \includegraphics[width=\linewidth]{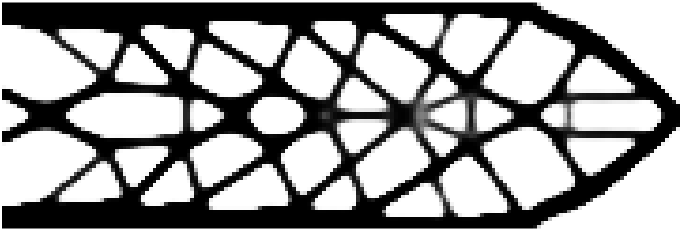}
    \includegraphics[width=\linewidth]{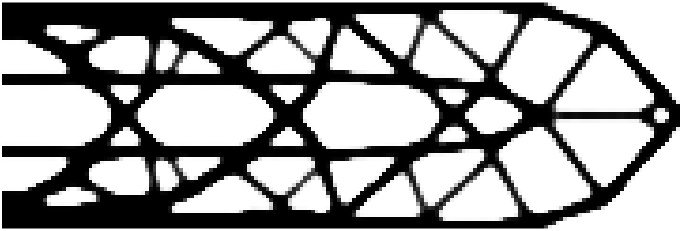}
    \caption{Final designs for beam optimization (see~\Cref{fig:beam_setup}). Top: Standard approach (obj.: 496.84). Mid: Oversampling boundary (obj.: 383.39). Bottom: Damage case reduction according to~\Cref{fig:Beam_damagegrid} (obj.: 490.58).}
    \label{fig:Beam_design}
\end{figure}

After adjusting the setup in the code, we also enforce symmetry with respect to the horizontal line at $y=\tfrac{\ell}{2}$ and employ the modified norm as explained in~\Cref{sec:Symmetry}. Fixing $\texttt{move=2.5e-3}$ and choosing all other parameters identical as in~\Cref{sec:ClampNonsym}, we obtain the final design shown in~\Cref{fig:Beam_design} (top).

The corresponding true objective function value ($496.84$) is poorly approximated by the nearest neighbor model ($1,262.74$), indicating that the chosen settings should be reconsidered. When investigating the final design, it is easy to see why that is the case: Locating the damage region in such a way that one of the main supports is cut results in very large compliance values of the damaged structure. Thus, in comparison to the examples above, the variation in compliance with respect to the random parameter is much larger. Specifically, ``edge cases", in which the damage parameter is located at the boundaries of $\Xi$, require more emphasis during the sampling process and are simply ``not seen enough" (only $1,500$ random samples from $5,811$ total possibilities).

To deal with this problem, we propose two different approaches. First, as mentioned above, we know beforehand that individual edge cases pose problems. Thus, we may purposefully oversample these critical regions. For example, after generating a (uniform) sampling sequence \texttt{X}, we include the following two lines of code:
\begin{minipage}{\linewidth}
\begin{lstlisting}[frame=single,style=Matlab-editor,basicstyle={\fontsize{7}{8}\ttfamily}]
X(1,1:15:end) = 1;
X(2,1:10:end) = 1;
\end{lstlisting}
\end{minipage}
Note that this does not impact the construction of our nearest neighbor model, but forces more samples to lie on the boundary. The corresponding final design for the modified approach can be found in~\Cref{fig:Beam_design} (mid). This time, the true final objective function value ($383.39$) is greatly reduced and much better approximated by the nearest neighbor model ($460.87$).

\begin{figure}
    \centering
    \includegraphics[width=\linewidth]{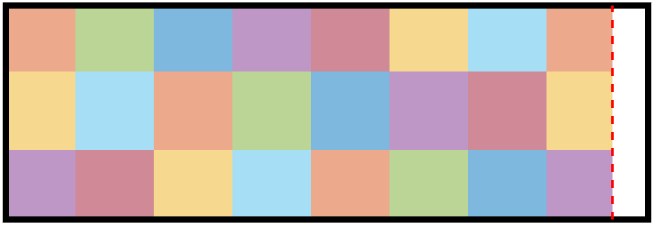}
    \includegraphics[width=\linewidth]{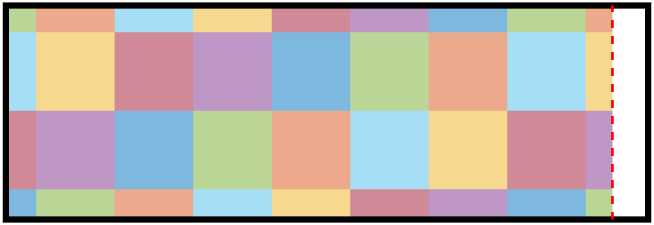}
    \caption{Visualization of all 60 considered damage cases in the reduced model.}
    \label{fig:Beam_damagegrid}
\end{figure}

Although these adjustments may be reasonable for this specific setup, it is unclear whether a strategic adjustment of the sampling sequence can be made with only a priori information in general. For this reason,~\cite{Zhou2016} proposed a generalized approach for replacing the full damage model with a coarse discretization grid. Following this methodology, we construct a grid of $60$ possible damage cases to cover the design region. All possible realizations are shown in~\Cref{fig:Beam_damagegrid}.

Using this model reduction with the same parameters, we obtain the third final design in~\Cref{fig:Beam_design} (bottom). Due to the drastic reduction in possible damage cases, the reduced final expected compliance over the 60 damage cases ($380.03$) is approximated almost exactly by the method ($379.96$). However, note that we no longer solve the original problem, which becomes apparent when looking at the true expected compliance over all possible damage cases ($490.58$).

A comparison of the approaches can be found in~\Cref{tab:Beam_comparison}.
\begin{table}
    \centering
    \begin{tabular}{|l||r|r|r|}
        \hline
        \textbf{method} & $J_{\text{final}}$ & \textbf{obj. (full)} & \textbf{obj. (red.)}  \\
        \hline
         standard & 1,262.74 & 496.84 & -- \\
         \hline
         oversample & 460.87 & 383.39 & -- \\
         \hline
         reduction & 379.96 & 490.58 & 380.03\\
         \hline
    \end{tabular}
    \caption{Final objective approximations $J_\text{final}$ (left) and objective function values ${\mathbb{E}_\xi\big[c(\hat{\mx}_\text{final},\xi)\big]}$ (mid) for all three approaches. Note that, when using the damage case reduction, the method no longer approximates the true objective, but the reduced expected compliance ${\tfrac{1}{60}\sum_{i=1}^{60}c(\hat{\mx}_\text{final},\xi_i)}$, given in the right column.}
    \label{tab:Beam_comparison}
\end{table}
\subsection{Uncertainty in force vector}\label{sec:load}
\begin{figure}
    \centering
    \input{Graphics/Load/load_setup}
    \caption{Design domain with supports at lower left and right boundary. The position of the downward force acting on the structure is randomly chosen on the top of the domain (red).}
    \label{fig:load_setup}
\end{figure}
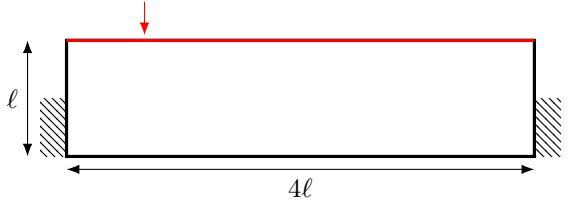

So far, we have only considered the case that uncertainty enters the system via the parameter $\xi$ in the stiffness matrix~\eqref{eq:StateEquation}. Thus, for the next example, we want to shift our attention to the probabilistic parameter $\psi$ in our force vector. For that, we consider a setup similar to~\Cref{sec:ExmplSetup}. This time, we do not model any local material failure, but instead assume that the downward force is applied at a random node at the top of our design domain. An illustration of the setup is given in~\Cref{fig:load_setup}. Furthermore, we want to use this example to showcase how to adapt the method to data-driven optimization, where the probability distribution is unknown and only a collection of sample points is available.

To simulate this situation, we construct a fixed set of samples by the following lines of code:
\begin{minipage}{\linewidth}
\begin{lstlisting}[frame=single,style=Matlab-editor,basicstyle={\fontsize{7}{8}\ttfamily}]
rng('default');
p1 = makedist('Normal','mu',0.25,'sigma',0.1);
p1 = truncate(p1,0,1);
p2 = makedist('Normal','mu',0.6,'sigma',0.2);
p2 = truncate(p2,0,1);
r1 = random(p1,1,3e5);
r2 = random(p2,1,1e5);
r = [r1,r2];
r = r(randperm(4e5));
r = round(nelx*r)+1;
save('Data.mat','r')
\end{lstlisting}
\end{minipage}
Throughout this section, we consider \texttt{r} to be an arbitrary given dataset of load cases, where $\texttt{r[i]}$ indicates at which top node the force is applied. The full dataset for ${\texttt{nelx} = 360}$ is shown in~\Cref{fig:Load_histogram}.
\begin{figure}
    \centering
    \input{Graphics/Load/Distribution_histogram}
    \caption{Probability of all 361 possible load cases in the dataset for ${\texttt{nelx} = 360}$.}
    \label{fig:Load_histogram}
\end{figure}
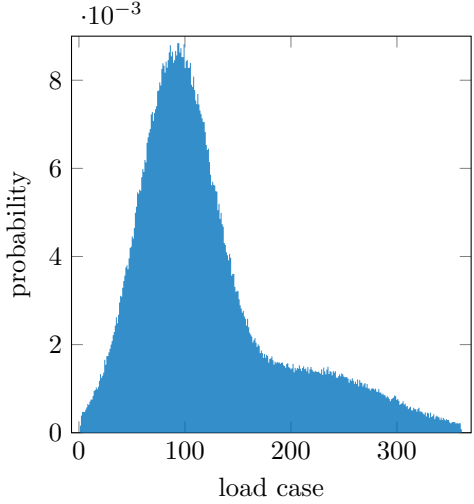
\begin{remark}\label{rem:Deterministic}
    Note that the implementation uses a decomposition technique to solve the state equation~\eqref{eq:StateEquation}. Combined with the fact that the uncertainty impacts only the force vector, we can reuse the decomposition to solve for multiple load cases simultaneously. Therefore, a deterministic approach, in which we solve for each possible load case in every iteration, is still numerically feasible.
\end{remark}
When setting up the nearest neighbor model, we adjust the integration points to consist of all individual possibilities observed in the full dataset
\begin{lstlisting}[frame=single,style=Matlab-editor,basicstyle={\fontsize{7}{8}\ttfamily}]
y_all = load('Data.mat','r'); y_all = y_all.r; 
y = unique(y_all); n_disc = size(y,2);
\end{lstlisting}
Furthermore, we need the probabilities $p_i$ of all possible scenarios to calculate the expected compliance
\begin{equation*}
    \mathbb{E}_\psi\big[c(\hat{\mx},\psi)\big] = \sum_{i}p_ic(\hat{\mx},\psi_i).
\end{equation*}
Since these are assumed to be unknown, they are reconstructed empirically from the dataset in line 43:
\begin{lstlisting}[frame=single,style=Matlab-editor,basicstyle={\fontsize{7}{8}\ttfamily}]
dist_w = sum(y_all' == y)/size(y_all,2);
\end{lstlisting}
Lastly, we need to change the integration weight calculation in line 105 accordingly (compare~\eqref{eq:weights_pseudo}):
\begin{lstlisting}[frame=single,style=Matlab-editor,basicstyle={\fontsize{7}{8}\ttfamily}]
weights = sum((csw==1:ulim).*dist_w');
\end{lstlisting}
A key advantage of the nearest neighbor model over standard stochastic procedures from literature is the fact that the sampling sequence is decoupled from the probability distribution of random parameters. In other words, we can choose between different approaches on how to generate the sampling sequence:
\begin{enumerate}
    \item Generate the sampling sequence by randomly drawing from the dataset. This results in a sampling sequence following the probability distribution of the random parameters and is the usual approach for standard stochastic approaches.
    \item Generate the sampling sequence in any other fashion, e.g., uniformly distributed over all possible cases. Although this seems counterintuitive at first, we will later see why this may significantly improve the performance in some cases.
\end{enumerate}
In the \texttt{topS140\_load} code, the construction of our random sequence is controlled by the parameter \texttt{type} and performed in lines 44-49
\begin{lstlisting}[frame=single,style=Matlab-editor,basicstyle={\fontsize{7}{8}\ttfamily}]
switch type
    case 'distribution' % sample according to distribution
        X = y_all(randi(numel(y_all),1,maxit));
    case 'uniform'  % sample uniformly
        X = y(randi(n_disc,1,maxit));
end
\end{lstlisting}
For our analysis, we run
\begin{lstlisting}[frame=single,style=Matlab-editor,basicstyle={\fontsize{7}{8}\ttfamily}]
topS140_load(360,90,0.4,3,6.4,2,'N',0.5,2, 
    1e-2,1,1500,type)
\end{lstlisting}
once for each \texttt{type} and compare the results to the solution found by a deterministic scheme using the exact gradient by simulating each possible load case in every iteration, see~\Cref{rem:Deterministic}.

\begin{figure}
    \centering
    \includegraphics[width=\linewidth]{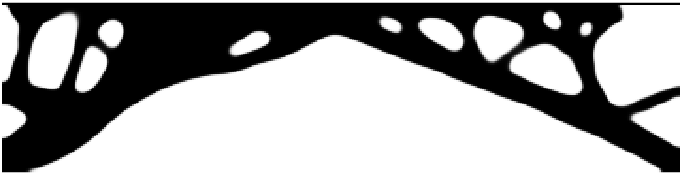}
    \includegraphics[width=\linewidth]{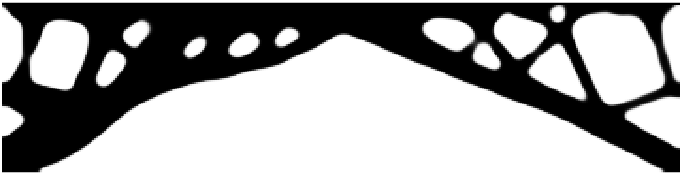}
    \includegraphics[width=\linewidth]{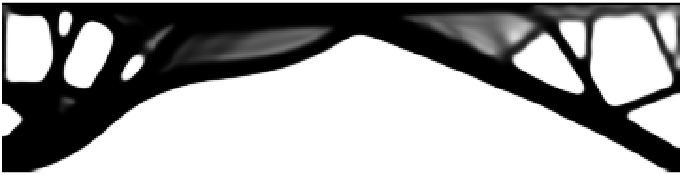}
    \includegraphics[width=\linewidth]{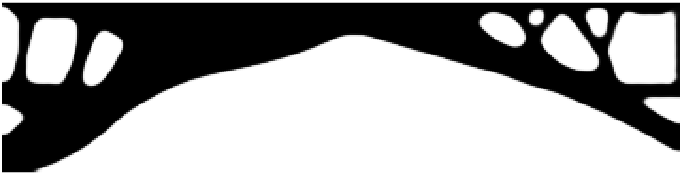}
    \caption{Final designs for probabilistic force position (see~\Cref{fig:load_setup}). Top: \texttt{type='distribution'} (obj.: 203.37). Second: \texttt{type='uniform'} (obj.: 23.63). Third: Exact deterministic approach with $p=3$ (obj.: 23.95). Bottom: Exact deterministic approach with $p=6$ (obj.: 25.14).}
    \label{fig:Load_design}
\end{figure}

The final physical designs of all methods are shown in~\Cref{fig:Load_design}. Interestingly, the deterministic approach struggles to get rid of intermediate material. Thus, we also run it with an increased SIMP parameter of ${p=6}$.

Of all approaches, sampling according to the probability distribution actually yields by far the largest final objective function value (203.37). When looking at the corresponding design in~\Cref{fig:Load_design}, it becomes clear why this is the case: A large portion at the right upper boundary is not sufficiently supported by material (only passive solid elements are present here). This is a direct consequence of the small probability associated to these force positions, see~\Cref{fig:Load_histogram}. Thus, during the optimization process, this area is sampled very sparsely, leading to inaccurate approximations of the compliance and gradient. Simply speaking: By drawing samples according to the distribution we create a blind spot in the optimizer.

In contrast, choosing \texttt{type='uniform'} results in a much better performance (final objective function value: 23.63), as the quality of the constructed nearest neighbor model no longer deteriorates in region with low probability. It should be noted that sampling uniformly over all observed cases also has some downsides. For example, if there exist regions in which the probability density is so small that $p_i\nabla_\mx c(\hat{\mx},\psi_i)$ is close to zero, this procedure might produce ``unnecessary" information.

Surprisingly, the deterministic scheme yields a slightly larger final objective function value (23.95). The reason for this is the large amount of intermediate material present in the final design. Increasing the SIMP parameter to $p=6$ solves this issue (see~\Cref{fig:Load_design}) and results in a final objective function value of 24.61. An overview of all final expected compliance values, evaluated for $p=3$ and $p=6$, can be found in~\Cref{tab:load_comparison}.
\begin{table}
    \centering
    \begin{tabular}{|l||r|r|}
        \hline
        \textbf{method} & SIMP $p=3$ & SIMP $p=6$\\
        \hline
         \texttt{distribution} ($p=3$) & 203.37 & 209.84 \\
         \hline
         \texttt{uniform} ($p=3$) & 23.63 & 24.50 \\
         \hline
         deterministic ($p=3)$ & 23.95 & 27.77\\
         \hline
         deterministic ($p=6$) & 24.61 & 25.14\\
         \hline
    \end{tabular}
    \caption{Expected compliance values of final designs obtained by all different approaches. Each design has been evaluated for a SIMP parameter of $p=3$ (left column) and $p=6$ (right column).}
    \label{tab:load_comparison}
\end{table}
\subsection{Mini-batching}
\begin{figure}
    \centering
    \input{Graphics/Dynamic/dyna_setup}
    \caption{Design domain with support at the bottom right corner and symmetry boundary conditions at left boundary. The applied force (red) is assumed to be time-harmonic with frequency ${\omega\in[\omega_\text{min},\omega_\text{max}]}$.}
    \label{fig:dynamic_setup}
\end{figure}
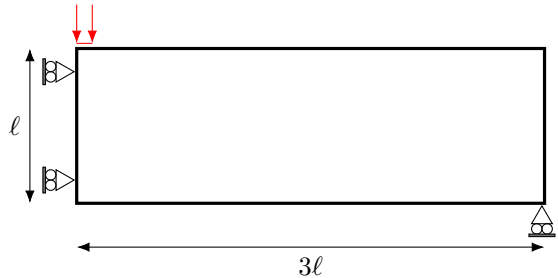

Consider the famous MBB beam setup, i.e., a rectangular design domain of height $\ell$ and width $6\ell$, which is fixed (in $y$-direction) at the bottom corners and loaded at a small area in the middle of the top boundary. A symmetry-adapted setup of this problem is shown in~\Cref{fig:dynamic_setup}. For this example, we no longer model local material failure. Instead, we consider the applied force to oscillate harmonically with uncertain frequency ${\omega\in[\omega_\text{min},\omega_\text{max}]}$. A common objective in this setup is to maximize the structure's fundamental eigenvalue~\cite{Andreassen2017}, which indirectly minimizes the expected dynamic compliance ${\mathbb{E}_\omega\big[c(\hat{\mx},\omega)\big]}$ of the structure. Here, we instead aim to minimize ${\mathbb{E}_\omega\big[c(\hat{\mx},\omega)\big]}$ directly.

In addition to the changes required to solve the dynamic state equation, we need to adjust the setup of the stochastic model. For example, choosing $\omega_\text{min}=0$ (static load), $\omega_\text{max}=0.2$ and a uniform distribution for $\omega$, we set
\begin{lstlisting}[frame=single,style=Matlab-editor,basicstyle={\fontsize{7}{8}\ttfamily}]
wInt = [0, 0.2];
n_disc = 3000;
y = linspace(wInt(1),wInt(2),n_disc);
\end{lstlisting}
Additionally, we want to increase the batch size of our sampling strategy, i.e., the gradient is calculated not only for a single random frequency, but for a (small) number of random frequencies. As an example, we assume that a processor with four cores may be used to solve four state equations in parallel. For our sample management, this means that we also have to remove four old samples from memory (instead of one), once \texttt{maxsmpl} samples have been evaluated. The corresponding adjustment of the sampling sequence \texttt{X} as well as the auxiliary variables \texttt{x\_birth} and \texttt{leavers} looks as follows:
\begin{lstlisting}[frame=single,style=Matlab-editor,basicstyle={\fontsize{7}{8}\ttfamily}]
bsz = 4;
x_birth = sort(repmat(1:maxsmpl/bsz,1,bsz));
leavers = 1:bsz;
X = wInt(1)+rand(1,maxit*bsz)*(wInt(2)-wInt(1));
\end{lstlisting}
In each iteration, determining the old samples to be removed from storage is realized by replacing lines 106-113 with
\begin{minipage}{\linewidth}
\begin{lstlisting}[frame=single,style=Matlab-editor,basicstyle={\fontsize{7}{8}\ttfamily}]
if (loop+1)*bsz <= maxsmpl
    leavers = loop*bsz+1 : (loop+1)*bsz;
else
    minweight = sort(weights);
    ind_can = find(weights-minweight(bsz)<1e-8);
    can_birth = x_birth(ind_can);
    [~,iind] = sort(can_birth);
    leavers = ind_can(iind(1:bsz));
    x_ind(leavers) = loop*bsz+1 : (loop+1)*bsz;
    x_birth(leavers) = loop+1;
end
\end{lstlisting}
\end{minipage}
Due to dynamic resonances and other effects, the compliance gradient is no longer guaranteed to have non-positive entries. Thus, the OC update is changed to
\begin{lstlisting}[frame=single,style=Matlab-editor,basicstyle={\fontsize{7}{8}\ttfamily}]
ocP = xT .* real( sqrt( max(1e-10,-dc( act )) ./ dg1(act) ) );
\end{lstlisting}

\begin{figure}
    \centering
    \includegraphics[width=\linewidth]{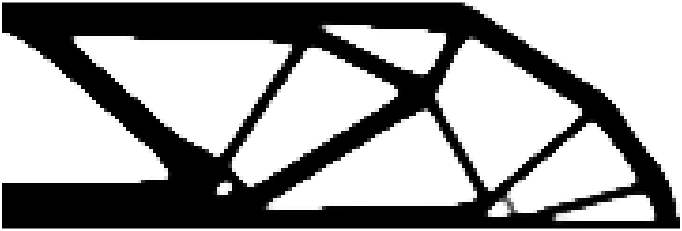}
    \includegraphics[width=\linewidth]{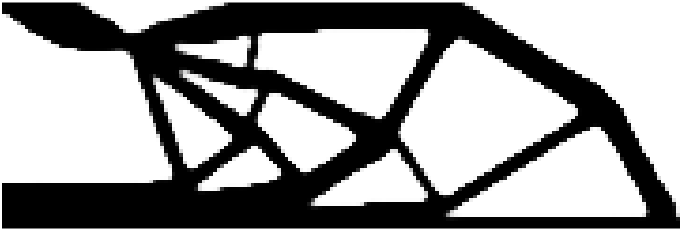}
    \includegraphics[width=\linewidth]{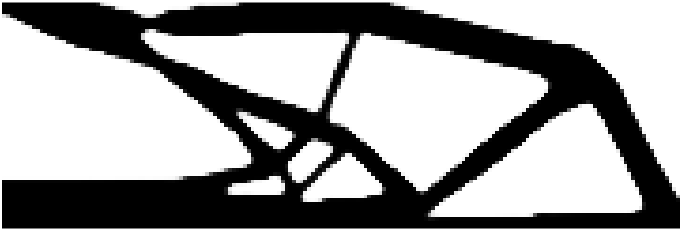}
    \caption{Final designs for mini-batch example and different $P$-norms. Top: Static load. Mid: Dynamic load with $\Vert\cdot\Vert_1$. Bottom: Dynamic load with $\Vert\cdot\Vert_{10}$.}
    \label{fig:dynamic_design}
\end{figure}

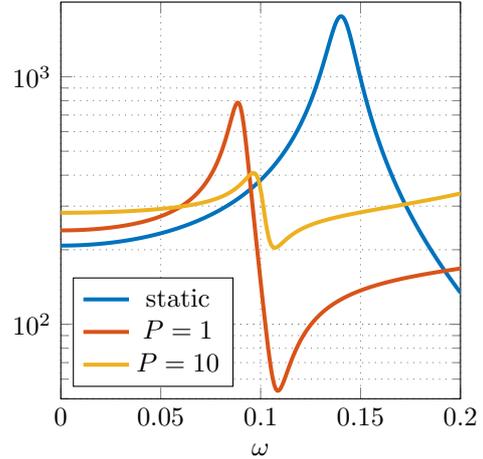
\begin{figure}
    \centering
    \input{Graphics/Dynamic/Dynamic_compliances}
    \caption{Frequency-dependent compliance values ${c(\hat{\mx}_\text{final},\omega)}$ for final designs obtained by the static load case (blue) as well as $P$-norm approaches with $P=1$ (orange) and $P=10$ (yellow) over the range of considered frequencies.}
    \label{fig:dynamic_compl}
\end{figure}

Choosing the same parameters as in~\Cref{sec:ClampNonsym}, we consider the $P$-norm parameters $P=1$ and $P=10$. For comparison, we also optimize the design for the static load case $\omega=0$ only. The corresponding final designs are shown in~\Cref{fig:dynamic_design}.
For each final design, the frequency-dependent dynamic compliance values can be found in~\Cref{fig:dynamic_compl}. As expected, we observe the following:
\begin{enumerate}
    \item The design optimized for $\omega=0$ has the lowest static compliance.
    \item The design optimized for $P=1$ has the lowest expected compliance.
    \item The design optimized for $P=10$ has the lowest maximum compliance.
\end{enumerate}
The corresponding values are given in~\Cref{tab:dyna_comparison}.

\begin{table}
    \centering
    \begin{tabular}{|l||r|r|r|}
        \hline
        \textbf{method} & $c(\hat{\mx},0)$ & $\mathbb{E}_\omega\big[c(\hat{\mx},\omega)\big]$ & $\max_{\omega}\,c(\hat{\mx},\omega)$  \\
        \hline
         static & 207.83 & 444.05 & 1754.30 \\
         \hline
         $P=1$ & 239.42 & 228.84 & 786.05 \\
         \hline
         $P=10$ & 282.02 & 294.00 & 408.53\\
         \hline
    \end{tabular}
    \caption{Dynamic compliance values for final designs. Left: Static compliance. Mid: Expected dynamic compliance over frequency range. Right: Maximum dynamic compliance over frequency range.}
    \label{tab:dyna_comparison}
\end{table}
\section{Conclusion}
This paper presents a simple yet efficient MATLAB code for topology optimization with uncertain parameters. The code is constructed as an optional add-on to the most recent descendant of the famous 88 line code, \texttt{top99neo}~\cite{99neoCode}. The proposed implementation does not represent cutting-edge research, nor is it intended as hyper-optimized code for high performance computing. It is published to provide an accessible tool for students, newcomers and other non-specialists, that is easy to modify and can be run on basically any hardware. Nonetheless, our experiments show that the method is capable of solving numerically expensive probabilistic optimization problems in short computing time.

We have given detailed motivation and theoretical background for each new step in the code. Additionally, we have provided extensive numerical results concerning several key aspects of the method, like the impact of different random sample sequences or adjustments of the stochastic model in case of symmetry. Moreover, several generalizations have been discussed. Specifically, we have considered tinkering with the sampling strategy, adapting a model reduction and incorporating mini-batch techniques.

\bmhead{Acknowledgements} A. Uihlein and M. Stingl acknowledge funding by the Deutsche Forschungsgemeinschaft (DFG, German Research Foundation) under Project-ID 416229255 (CRC 1411). O. Sigmund acknowledges financial support from the Villum Foundation through the Villum Investigator Project Amstrad (VIL54487).
\bmhead{Replication of results} MATLAB codes are given in the Appendix and are available at \url{www.topopt.dtu.dk} as well as~\cite{Zenodo}. They use the \texttt{stenglib} package for fast sparse matrix operations, which can be downloaded from \url{https://github.com/stefanengblom/stenglib}.
\bmhead{Conflict of interest} The authors declare that they have no conflict of interest.
\bibliography{sn-bibliography}
\appendix
\onecolumn
\newpage
\section{\texttt{topS140} code}\label{sec:AppendixCode}
\begin{lstlisting}[numbers=left,style=Matlab-editor,basicstyle={\fontsize{7}{8}\ttfamily}]
function topS140(nelx,nely,volfrac,penal,rmin,ft,ftBC,eta,beta,move,pnorm,maxit)
% ---------------------------- PRE. 1) MATERIAL AND CONTINUATION PARAMETERS
E0 = 1;                                                                     % Young modulus of solid
Emin = 1e-9;                                                                % Young modulus of "void"
nu = 0.3;                                                                   % Poisson ratio
penalCnt = { 3,  5, 50, 0 };                                                % continuation penal
betaCnt  = { 2,  16, 50, 1 };                                               % continuation beta
if ftBC == 'N', bcF = 'symmetric'; else, bcF = 0; end                       % filter BC selector
% ----------------------------------------- PRE. 2) DISCRETIZATION FEATURES
nEl = nelx * nely;                                                          % number of elements
nodeNrs = int32( reshape( 1 : (1 + nelx) * (1 + nely), 1+nely, 1+nelx ) );  % nodes numbers
cVec = reshape( 2 * nodeNrs( 1 : end - 1, 1 : end - 1 ) + 1, nEl, 1 );
cMat = cVec + int32( [ 0, 1, 2 * nely + [ 2, 3, 0, 1 ], -2, -1 ] );         % connectivity matrix
nDof = ( 1 + nely ) * ( 1 + nelx ) * 2;                                     % total number of DOFs
[ sI, sII ] = deal( [ ] );
for j = 1 : 8
    sI = cat( 2, sI, j : 8 );
    sII = cat( 2, sII, repmat( j, 1, 8 - j + 1 ) );
end
[ iK , jK ] = deal( cMat( :,  sI )', cMat( :, sII )' );
Iar = sort( [ iK( : ), jK( : ) ], 2, 'descend' ); clear iK jK
c1 = [12;3;-6;-3;-6;-3;0;3;12;3;0;-3;-6;-3;-6;12;-3;0;-3;-6;3;12;3;...
    -6;3;-6;12;3;-6;-3;12;3;0;12;-3;12];
c2 = [-4;3;-2;9;2;-3;4;-9;-4;-9;4;-3;2;9;-2;-4;-3;4;9;2;3;-4;-9;-2;...
    3;2;-4;3;-2;9;-4;-9;4;-4;-3;-4];
Ke = 1/(1-nu^2)/24*( c1 + nu .* c2 );                                       % lower sym. part
Ke0( tril( ones( 8 ) ) == 1 ) = Ke';
Ke0 = reshape( Ke0, 8, 8 );
Ke0 = Ke0 + Ke0' - diag( diag( Ke0 ) );                                     % recover full matrix
% ----------------------------- PRE. 3) LOADS, SUPPORTS AND PASSIVE DOMAINS
nrlo = union(nodeNrs(:,1),nodeNrs(:,end));
fixed = union(2*nrlo,2*nrlo-1);                                             % fix in both directions
[ pasS, pasV ] = deal([],[]);
free = setdiff( 1 : nDof, fixed );                                          % set of free DOFs
act = setdiff( (1 : nEl )', union( pasS, pasV ) );                          % set of active d.v.
lcDof = 2 * nodeNrs( 1, : );
F = fsparse( lcDof', 1, -1/numel(lcDof), [ nDof, 1 ] );                     % define load vector
F(lcDof(1), :)   = 0.5*F(lcDof(1), :);
F(lcDof(end), :) = 0.5*F(lcDof(end), :);
% ---------------------------------------------------------- PRE. 4) DAMAGE
[L, nonD, dmg_fac] = deal(20, 5, 1);
% ------------------------------------------------ PRE. 5) STOCHASTIC MODEL
rng('default');                                                             % reset random parameter
com0 = 100;                                                                 % initial guess (scaling)
[y1, y2] = meshgrid( 1:(nelx-L+1) , 1:(nely-L+1-nonD)  ); 
y = [y1(:),y2(:)]';                                                         % integration points
n_disc = size(y,2);                                                         % number of points
X = [randi(nelx-L+1,1,maxit) ; randi(nely-L+1-nonD,1,maxit)];               % sample sequence
maxsmpl = 2000;                                                             % max stored samples
[x_birth,x_ind,leavers] = deal(1:maxsmpl,1:maxsmpl,1);
y_weight = volfrac*sqrt(nEl);                                               % weighting norm
y_diff = pdist2(y',X');
y_diff = y_diff/max(max(y_diff,1e-10),[],'all');                            % normalization
[Gra,DesH,ComH] = deal( zeros(nEl,maxsmpl), zeros(nEl,maxsmpl), zeros(1,maxsmpl) );
% --------------------------------------- PRE. 6) DEFINE IMPLICIT FUNCTIONS
prj = @(v,eta,beta) (tanh(beta*eta)+tanh(beta*(v(:)-eta)))./...
    (tanh(beta*eta)+tanh(beta*(1-eta)));                                    % projection
deta = @(v,eta,beta) - beta * csch( beta ) .* sech( beta * ( v( : ) - eta ) ).^2 .* ...
    sinh( v( : ) * beta ) .* sinh( ( 1 - v( : ) ) * beta );
dprj = @(v,eta,beta) beta*(1-tanh(beta*(v-eta)).^2)./(tanh(beta*eta)+tanh(beta*(1-eta)));
cnt = @(v,vCnt,l) v+(l>=vCnt{1})*(v<vCnt{2})*(mod(l,vCnt{3})==0)*vCnt{4};   % apply continuation
% -------------------------------------------------- PRE. 7) PREPARE FILTER
[dy,dx] = meshgrid(-ceil(rmin)+1:ceil(rmin)-1,-ceil(rmin)+1:ceil(rmin)-1);
h = max( 0, rmin - sqrt( dx.^2 + dy.^2 ) );                                 % conv. kernel
Hs = imfilter( ones( nely, nelx ), h, bcF );                                % filter weights
dHs = Hs;
% ------------------------ PRE. 8) ALLOCATE AND INITIALIZE OTHER PARAMETERS
[ x, dsK, dV ] = deal( zeros( nEl, 1 ) );                                   % initialize vectors
dV( act, 1 ) = 1/nEl/volfrac;                                               % derivative of volume
x( act ) = ( volfrac*( nEl - length(pasV) ) - length(pasS) )/length( act ); % volume fraction
x( pasS ) = 1;                                                              % set x = 1 on pasS set
[ xPhys, loop, U ] = deal( x, 0, zeros( nDof, 1 ) );                        % it. counter, U
% ================================================= START OPTIMIZATION LOOP
while loop < maxit 
    loop = loop + 1;                                                        % update iter. counter
    % --------- RL. 1) COMPUTE PHYSICAL DENSITY FIELD (AND ETA IF PROJECT.)
    xTilde = imfilter( reshape( x, nely, nelx ), h, bcF ) ./ Hs;
    xPhys( act ) = xTilde( act );  
    if ft > 1                              % compute optimal eta* with Newton
        f = ( mean( prj( xPhys, eta, beta ) ) - volfrac ) * ( ft == 3 );    % function (volume)
        while abs( f ) > 1e-6           % Newton process for finding opt. eta
            eta = eta - f / mean( deta( xPhys( : ), eta, beta ) );
            f = mean( prj( xPhys, eta, beta ) ) - volfrac;
        end
        dHs = Hs ./ reshape( dprj( xTilde, eta, beta ), nely, nelx );       % modification sensitivity
        xPhys = prj( xPhys, eta, beta );                                    % projected (phys.) field
    end
    % ------------------------ RL. 2) SETUP AND SOLVE EQUILIBRIUM EQUATIONS
    D = zeros(nely,nelx);
    D( nely+1-(X(2,loop):X(2,loop)+L-1) , X(1,loop):X(1,loop)+L-1 ) = 1;    % damage cell
    x_dmg = max(0,min(1,xPhys-dmg_fac*D(:)));
    sK = ( Emin + x_dmg.^penal * ( E0 - Emin ) );                           % stiffness interpolation
    dsK( act ) = -penal * ( E0 - Emin ) * x_dmg( act ) .^ ( penal - 1 );
    sK = reshape( Ke( : ) * sK', length( Ke ) * nEl, 1 );
    K = fsparse( Iar( :, 1 ), Iar( :, 2 ), sK, [ nDof, nDof ] );            % stiffness matrix
    U(free) = decomposition( K( free, free ), 'chol','lower' ) \ F(free);   % solve equilibrium system
    % ---------------------------------------- RL. 3) COMPUTE SENSITIVITIES
    dV0 = imfilter( reshape( dV, nely, nelx ) ./ dHs, h, bcF );
    dc = reshape(dsK .* sum( (U( cMat ) * Ke0 ) .* U( cMat ), 2 ), nely, nelx);
    dc = imfilter( dc ./ dHs, h, bcF );
    % -------------------- RL. 4) SAMPLE MANAGEMENT AND INTEGRATION WEIGHTS
    ulim = min(maxsmpl,loop);
    [Gra(:,leavers), ComH(:,leavers), DesH(:,leavers)] = deal( dc(:), F'*U, xPhys );
    [~,csw] = min( vecnorm(xPhys-DesH(:,1:ulim),2,1) + y_weight*y_diff(:,x_ind(1:ulim)), [], 2);
    weights = sum(csw==1:ulim)/n_disc;                                      % integration weights
    if (loop+1) <= maxsmpl
        leavers = loop+1;
    else
        ind_can = find(weights-min(weights)<1e-8);                          % small weights
        [~,iind] = min(x_birth(ind_can));                                   % oldest of these samples
        leavers = ind_can(iind);
        [x_ind(leavers), x_birth(leavers)] = deal(loop+1,loop+1);
    end
    % ------------------------------ RL. 5) NEAREST NEIGHBOR APPROXIMATIONS
    Compl = sum(weights.*ComH(:,1:ulim),2);Cp = (sum(weights.*(ComH(:,1:ulim).^pnorm),2))^(1/pnorm);
    if mod(loop,25) == 0
        com0 = Compl;                                                       % adjust normalization
    end
    dc = 1/com0*sum(weights.*(((ComH(:,1:ulim)/com0).^(pnorm-1)).*Gra(:,1:ulim)),2);
    % --------------- RL. 6) UPDATE DESIGN VARIABLES AND APPLY CONTINUATION
    xT = x( act );
    [ xU, xL ] = deal( xT + move, xT - move );                              % current bounds
    ocP = xT .* real( sqrt( -dc( act ) ./ dV0( act ) ) );
    l = [ 0, (mean( ocP ) / volfrac) ];                                     % initial estimate for LM
    while ( l( 2 ) - l( 1 ) ) / ( l( 2 ) + l( 1 ) ) > 1e-8  && l(2) > 1e-40 % OC resizing rule
        lmid = 0.5 * ( l( 1 ) + l( 2 ) );
        x( act ) = max( max( min( min( ocP / lmid, xU ), 1 ), xL ), 0 );
        xTilde = imfilter( reshape( x, nely, nelx ), h, bcF ) ./ Hs;
        xPhys( act ) = xTilde( act );
        xPhys = prj( xPhys, eta, beta );
        if mean( xPhys ) > volfrac, l( 1 ) = lmid; else, l( 2 ) = lmid; end
    end
    [penal,beta] = deal(cnt(penal,penalCnt,loop), cnt(beta,betaCnt,loop));  % apply conitnuation
    % ------------------------ RL. 7) PRINT CURRENT RESULTS AND PLOT DESIGN
    fprintf( 'It.:%5i C:%7.4f Cp:%7.4f V:%7.3f penal:%7.2f beta:%7.1f eta:%7.2f \n', ...
        loop, Compl, Cp, mean(xPhys), penal, beta, eta);
    colormap( gray ); imagesc( 1 - reshape( xPhys, nely, nelx ) );
    caxis([0 1]); axis equal off; drawnow;
end
end
\end{lstlisting}
\newpage
\section{\texttt{topS140\_load} code}\label{sec:AppendixCode_load}
\begin{lstlisting}[numbers=left,style=Matlab-editor,basicstyle={\fontsize{7}{8}\ttfamily}]
function topS140_load(nelx,nely,volfrac,penal,rmin,ft,ftBC,eta,beta,move,pnorm,maxit,type)
% ---------------------------- PRE. 1) MATERIAL AND CONTINUATION PARAMETERS
E0 = 1;                                                                     % Young modulus of solid
Emin = 1e-9;                                                                % Young modulus of "void"
nu = 0.3;                                                                   % Poisson ratio
penalCnt = { 3,  5, 50, 0 };                                                % continuation penal
betaCnt  = { 2,  16, 75, 1 };                                               % continuation beta
if ftBC == 'N', bcF = 'symmetric'; else, bcF = 0; end                       % filter BC selector
% ----------------------------------------- PRE. 2) DISCRETIZATION FEATURES
nEl = nelx * nely;                                                          % number of elements
nodeNrs = int32( reshape( 1 : (1 + nelx) * (1 + nely), 1+nely, 1+nelx ) );  % nodes numbers
cVec = reshape( 2 * nodeNrs( 1 : end - 1, 1 : end - 1 ) + 1, nEl, 1 );
cMat = cVec + int32( [ 0, 1, 2 * nely + [ 2, 3, 0, 1 ], -2, -1 ] );         % connectivity matrix
nDof = ( 1 + nely ) * ( 1 + nelx ) * 2;                                     % total number of DOFs
[ sI, sII ] = deal( [ ] );
for j = 1 : 8
    sI = cat( 2, sI, j : 8 );
    sII = cat( 2, sII, repmat( j, 1, 8 - j + 1 ) );
end
[ iK , jK ] = deal( cMat( :,  sI )', cMat( :, sII )' );
Iar = sort( [ iK( : ), jK( : ) ], 2, 'descend' ); clear iK jK
c1 = [12;3;-6;-3;-6;-3;0;3;12;3;0;-3;-6;-3;-6;12;-3;0;-3;-6;3;12;3;...
    -6;3;-6;12;3;-6;-3;12;3;0;12;-3;12];
c2 = [-4;3;-2;9;2;-3;4;-9;-4;-9;4;-3;2;9;-2;-4;-3;4;9;2;3;-4;-9;-2;...
    3;2;-4;3;-2;9;-4;-9;4;-4;-3;-4];
Ke = 1/(1-nu^2)/24*( c1 + nu .* c2 );                                       % lower sym. part
Ke0( tril( ones( 8 ) ) == 1 ) = Ke';
Ke0 = reshape( Ke0, 8, 8 );
Ke0 = Ke0 + Ke0' - diag( diag( Ke0 ) );                                     % recover full matrix
% ------------------------------------ PRE. 3) SUPPORTS AND PASSIVE DOMAINS
nrlo = union(nodeNrs(round(nely/2):end,1),nodeNrs(round(nely/2):end,end));
fixed = union(2*nrlo,2*nrlo-1);                                             % fix in both directions
[ pasS, pasV ] = deal(1:nely:nelx*nely,[]);
free = setdiff( 1 : nDof, fixed );                                          % set of free DOFs
act = setdiff( (1 : nEl )', union( pasS, pasV ) );                          % set of active d.v.
% ----------------------------------------------------------- PRE. 4) LOADS
lcDof = 2 * nodeNrs(1,:);
% ------------------------------------------------ PRE. 5) STOCHASTIC MODEL
rng('default');                                                             % reset random parameter
com0 = 100;                                                                 % initial guess (scaling)
y_all = load('Data.mat','r'); y_all = y_all.r; 
y = unique(y_all); n_disc = size(y,2);                                      % integration points
dist_w = sum(y_all' == y)/size(y_all,2);                                    % weights quadrature rule
switch type
    case 'distribution' % sample according to distribution
        X = y_all(randi(numel(y_all),1,maxit));
    case 'uniform'  % sample uniformly
        X = y(randi(n_disc,1,maxit));
end
maxsmpl = 2000;                                                             % max stored samples
[x_birth,x_ind,leavers] = deal(1:maxsmpl,1:maxsmpl,1);
y_weight = 5*volfrac*sqrt(nEl);                                             % weighting norm
y_diff = pdist2(y',X');
y_diff = y_diff/max(max(y_diff,1e-10),[],'all');                            % normalization
X = int32(X);
[Gra,DesH,ComH] = deal( zeros(nEl,maxsmpl), zeros(nEl,maxsmpl), zeros(1,maxsmpl) );
% --------------------------------------- PRE. 6) DEFINE IMPLICIT FUNCTIONS
prj = @(v,eta,beta) (tanh(beta*eta)+tanh(beta*(v(:)-eta)))./...
    (tanh(beta*eta)+tanh(beta*(1-eta)));                                    % projection
deta = @(v,eta,beta) - beta * csch( beta ) .* sech( beta * ( v( : ) - eta ) ).^2 .* ...
    sinh( v( : ) * beta ) .* sinh( ( 1 - v( : ) ) * beta );
dprj = @(v,eta,beta) beta*(1-tanh(beta*(v-eta)).^2)./(tanh(beta*eta)+tanh(beta*(1-eta)));
cnt = @(v,vCnt,l) v+(l>=vCnt{1})*(v<vCnt{2})*(mod(l,vCnt{3})==0)*vCnt{4};   % apply continuation
% -------------------------------------------------- PRE. 7) PREPARE FILTER
[dy,dx] = meshgrid(-ceil(rmin)+1:ceil(rmin)-1,-ceil(rmin)+1:ceil(rmin)-1);
h = max( 0, rmin - sqrt( dx.^2 + dy.^2 ) );                                 % conv. kernel
Hs = imfilter( ones( nely, nelx ), h, bcF );                                % filter weights
dHs = Hs;
% ------------------------ PRE. 8) ALLOCATE AND INITIALIZE OTHER PARAMETERS
[ x, dsK, dV ] = deal( zeros( nEl, 1 ) );                                   % initialize vectors
dV( act, 1 ) = 1/nEl/volfrac;                                               % derivative of volume
x( act ) = ( volfrac*( nEl - length(pasV) ) - length(pasS) )/length( act ); % volume fraction
x( pasS ) = 1;                                                              % set x = 1 on pasS set
[ xPhys, loop, U ] = deal( x, 0, zeros( nDof, 1 ) );                        % it. counter, U
% ================================================= START OPTIMIZATION LOOP
while loop < maxit 
    loop = loop + 1;                                                        % update iter. counter
    % --------- RL. 1) COMPUTE PHYSICAL DENSITY FIELD (AND ETA IF PROJECT.)
    xTilde = imfilter( reshape( x, nely, nelx ), h, bcF ) ./ Hs;
    xPhys( act ) = xTilde( act );  
    if ft > 1                              % compute optimal eta* with Newton
        f = ( mean( prj( xPhys, eta, beta ) ) - volfrac ) * ( ft == 3 );    % function (volume)
        while abs( f ) > 1e-6           % Newton process for finding opt. eta
            eta = eta - f / mean( deta( xPhys( : ), eta, beta ) );
            f = mean( prj( xPhys, eta, beta ) ) - volfrac;
        end
        dHs = Hs ./ reshape( dprj( xTilde, eta, beta ), nely, nelx );       % modification sensitivity
        xPhys = prj( xPhys, eta, beta );                                    % projected (phys.) field
    end
    % ------------------------ RL. 2) SETUP AND SOLVE EQUILIBRIUM EQUATIONS
    F = fsparse( lcDof(X(:,loop)), 1, -1, [ nDof, 1 ] );
    sK = ( Emin + xPhys.^penal * ( E0 - Emin ) );                           % stiffness interpolation
    dsK( act ) = -penal * ( E0 - Emin ) * xPhys( act ) .^ ( penal - 1 );
    sK = reshape( Ke( : ) * sK', length( Ke ) * nEl, 1 );
    K = fsparse( Iar( :, 1 ), Iar( :, 2 ), sK, [ nDof, nDof ] );            % stiffness matrix
    U(free) = decomposition( K( free, free ), 'chol','lower' ) \ F(free);   % solve equilibrium system
    % ---------------------------------------- RL. 3) COMPUTE SENSITIVITIES
    dV0 = imfilter( reshape( dV, nely, nelx ) ./ dHs, h, bcF );
    dc = reshape(dsK .* sum( (U( cMat ) * Ke0 ) .* U( cMat ), 2 ), nely, nelx);
    dc = imfilter( dc ./ dHs, h, bcF );
    % -------------------- RL. 4) SAMPLE MANAGEMENT AND INTEGRATION WEIGHTS
    ulim = min(maxsmpl,loop);
    [Gra(:,leavers), ComH(:,leavers), DesH(:,leavers)] = deal( dc(:), F'*U, xPhys );
    [~,csw] = min( vecnorm(xPhys-DesH(:,1:ulim),2,1) + y_weight*y_diff(:,x_ind(1:ulim)), [], 2);
    weights = sum((csw==1:ulim).*dist_w');                                  % integration weights
    if (loop+1) <= maxsmpl
        leavers = loop+1;
    else
        ind_can = find(weights-min(weights)<1e-8);                          % small weights
        [~,iind] = min(x_birth(ind_can));                                   % oldest of these samples
        leavers = ind_can(iind);
        [x_ind(leavers), x_birth(leavers)] = deal(loop+1,loop+1);
    end
    % ------------------------------ RL. 5) NEAREST NEIGHBOR APPROXIMATIONS
    Compl = sum(weights.*ComH(:,1:ulim),2);Cp = (sum(weights.*(ComH(:,1:ulim).^pnorm),2))^(1/pnorm);
    if mod(loop,25) == 0
        com0 = Compl;                                                       % adjust normalization
    end
    dc = 1/com0*sum(weights.*(((ComH(:,1:ulim)/com0).^(pnorm-1)).*Gra(:,1:ulim)),2);
    % --------------- RL. 6) UPDATE DESIGN VARIABLES AND APPLY CONTINUATION
    xT = x( act );
    [ xU, xL ] = deal( xT + move, xT - move );                              % current bounds
    ocP = xT .* real( sqrt( -dc( act ) ./ dV0( act ) ) );
    l = [ 0, 10*(mean( ocP ) / volfrac) ];                                  % initial estimate for LM
    while ( l( 2 ) - l( 1 ) ) / ( l( 2 ) + l( 1 ) ) > 1e-8  && l(2) > 1e-40 % OC resizing rule
        lmid = 0.5 * ( l( 1 ) + l( 2 ) );
        x( act ) = max( max( min( min( ocP / lmid, xU ), 1 ), xL ), 0 );
        xTilde = imfilter( reshape( x, nely, nelx ), h, bcF ) ./ Hs;
        xPhys( act ) = xTilde( act );
        xPhys = prj( xPhys, eta, beta );
        if mean( xPhys ) > volfrac, l( 1 ) = lmid; else, l( 2 ) = lmid; end
    end
    [penal,beta] = deal(cnt(penal,penalCnt,loop), cnt(beta,betaCnt,loop));  % apply conitnuation
    % ------------------------ RL. 7) PRINT CURRENT RESULTS AND PLOT DESIGN
    fprintf( 'It.:%5i C:%7.4f Cp:%7.4f V:%7.3f penal:%7.2f beta:%7.1f eta:%7.2f \n', ...
        loop, Compl, Cp, mean(xPhys), penal, beta, eta);
    colormap( gray ); imagesc( 1 - reshape( xPhys, nely, nelx ) );
    caxis([0 1]); axis equal off; drawnow
end
end
\end{lstlisting}
\end{document}

%% file: Graphics/NN_exmpl.tex
\definecolor{mycolor1}{rgb}{0.00000,0.44700,0.74100}%
\definecolor{mycolor2}{rgb}{0.6350, 0.0780, 0.1840}%
\definecolor{mycolor3}{rgb}{0.92900,0.69400,0.12500}%
\definecolor{mycolor4}{rgb}{0.4940, 0.1840, 0.5560}%
\definecolor{mycolor5}{rgb}{0.4660, 0.6740, 0.1880}%

\begin{tikzpicture}
\begin{axis}[%
width=.9\linewidth,
height=.9\linewidth,
xmin=0,
xmax=2,
ymin=0,
ymax=1.5,
xlabel = {$\Xi$},
xtick = {.2,.65,1.1,1.65},
ytick = {.75,1,.45,.3},
yticklabels = {{$g_3$},{$g_2$},{$g_1$},{$g_4$}},
xticklabels = {{\color{mycolor2}$\mathcal{V}_{3,4}$},{\color{mycolor3}$\mathcal{V}_{2,4}$},{\color{mycolor4}$\mathcal{V}_{1,4}$},{\color{mycolor5}$\mathcal{V}_{4,4}$}},
xtick style = {draw=none},
ytick style = {draw=none},
after end axis/.code={
   \draw[thick,latex-latex,mycolor2] (axis cs:0,0) -- (axis cs:0.4,0);
   \draw[thick,latex-latex,mycolor3] (axis cs:0.4,0) -- (axis cs:0.9,0);
   \draw[thick,latex-latex,mycolor4] (axis cs:0.9,0) -- (axis cs:1.3,0);
   \draw[thick,latex-latex,mycolor5] (axis cs:1.3,0) -- (axis cs:2.0,0);
                    }
]
\filldraw[draw=none,fill=mycolor2!15] (0,0) rectangle (0.4,0.75);
\filldraw[draw=none,fill=mycolor3!15] (0.4,0) rectangle (0.9,1);
\filldraw[draw=none,fill=mycolor4!15] (0.9,0) rectangle (1.3,.45);
\filldraw[draw=none,fill=mycolor5!15] (1.3,0) rectangle (2,.3);
\draw[dotted,line width=.75pt] (0,1) -- (0.4,1);
\draw[dotted,line width=.75pt] (0,.45) -- (0.9,.45);
\draw[dotted,line width=.75pt] (0,.3) -- (1.3,.3);

\addplot[domain=0:2,smooth,samples=101,mycolor1,line width=1.75pt]{exp(-4/5*x)*sin(3*deg(x))+.5};
\addlegendentry{$\nabla_\mx c(\hat{\mx}_4,\xi)$};
\addplot[domain=0:0.4,smooth,samples=5,mycolor2,line width=1.5pt,line join=round,forget plot]{.75};
\addplot[domain=0.4:0.9,smooth,samples=5,mycolor3,line width=1.5pt,line join=round,forget plot]{1};
\addplot[domain=0.9:1.3,smooth,samples=5,mycolor4,line width=1.5pt,line join=round,forget plot]{.45};
\addplot[domain=1.3:2,smooth,samples=5,mycolor5,line width=1.5pt,line join=round,forget plot]{.3};
\end{axis}    
\end{tikzpicture}

%% file: Graphics/voronoi.tex
\definecolor{mycolor6}{rgb}{0.30100,0.74500,0.93300}%
\definecolor{mycolor3}{rgb}{0.46600,0.67400,0.18800}%
\definecolor{mycolor4}{rgb}{0.85000,0.32500,0.09800}%
\definecolor{mycolor2}{rgb}{0.00000,0.44700,0.74100}%
\definecolor{mycolor5}{rgb}{0.49400,0.18400,0.55600}%
\definecolor{mycolor1}{rgb}{0.92900,0.69400,0.12500}%
\definecolor{mycolor7}{rgb}{0.6350 0.0780 0.1840}%
\begin{tikzpicture}

\begin{axis}[%
width=\linewidth,
height=\linewidth,
at={(0.758in,0.481in)},
xmin=0,
xmax=1,
xlabel={$\Xi$},
ymin=0,
ymax=1,
ylabel={$\Psi$},
ticks = none,
axis background/.style={fill=white},
x post scale=-1, 
y post scale=1
]
\addplot[area legend, draw=white!40!black, fill=olive, fill opacity=0.3, forget plot]
table[row sep=crcr] {%
x	y\\
0.191313138760448	0.443152008253133\\
0.205760196793843	0.529323355692337\\
-4.4257077073758	5.89955236009779\\
-9.97380095830437	6.66133814775094e-16\\
-8.38886821693578	-1.4983562408777\\
0.130233581282631	0.376603751855166\\
}--cycle;

\addplot[area legend, draw=white!40!black, fill=brown, fill opacity=0.3, forget plot]
table[row sep=crcr] {%
x	y\\
0.490039018510314	0.370132333705433\\
0.648679251655518	0.434941346019338\\
0.68316377793438	0.552884969202744\\
0.456283037232469	0.804394918432142\\
0.205760196793843	0.529323355692337\\
0.191313138760448	0.443152008253133\\
}--cycle;

\addplot[area legend, draw=white!40!black, fill=mycolor7, fill opacity=0.3, forget plot]
table[row sep=crcr] {%
x	y\\
0.68316377793438	0.552884969202744\\
1.15606421903481	0.714449324866466\\
7.19135423645013	3.57567636466009\\
-1.11022302462516e-16	10.827058656007\\
-0.931615644360846	9.73666998860743\\
0.456283037232469	0.804394918432142\\
}--cycle;

\addplot[area legend, draw=white!40!black, fill=teal, fill opacity=0.3, forget plot]
table[row sep=crcr] {%
x	y\\
0.205760196793843	0.529323355692337\\
-4.4257077073758	5.89955236009779\\
-0.931615644360846	9.73666998860743\\
0.456283037232469	0.804394918432142\\
}--cycle;

\addplot[area legend, draw=white!40!black, fill=mycolor4, fill opacity=0.3, forget plot]
table[row sep=crcr] {%
x	y\\
0.425660854787906	0.187523202803647\\
0.837927013728227	-0.0237491576774736\\
0.797079135971164	0.12493449036376\\
0.648679251655518	0.434941346019338\\
0.490039018510314	0.370132333705433\\
}--cycle;

\addplot[area legend, draw=white!40!black, fill=mycolor3, fill opacity=0.3, forget plot]
table[row sep=crcr] {%
x	y\\
0.425660854787906	0.187523202803647\\
0.837927013728227	-0.0237491576774736\\
3.86402170560819	-6.25515806856365\\
5.55111512312578e-17	-9.90480035640677\\
-0.943269227939864	-8.92497930939605\\
0.290400057561397	0.0985642075814856\\
}--cycle;

\addplot[area legend, draw=white!40!black, fill=mycolor6, fill opacity=0.3, forget plot]
table[row sep=crcr] {%
x	y\\
0.797079135971164	0.12493449036376\\
0.837927013728227	-0.0237491576774736\\
3.86402170560819	-6.25515806856365\\
11.024258288727	0\\
7.19135423645013	3.57567636466009\\
1.15606421903481	0.714449324866466\\
}--cycle;

\addplot[area legend, draw=white!40!black, fill=mycolor5, fill opacity=0.3, forget plot]
table[row sep=crcr] {%
x	y\\
0.648679251655518	0.434941346019338\\
0.797079135971164	0.12493449036376\\
1.15606421903481	0.714449324866466\\
0.68316377793438	0.552884969202744\\
}--cycle;

\addplot[area legend, draw=white!40!black, fill=mycolor2, fill opacity=0.3, forget plot]
table[row sep=crcr] {%
x	y\\
0.425660854787906	0.187523202803647\\
0.490039018510314	0.370132333705433\\
0.191313138760448	0.443152008253133\\
0.130233581282631	0.376603751855166\\
0.290400057561397	0.0985642075814856\\
}--cycle;

\addplot[area legend, draw=white!40!black, fill=mycolor1, fill opacity=0.3, forget plot]
table[row sep=crcr] {%
x	y\\
0.290400057561397	0.0985642075814856\\
-0.943269227939864	-8.92497930939605\\
-8.38886821693578	-1.4983562408777\\
0.130233581282631	0.376603751855166\\
}--cycle;

\addplot[black,dashed]
table[row sep=crcr] {%
x y\\
1 1\\
0 0\\
};
\addplot[black,dashed]
table[row sep=crcr] {%
x y\\
0.75 1\\
0 0.25\\
};
\addplot[black,dashed]
table[row sep=crcr] {%
x y\\
0.5 1\\
0 0.5\\
};
\addplot[black,dashed]
table[row sep=crcr] {%
x y\\
0.25 1\\
0 0.75\\
};
\addplot[black,dashed]
table[row sep=crcr] {%
x y\\
1 0.75\\
0.25 0\\
};
\addplot[black,dashed]
table[row sep=crcr] {%
x y\\
1 0.5\\
0.5 0\\
};
\addplot[black,dashed]
table[row sep=crcr] {%
x y\\
1 0.25\\
0.75 0\\
};
\addplot[black,dashed]
table[row sep=crcr] {%
x y\\
1 0\\
0 1\\
};
\addplot[black,dashed]
table[row sep=crcr] {%
x y\\
1 0.25\\
0.25 1\\
};
\addplot[black,dashed]
table[row sep=crcr] {%
x y\\
1 0.5\\
0.5 1\\
};
\addplot[black,dashed]
table[row sep=crcr] {%
x y\\
1 0.75\\
0.75 1\\
};
\addplot[black,dashed]
table[row sep=crcr] {%
x y\\
0.75 0\\
0 0.75\\
};
\addplot[black,dashed]
table[row sep=crcr] {%
x y\\
0.5 0\\
0 0.5\\
};
\addplot[black,dashed]
table[row sep=crcr] {%
x y\\
0.25 0\\
0 0.25\\
};

\addplot [black,only marks, mark=*, mark size = 2.5, mark options={solid, fill=mycolor1}, forget plot]
  table[row sep=crcr]{%
0 0\\
0 0.25\\
0.25 0\\
0.125 0.125\\
};
\addplot [black,only marks, mark=*, mark size = 2.5, mark options={solid, fill=mycolor3}, forget plot]
  table[row sep=crcr]{%
0.5 0\\
0.75 0\\
0.375 0.125\\
};
\addplot [black,only marks, mark=*, mark size = 2.5, mark options={solid, fill=mycolor6}, forget plot]
  table[row sep=crcr]{%
1 0\\
1 0.25\\
0.875 0.125\\
};
\addplot [black,only marks, mark=*, mark size = 2.5, mark options={solid, fill=mycolor5}, forget plot]
  table[row sep=crcr]{%
1 0.5\\
0.75 0.5\\
0.75 0.25\\
0.875 0.375\\
};
\addplot [black,only marks, mark=*, mark size = 2.5, mark options={solid, fill=mycolor4}, forget plot]
  table[row sep=crcr]{%
0.5 0.25\\
0.625 0.125\\
0.625 0.375\\
};
\addplot [black,only marks, mark=*, mark size = 2.5, mark options={solid, fill=mycolor2}, forget plot]
  table[row sep=crcr]{%
0.25 0.25\\
0.375 0.375\\
};
\addplot [black,only marks, mark=*, mark size = 2.5, mark options={solid, fill=olive}, forget plot]
  table[row sep=crcr]{%
0 0.5\\
0 0.75\\
0.125 0.375\\
};
\addplot [black,only marks, mark=*, mark size = 2.5, mark options={solid, fill=brown}, forget plot]
  table[row sep=crcr]{%
0.5 0.5\\
0.25 0.5\\
0.5 0.75\\
0.375 0.625\\
};
\addplot [black,only marks, mark=*, mark size = 2.5, mark options={solid, fill=teal}, forget plot]
  table[row sep=crcr]{%
0 1\\
0.25 1\\
0.25 0.75\\
0.375 0.875\\
0.125 0.875\\
0.125 0.625\\
};
\addplot [black,only marks, mark=*, mark size = 2.5, mark options={solid, fill=mycolor7}, forget plot]
  table[row sep=crcr]{%
1 0.75\\
1 1\\
0.75 0.75\\
0.75 1\\
0.5 1\\
0.875 0.875\\
0.625 0.875\\
0.875 0.625\\
0.625 0.625\\
};
\end{axis}
\end{tikzpicture}%

%% file: Graphics/Clamp/clamp_setup.tex
\begin{tikzpicture}[x=0.9\linewidth/200,y=0.9\linewidth/200]
\draw[very thick] (10,5) -- (190,5) -- (190,50) -- (10,50) -- cycle;
\fill[pattern=north west lines] (0,5) rectangle (10,50);
\fill[pattern=north west lines] (190,5) rectangle (200,50);
\fill[pattern=north east lines,pattern color = blue] (35,10) rectangle (55,30);
\draw[blue,thick] (35,10) rectangle (55,30);
\draw[Latex-Latex] (10,0)--node[midway,below]{$4\ell$}(190,0);
\draw[Latex-Latex] (-5,5)--node[left,midway]{$\ell$}(-5,50);
\foreach \idx in {2,...,38}
{
    \draw[Latex-,red] (5*\idx,52) -- (5*\idx,65);
}
\draw[thin,red] (10,52) -- (190,52);
\end{tikzpicture}

%% file: Graphics/Clamp/Reference_obj.tex
%
%
\definecolor{mycolor1}{rgb}{0.00000,0.44700,0.74100}%
\definecolor{mycolor2}{rgb}{0.85000,0.32500,0.09800}%
\begin{tikzpicture}

\begin{axis}[%
width=.9\linewidth,
height=.9\linewidth,
at={(0.758in,0.481in)},
xmin=0,
xmax=1500,
ymin=5,
ymax=55,
xlabel={Iteration},
axis background/.style={fill=white},
xmajorgrids,
ymajorgrids
]
\addplot [color=mycolor1, line width = 1pt, line join=round]
  table[row sep=crcr]{%
1	49.9027306485921\\
2	50.4201310954899\\
3	46.7034546855003\\
4	45.760813770353\\
5	39.488363695647\\
6	39.0821726652808\\
7	37.630045176397\\
8	35.9278520020072\\
9	33.6706800956268\\
10	33.2928124695405\\
11	31.5932417387744\\
12	30.7075088603052\\
13	30.7318405946205\\
14	27.3027070847228\\
15	26.6481736777287\\
16	25.3781467957875\\
17	25.6382050147258\\
18	25.8841471717047\\
19	25.139923388993\\
20	25.208482387779\\
21	24.5472374926041\\
22	25.7155612717482\\
23	27.8252609041464\\
24	26.7304536244628\\
25	27.5363218429388\\
26	26.0154818203097\\
27	25.9017425530747\\
28	22.5913399236615\\
29	21.5984615347311\\
30	22.0321363271273\\
31	20.473779396016\\
32	20.9655976622473\\
33	20.086212293403\\
34	20.473822466289\\
35	19.9000207380029\\
36	24.4478608942389\\
37	27.3329503787088\\
38	26.7759651013753\\
39	26.334947841752\\
40	25.1762663612581\\
41	22.1492189653336\\
42	21.9631302143327\\
43	19.3790492513636\\
44	18.0387161335574\\
45	20.8433179119935\\
46	19.7576829795284\\
47	19.2408664839872\\
48	18.6663990792199\\
49	19.1451103109751\\
50	17.0557528105543\\
51	13.6022795356252\\
52	11.3358685411162\\
53	12.2787111429064\\
54	13.8811812371645\\
55	14.3800457002078\\
56	13.3293922863684\\
57	13.0321722207924\\
58	12.4578067970728\\
59	11.9725590170085\\
60	10.4937630506097\\
61	11.6282951919304\\
62	17.5068154247429\\
63	14.9840428054059\\
64	13.7124295752068\\
65	17.4021196650025\\
66	17.8075273215738\\
67	17.4994579337694\\
68	15.3613183894854\\
69	16.1242130528906\\
70	15.9016808988248\\
71	15.3465168166242\\
72	13.8729735930234\\
73	12.8504120434624\\
74	12.225785509523\\
75	11.2426295773689\\
76	12.0039256924387\\
77	12.2784614434125\\
78	13.6321917872159\\
79	12.5057928204006\\
80	11.9297755689179\\
81	11.1009372582516\\
82	10.530724381636\\
83	10.5095663880851\\
84	10.5211191436394\\
85	15.8461679233685\\
86	14.3520535742941\\
87	13.2909786019341\\
88	12.4804210103546\\
89	12.6366094327737\\
90	12.9383878412571\\
91	12.8218331301969\\
92	12.7755381206076\\
93	12.7133246623626\\
94	12.4766741274515\\
95	13.2156057479958\\
96	13.3789121435233\\
97	12.7606550745973\\
98	12.5763333660248\\
99	12.5115160376184\\
100	12.051714147807\\
101	12.0860327713298\\
102	10.1791906897246\\
103	10.9357965271215\\
104	10.4950089998673\\
105	11.5033023847582\\
106	11.1996755502219\\
107	10.7566453911343\\
108	11.1097367945438\\
109	11.212251284313\\
110	9.52323715133533\\
111	9.35376148337215\\
112	10.2260986761484\\
113	13.6994985529107\\
114	13.504517087502\\
115	12.3103444011178\\
116	13.5682882795487\\
117	12.8184563591271\\
118	12.0080423045299\\
119	11.033992661944\\
120	10.4229653603591\\
121	10.469187806243\\
122	11.9589067174481\\
123	11.6667286287004\\
124	10.2961665501861\\
125	10.1130930573127\\
126	11.6929964566645\\
127	12.4551510866026\\
128	11.6963185566241\\
129	14.0335518944306\\
130	13.0360136411279\\
131	12.6057123430454\\
132	12.9411471871674\\
133	10.6735473832561\\
134	10.6584037760906\\
135	11.4450236845261\\
136	11.3277815014624\\
137	12.5354338865313\\
138	11.4192251968821\\
139	11.3887864790179\\
140	12.0333089895438\\
141	11.3129526391528\\
142	10.9820362572561\\
143	10.1568372339365\\
144	10.2927242419946\\
145	10.6419891692875\\
146	10.922692431985\\
147	11.0678771029069\\
148	10.6585787497065\\
149	10.6210668876995\\
150	10.8230023380922\\
151	10.1934517022894\\
152	10.2494131306464\\
153	10.8583067505107\\
154	10.5277536495069\\
155	10.2777313854598\\
156	10.1174965640626\\
157	9.79626537654742\\
158	13.2075323214883\\
159	12.7098136440677\\
160	12.4575366858629\\
161	11.6681653334518\\
162	11.8412801420596\\
163	11.4593705030392\\
164	11.1719461568825\\
165	12.1691765340535\\
166	11.7304041875477\\
167	10.2447384214183\\
168	10.4634912012652\\
169	9.771513959026\\
170	9.7093444138481\\
171	11.8050628896437\\
172	11.2293958656086\\
173	11.5947683927614\\
174	11.1294145310954\\
175	11.4692769792446\\
176	11.0840052111032\\
177	10.3082724681932\\
178	9.49852759478016\\
179	10.35642108761\\
180	9.77896884652748\\
181	9.8525301361288\\
182	11.9119124039618\\
183	13.0731072436028\\
184	11.9091491507637\\
185	10.580686292173\\
186	10.7058288491422\\
187	10.3044675515873\\
188	10.3402246053072\\
189	9.03963369137552\\
190	9.5867178150148\\
191	11.8097112876447\\
192	13.2837966273433\\
193	12.8011055900092\\
194	12.3334861317217\\
195	12.1862240284316\\
196	12.4693138658621\\
197	12.2787319170474\\
198	12.1477868066528\\
199	11.2227122723097\\
200	10.3176166565684\\
201	10.7854168535197\\
202	10.0365770778851\\
203	12.3548602661412\\
204	11.8848657085194\\
205	11.3356637199979\\
206	10.6410880202387\\
207	11.6233065893766\\
208	11.0303862520305\\
209	10.5420313008119\\
210	10.6007850791888\\
211	13.2114099743512\\
212	12.4793507651585\\
213	12.1904757782435\\
214	10.8312254577907\\
215	10.6338391233241\\
216	11.9095449888192\\
217	11.8159215187422\\
218	10.9342415314276\\
219	10.4150922402622\\
220	10.3062785374769\\
221	11.2824072685656\\
222	11.769239877493\\
223	11.4254266434637\\
224	10.7467012715854\\
225	11.0968251888378\\
226	12.2422384206751\\
227	11.806353142897\\
228	13.187658741144\\
229	12.4319491372266\\
230	12.2878743750914\\
231	11.4257755540618\\
232	12.5418891444681\\
233	12.2063731791442\\
234	12.8240971891971\\
235	11.6242419722901\\
236	12.3649281400185\\
237	12.1939123698777\\
238	12.1736958994351\\
239	11.8019245771022\\
240	11.6789617628692\\
241	11.512807065184\\
242	13.1993053002157\\
243	10.4612704904692\\
244	9.82614542564353\\
245	10.8854477598122\\
246	11.895847660845\\
247	11.3816450715894\\
248	10.8666064631925\\
249	10.8455553648949\\
250	10.5092376928265\\
251	9.52656160056512\\
252	9.33660206929456\\
253	9.3653117757892\\
254	9.31542178753206\\
255	9.12534717516398\\
256	8.53365604972396\\
257	8.47803368071528\\
258	8.49310205913591\\
259	8.17739015985681\\
260	8.16041204386284\\
261	7.62309648704231\\
262	9.53589451087447\\
263	9.46120746463491\\
264	9.69486613992983\\
265	10.8913879183062\\
266	11.5486141966046\\
267	10.9390487973092\\
268	11.7713484676976\\
269	11.0201477168114\\
270	9.80936654928853\\
271	10.9073701757143\\
272	10.278398272594\\
273	10.1173518163034\\
274	10.0208546391443\\
275	10.3219444144085\\
276	11.1820699865753\\
277	11.6034826580283\\
278	12.3904196392076\\
279	12.696693055407\\
280	12.4752197225909\\
281	11.8277552864721\\
282	10.4419192469219\\
283	9.87312249890009\\
284	9.68811821346896\\
285	10.1615044008603\\
286	9.99297732078509\\
287	10.2746153605931\\
288	10.5402210339519\\
289	10.8033446329647\\
290	10.6362595303229\\
291	12.2688924700337\\
292	10.9285046064228\\
293	11.5507112947981\\
294	11.2541677368779\\
295	11.1043806083141\\
296	10.8314516092092\\
297	10.532461978252\\
298	10.5316473184423\\
299	10.3682048426984\\
300	10.374371313981\\
301	10.5031879063599\\
302	9.37584675244492\\
303	9.05182123060943\\
304	9.04282821958735\\
305	9.02970080419568\\
306	10.0073354459841\\
307	9.20608903658806\\
308	12.1522309759405\\
309	10.8654431337885\\
310	10.605899596541\\
311	10.3824825850182\\
312	10.0953230776863\\
313	8.70302961655844\\
314	8.87214292460893\\
315	9.30368242213028\\
316	9.01540323806886\\
317	9.61513798569595\\
318	10.1987399041796\\
319	10.1533616155191\\
320	10.1384913054648\\
321	9.91357016753858\\
322	10.0952228798049\\
323	10.5530872430775\\
324	10.5418073029693\\
325	10.4107117412537\\
326	9.48409871960696\\
327	9.06843244467868\\
328	8.99331989536759\\
329	8.51362126216385\\
330	8.52565310917259\\
331	8.56929186292182\\
332	8.52335904808766\\
333	10.1845304999305\\
334	9.37463170818888\\
335	9.62214479097455\\
336	10.0368522836882\\
337	9.8055086282392\\
338	9.54187554227975\\
339	9.6311090135303\\
340	9.39905849049427\\
341	10.0779930899745\\
342	11.3724468152066\\
343	10.153642840025\\
344	9.90616908735186\\
345	10.0325084593528\\
346	9.82234918120398\\
347	9.72827610818622\\
348	9.5178795190048\\
349	9.52529390283739\\
350	9.09253013072845\\
351	8.80893039946409\\
352	9.212227877232\\
353	9.31304780821977\\
354	9.13586079357024\\
355	8.57379004842121\\
356	8.50192914735903\\
357	8.28057470587115\\
358	8.68350255380076\\
359	8.21320042010664\\
360	9.01646127779101\\
361	8.47977462327502\\
362	8.29085327410003\\
363	10.4418490134678\\
364	9.71302352767361\\
365	10.6730820804183\\
366	10.118357525404\\
367	9.85493642558288\\
368	8.67836245980951\\
369	11.0009971142715\\
370	12.188006110455\\
371	11.471013745623\\
372	10.8557668654826\\
373	9.49187532866813\\
374	9.54352721144934\\
375	9.32642822986803\\
376	9.63750049400161\\
377	9.46934126698882\\
378	9.02448194824271\\
379	8.7777879925998\\
380	8.56539599858747\\
381	8.39646872341568\\
382	8.62443243843157\\
383	9.03715017913934\\
384	9.73411453044026\\
385	10.1616716304227\\
386	10.008441338636\\
387	9.75459342537472\\
388	11.8362327793453\\
389	10.8346542895118\\
390	11.3779198123808\\
391	11.0338856240868\\
392	11.0584391715184\\
393	11.3074335724068\\
394	11.9099440671511\\
395	11.1613154488547\\
396	11.017129101347\\
397	10.1596010979705\\
398	10.8161062496027\\
399	11.0696653409482\\
400	10.4621380018951\\
401	10.0369363378584\\
402	9.71261242124618\\
403	9.24980501483212\\
404	9.22361693244581\\
405	8.37736177436343\\
406	9.28726723451563\\
407	9.16512141295118\\
408	8.55475234448232\\
409	9.21908964844317\\
410	8.2879185668149\\
411	8.25492905302264\\
412	7.60575760161352\\
413	8.11114635496104\\
414	9.48090676366316\\
415	9.59300986175577\\
416	9.53090640349086\\
417	9.38086459438099\\
418	9.1699486482179\\
419	10.4742950702808\\
420	10.3464903150074\\
421	10.7393717846365\\
422	10.4859298094804\\
423	9.46858858693761\\
424	10.6793162604856\\
425	10.4193594767465\\
426	11.1190741692563\\
427	10.5413260476208\\
428	10.4691410652625\\
429	10.350064652127\\
430	10.3906152968533\\
431	10.2439387853165\\
432	11.4501891255805\\
433	10.9367441200617\\
434	10.1172192571753\\
435	9.45725736887537\\
436	9.18449264194317\\
437	8.98919468895873\\
438	9.05679906707948\\
439	8.83927370018898\\
440	10.6735397782823\\
441	10.0206882533883\\
442	10.2046012800819\\
443	10.6416655103724\\
444	10.612160745209\\
445	10.7050528325525\\
446	10.9590867483008\\
447	11.4510198298239\\
448	9.81810338657197\\
449	9.18540168867145\\
450	8.93041986552381\\
451	8.86779852071575\\
452	8.922974826676\\
453	8.98227037181079\\
454	9.17528727227637\\
455	9.02999527993915\\
456	8.83567165572925\\
457	8.66491897001145\\
458	8.21563532637912\\
459	8.16227725281836\\
460	8.21303747071971\\
461	10.7523568961865\\
462	11.2082020495269\\
463	10.8503157035604\\
464	8.88222808026774\\
465	8.97452775523909\\
466	8.81162469984912\\
467	8.51416794315606\\
468	8.44636841998812\\
469	9.39610519058878\\
470	9.24695038851045\\
471	9.87739400046134\\
472	10.4768677678426\\
473	10.2888886189841\\
474	9.7761151213017\\
475	9.64060068593554\\
476	9.29157463445626\\
477	8.66498974188433\\
478	8.62530710943315\\
479	8.42611445684398\\
480	8.71389365453966\\
481	8.69308915752761\\
482	8.2995857563874\\
483	8.2108362954437\\
484	8.59781868578061\\
485	8.9513429481102\\
486	9.31261412605715\\
487	10.1291544717106\\
488	11.5963163563552\\
489	11.2638590326391\\
490	10.6621802743923\\
491	10.5783453553315\\
492	10.4512402273636\\
493	11.1835613705615\\
494	10.5754251136997\\
495	10.569403344422\\
496	10.5083252565174\\
497	10.7771123207211\\
498	10.9646439373413\\
499	10.5240240612481\\
500	10.3523127562421\\
501	10.391205543961\\
502	10.3661427323281\\
503	10.2245655896175\\
504	8.45650349943863\\
505	8.66831413817425\\
506	10.9649564293962\\
507	10.8759098935154\\
508	10.8916297671562\\
509	10.6841737077604\\
510	10.2051188950006\\
511	10.2260261835986\\
512	10.4731097767494\\
513	10.4338091626577\\
514	10.4770386597824\\
515	10.391462382585\\
516	10.83459224158\\
517	8.52638521929966\\
518	10.4803092051623\\
519	10.5950961151018\\
520	9.8959787722152\\
521	9.74477959048033\\
522	10.0337645018221\\
523	10.1313878853691\\
524	11.3889700801994\\
525	10.3771947633736\\
526	10.4851766012965\\
527	10.135020956074\\
528	9.86647057538781\\
529	10.6515162673337\\
530	11.3149543779838\\
531	11.2860594825036\\
532	10.3100002989119\\
533	10.0137125840343\\
534	9.9302702886949\\
535	11.2315023472408\\
536	10.8307843823458\\
537	10.2790787389939\\
538	9.36621298519043\\
539	9.38392646256929\\
540	9.39112043157998\\
541	9.35534066725542\\
542	10.9043538799656\\
543	10.8781218654827\\
544	10.8711733509762\\
545	10.0503978643325\\
546	9.18218376133416\\
547	9.26783655813199\\
548	9.11627702869946\\
549	9.05864406589282\\
550	8.95785178533098\\
551	9.2240539753654\\
552	9.14468092135474\\
553	9.74244977313284\\
554	9.1239985296692\\
555	9.24216855447842\\
556	8.25804465754137\\
557	8.06941536452578\\
558	7.88141063639402\\
559	8.76291049840537\\
560	10.0269676422564\\
561	10.8523988676888\\
562	10.6837937828491\\
563	10.4540832693698\\
564	10.1758057703126\\
565	12.1842351592832\\
566	11.4984928512847\\
567	10.5924580858603\\
568	10.2417874386939\\
569	10.2779247360825\\
570	10.1772832829822\\
571	10.9784625405975\\
572	10.5291974782133\\
573	10.1633376915325\\
574	9.20172277529426\\
575	8.95714394835323\\
576	10.4119461490601\\
577	10.1616651971481\\
578	12.1284208209155\\
579	11.5344495995904\\
580	11.2140981753657\\
581	11.0930709898574\\
582	10.6471741426957\\
583	10.5128096161322\\
584	9.12406408587405\\
585	9.3083291802454\\
586	9.03657728641893\\
587	8.87512421602003\\
588	8.97810075266423\\
589	9.06590136021509\\
590	9.77215289917005\\
591	9.7936227784509\\
592	11.0279717735793\\
593	9.96594180159511\\
594	10.8046998530639\\
595	11.1133735342722\\
596	10.8758985886768\\
597	11.54680813345\\
598	11.330270130143\\
599	11.6009886685761\\
600	10.9118781605227\\
601	12.9171566070354\\
602	12.3923649033377\\
603	11.3181693444524\\
604	10.8777504934345\\
605	10.7519623895259\\
606	10.5973448855811\\
607	10.6885300101103\\
608	10.6471693089265\\
609	10.126421723883\\
610	10.1258747460966\\
611	9.80748615347734\\
612	9.33767219887913\\
613	9.14219856894907\\
614	8.89171769924585\\
615	9.22360198057205\\
616	8.19538621636868\\
617	8.12720052322001\\
618	8.70441953478837\\
619	8.59489518059074\\
620	9.67034678937078\\
621	9.1690207718163\\
622	9.26299690603501\\
623	9.67906193778326\\
624	9.44303625184207\\
625	9.41683517881493\\
626	9.35945633854992\\
627	9.1819683366775\\
628	9.1239847401568\\
629	9.24795290507113\\
630	9.3082207106224\\
631	9.73869156074999\\
632	9.66907765673691\\
633	9.67513532429258\\
634	9.67854452301927\\
635	9.11950008070525\\
636	10.7511853817935\\
637	10.5256178014755\\
638	10.093876771025\\
639	10.0084103676572\\
640	9.40516878635028\\
641	9.37277221969815\\
642	9.13645666153068\\
643	8.72615681066584\\
644	8.76841764460091\\
645	8.82002086749519\\
646	9.60073163492604\\
647	9.45499808780034\\
648	9.9347505554651\\
649	9.7479443593757\\
650	9.374101942423\\
651	9.87009216484267\\
652	10.602675843643\\
653	10.3623636674309\\
654	9.85764985944887\\
655	10.1001242125914\\
656	9.9219623925856\\
657	9.62526765733688\\
658	9.52747420610319\\
659	9.51711639193142\\
660	9.29705438270891\\
661	9.33928502471213\\
662	9.39540972528875\\
663	9.93351198111609\\
664	9.88195322242677\\
665	10.1656994084531\\
666	9.93195918793197\\
667	9.87214287159411\\
668	9.62505182014895\\
669	9.16803167018361\\
670	8.61027301707693\\
671	8.33060876603996\\
672	8.37653907772645\\
673	8.76541086411926\\
674	9.42400977977919\\
675	9.09462036891249\\
676	8.97815079261422\\
677	9.00013526407755\\
678	9.30793021519616\\
679	9.73638661813627\\
680	9.2307296053515\\
681	9.69219714690126\\
682	9.01052654133799\\
683	9.5269839148257\\
684	9.35879100225527\\
685	9.09205037303514\\
686	9.20378136244719\\
687	9.538883733504\\
688	9.22317856490464\\
689	10.0204122172784\\
690	9.35567926450168\\
691	9.02738637752718\\
692	9.4380764186047\\
693	9.51879049120022\\
694	9.12587413433221\\
695	9.79832501008083\\
696	9.48388220890254\\
697	9.39437033919509\\
698	9.08849917220002\\
699	8.40748938753871\\
700	10.3927891752628\\
701	10.2733365961666\\
702	10.7886844684012\\
703	10.5891953091034\\
704	10.4248230110293\\
705	11.0232656531551\\
706	12.0524141383339\\
707	11.3814508635962\\
708	11.4415802457654\\
709	10.5395783875888\\
710	10.0296014771605\\
711	10.0641708677705\\
712	10.0381463970929\\
713	10.2099223711581\\
714	9.85536587443726\\
715	9.12137731297354\\
716	9.33384085509452\\
717	9.22611819560666\\
718	10.246234152868\\
719	8.9720455596585\\
720	9.90722764824249\\
721	9.56427852228913\\
722	9.70703870617257\\
723	9.13792010630676\\
724	9.03742873759447\\
725	8.67004736421719\\
726	8.71094203212816\\
727	9.25583525791679\\
728	9.58003308618791\\
729	8.89655371806476\\
730	8.8467029697244\\
731	8.88858769071174\\
732	8.77248960516717\\
733	9.00573666162405\\
734	9.11383451954619\\
735	9.00842125556479\\
736	10.407997077099\\
737	10.2106751615556\\
738	9.4708804105379\\
739	9.88633461098743\\
740	9.59882033218206\\
741	9.85825495559852\\
742	9.0759310022193\\
743	9.36064777065742\\
744	9.26724132245473\\
745	9.0588675203457\\
746	8.97604647710702\\
747	9.96882943993927\\
748	9.60603082401225\\
749	10.1638545300339\\
750	9.97259146559754\\
751	9.85427227059176\\
752	9.19239991361313\\
753	8.9382651653893\\
754	9.49228068463346\\
755	9.68592345005582\\
756	9.56354299651544\\
757	8.96588313987751\\
758	9.03724601617955\\
759	8.95590783836905\\
760	8.73219822662209\\
761	8.61626870825982\\
762	8.88477144862603\\
763	9.53344698701276\\
764	9.19660459759791\\
765	9.18055754938436\\
766	9.07601622399618\\
767	9.98544039814835\\
768	10.2886535346235\\
769	10.4135067114459\\
770	9.64080490710653\\
771	9.63139727475752\\
772	9.00194282449561\\
773	8.65106252223286\\
774	8.59189425252657\\
775	8.63697696096077\\
776	8.60367833948406\\
777	8.18412058694492\\
778	8.79305840421272\\
779	8.17758985931551\\
780	8.93626822695293\\
781	8.61692903484992\\
782	8.57047330821959\\
783	8.01912808502229\\
784	7.76190230647398\\
785	9.6083362637111\\
786	8.96433982562265\\
787	8.71081124263458\\
788	7.83088446322973\\
789	7.96998674052522\\
790	8.17368858054016\\
791	8.22851115220263\\
792	8.2481547756951\\
793	8.83945886383539\\
794	8.31962010824421\\
795	8.22413785609269\\
796	9.31904998022568\\
797	8.85248981413542\\
798	9.47705074227304\\
799	9.39341149277713\\
800	9.32684206380472\\
801	9.06953434505915\\
802	8.6857460262199\\
803	9.1672874359359\\
804	8.76870123110614\\
805	8.59808609341456\\
806	8.44605316310412\\
807	8.50515584818658\\
808	8.55934189056885\\
809	8.14511017924863\\
810	8.14381474083878\\
811	8.09085031943996\\
812	10.2650990561741\\
813	10.5767838907283\\
814	10.4448684625382\\
815	10.4617285362081\\
816	9.57696236817926\\
817	8.76483102975506\\
818	9.02297592984012\\
819	8.90824357963521\\
820	8.7875084997448\\
821	8.74113069020653\\
822	8.74344617971478\\
823	9.36632955005324\\
824	9.10632890770471\\
825	9.05258947826871\\
826	8.96069846360423\\
827	8.59493838586761\\
828	8.58768681171558\\
829	9.12503254696604\\
830	8.94542606453193\\
831	8.95753717518797\\
832	8.57769575855317\\
833	7.65108558633342\\
834	8.98561921654745\\
835	8.5996606318365\\
836	8.5241936443371\\
837	8.70259224569067\\
838	8.18314751174027\\
839	8.16394693446413\\
840	7.86880250238048\\
841	8.00279130211968\\
842	7.57559699383883\\
843	8.00211126080268\\
844	7.91883817452332\\
845	10.517332757104\\
846	10.6394199101769\\
847	11.0727510157491\\
848	10.7181390412344\\
849	10.3993282297246\\
850	10.2122956213907\\
851	10.8794652860647\\
852	9.94291731641168\\
853	9.76991805849551\\
854	9.42370974517669\\
855	8.86312526217872\\
856	8.84859652299826\\
857	9.30242725139984\\
858	9.24221381789395\\
859	9.1533320450019\\
860	8.601984722893\\
861	8.38971530563033\\
862	8.45349724649433\\
863	8.52183843147846\\
864	8.52912871552476\\
865	9.81719761652009\\
866	9.7739345083724\\
867	9.41670788103802\\
868	9.28988962399579\\
869	9.05719067510547\\
870	10.0956976726078\\
871	9.18594406171764\\
872	9.15264640183376\\
873	9.03761586746341\\
874	9.4443857128571\\
875	8.77859152211138\\
876	8.66110306898183\\
877	8.35085758400223\\
878	8.47493989794196\\
879	9.25061229528698\\
880	10.1597110484357\\
881	9.83594672463735\\
882	10.0686841003155\\
883	10.7917830517964\\
884	10.9819833035999\\
885	9.98802108401934\\
886	9.84541702938057\\
887	9.68604167650756\\
888	9.60174388280145\\
889	11.0852573042868\\
890	10.1506781234197\\
891	10.0216059131108\\
892	9.48855195725896\\
893	9.81358902501724\\
894	9.87883014826242\\
895	9.38442932175652\\
896	9.73675459949343\\
897	9.50044652793068\\
898	9.57833570302353\\
899	10.6857524774594\\
900	10.7020600592534\\
901	9.61508892346288\\
902	8.7361911797027\\
903	8.74812698452736\\
904	8.98915977903443\\
905	8.68948563351507\\
906	7.92192769968298\\
907	8.0609887431764\\
908	8.11710319363243\\
909	7.81725542957644\\
910	7.73694017907601\\
911	7.67554888303298\\
912	7.95629533484296\\
913	9.57992476555921\\
914	8.7259028407417\\
915	10.1171424461769\\
916	9.79257975141051\\
917	10.1226681183785\\
918	9.64722975949596\\
919	9.41395730794848\\
920	9.37780502957906\\
921	9.21008516453123\\
922	10.1300125978562\\
923	9.39226172147735\\
924	9.08353865301239\\
925	9.17411500128395\\
926	9.03042397263964\\
927	8.93943501155222\\
928	11.1667012084704\\
929	10.8515327463137\\
930	10.8180312924538\\
931	12.2300453857574\\
932	11.5035554191572\\
933	10.849900380418\\
934	11.2271605194777\\
935	9.56678441040591\\
936	9.55402153184496\\
937	9.56888086183466\\
938	9.42827681035796\\
939	8.66539589846828\\
940	8.95502262886491\\
941	8.41731363370052\\
942	8.74282623405358\\
943	8.48499991443881\\
944	8.52876101707443\\
945	8.89685546409681\\
946	8.55214764800692\\
947	9.35597274580686\\
948	9.16815596039516\\
949	11.2545031101427\\
950	12.2866503305186\\
951	11.4428427351631\\
952	10.3708937763272\\
953	9.34814992668861\\
954	9.22607195269732\\
955	8.94008267427605\\
956	8.55310260193456\\
957	8.94060606733064\\
958	8.98363133145883\\
959	9.01040566270926\\
960	8.95438392914114\\
961	9.20040432977674\\
962	8.98176740382775\\
963	9.06402234855062\\
964	8.88002625078435\\
965	8.79424985025568\\
966	8.7837902191268\\
967	8.62912891142277\\
968	8.62086095220154\\
969	9.53719067257377\\
970	9.53762477075873\\
971	10.0482898213986\\
972	10.5878246862839\\
973	10.3697044267388\\
974	9.76028066761118\\
975	9.18646169332401\\
976	8.92676376868969\\
977	9.31327022676387\\
978	9.22803482363609\\
979	9.15151064132058\\
980	8.797372748655\\
981	8.8473378536495\\
982	9.70121353578851\\
983	9.20178792842107\\
984	8.52505580129175\\
985	8.56249659622754\\
986	8.52886365430035\\
987	8.82145726273649\\
988	8.48683216232726\\
989	8.41437096608756\\
990	10.7231839591148\\
991	10.4486511611671\\
992	10.7452069194365\\
993	9.5197924327458\\
994	9.11456635690929\\
995	9.04378409065327\\
996	8.87901047691566\\
997	8.84016027928509\\
998	8.79650857718075\\
999	8.58877542914252\\
1000	8.56692555671208\\
1001	8.7470687064817\\
1002	8.84158246127182\\
1003	8.61038388733906\\
1004	8.84835483214711\\
1005	8.58531219888343\\
1006	8.49169676214753\\
1007	8.33340035479407\\
1008	8.43766142543408\\
1009	8.26985717392586\\
1010	8.49545180679698\\
1011	8.60410307589001\\
1012	8.68007476588211\\
1013	8.37190415372969\\
1014	8.59589262210795\\
1015	8.1952884908511\\
1016	8.2203209893107\\
1017	8.40318906254739\\
1018	8.82323905468691\\
1019	8.72118125120216\\
1020	8.79220631549522\\
1021	10.6238803770395\\
1022	10.0434083952439\\
1023	9.89448093466649\\
1024	10.0547380669608\\
1025	9.67363881251712\\
1026	10.1021730141398\\
1027	10.0220613336813\\
1028	10.0269290528949\\
1029	9.6235547906113\\
1030	9.71124346520454\\
1031	9.70976687381424\\
1032	9.11411171648257\\
1033	8.99580848900799\\
1034	10.2570472237981\\
1035	9.41785248130903\\
1036	9.44543609988481\\
1037	8.98248695756066\\
1038	8.9428898827407\\
1039	9.03537270188876\\
1040	9.06422467732793\\
1041	8.75337592999482\\
1042	9.33572266121593\\
1043	9.25003423393904\\
1044	10.5815479350951\\
1045	10.4898805503235\\
1046	10.250309053032\\
1047	9.32240825746984\\
1048	9.21879365835146\\
1049	9.83793190630946\\
1050	8.90001458954648\\
1051	8.69580949621737\\
1052	9.27105256574771\\
1053	9.93339528491611\\
1054	9.22626422412674\\
1055	9.64147394081674\\
1056	10.2068329187347\\
1057	9.89668487511111\\
1058	9.89888329413564\\
1059	9.79798693193671\\
1060	9.35200941924261\\
1061	9.03633350381471\\
1062	9.47978700058517\\
1063	9.35043964904542\\
1064	9.22285184765365\\
1065	9.69510872572724\\
1066	9.35740075863934\\
1067	9.52771544843059\\
1068	9.48524095148822\\
1069	9.3071793150127\\
1070	9.34106599070142\\
1071	9.31904938608908\\
1072	8.90511093397076\\
1073	8.96346429473709\\
1074	9.3347696461897\\
1075	9.82165362844156\\
1076	9.6606877501622\\
1077	8.7268114721926\\
1078	9.42613696523671\\
1079	9.15044643886453\\
1080	9.06770908903706\\
1081	8.93464142150627\\
1082	8.94880305783344\\
1083	9.25021815519875\\
1084	9.21816154788822\\
1085	9.21125893223823\\
1086	9.20769404880818\\
1087	9.14004574842908\\
1088	8.75867754529589\\
1089	8.6530893255275\\
1090	8.50623032858267\\
1091	7.78034697941745\\
1092	7.58107156039797\\
1093	7.87136856968219\\
1094	8.84184144471537\\
1095	8.77069718659531\\
1096	8.80175920484949\\
1097	9.50192892450538\\
1098	9.46748215223482\\
1099	9.75800439265274\\
1100	9.72118105212091\\
1101	9.04275597040529\\
1102	9.19957456265368\\
1103	8.79515801570849\\
1104	9.23438755841621\\
1105	9.57169357828448\\
1106	9.73217268103279\\
1107	9.5594195268205\\
1108	9.39979037793761\\
1109	8.88733316548493\\
1110	8.77310190961812\\
1111	9.14579769069415\\
1112	8.88727407650548\\
1113	9.4039681480687\\
1114	8.99666863491967\\
1115	8.6846792118006\\
1116	8.75367870433798\\
1117	8.81999714619991\\
1118	11.2857599903271\\
1119	10.0686353490707\\
1120	9.23957618874989\\
1121	10.486821950852\\
1122	10.2936441751853\\
1123	10.2869458909967\\
1124	10.0379426053237\\
1125	10.2961961194124\\
1126	9.79352753720294\\
1127	10.209264084778\\
1128	10.2747940424305\\
1129	10.1462382307827\\
1130	10.5357918898368\\
1131	9.82740594851904\\
1132	9.90071742774089\\
1133	10.3402549534063\\
1134	9.99279520919731\\
1135	9.70534990982467\\
1136	9.62282865050387\\
1137	9.32296208606154\\
1138	8.86992113882368\\
1139	8.83808506750829\\
1140	8.96614172098988\\
1141	8.90766130608158\\
1142	8.80761226036956\\
1143	8.85854869292886\\
1144	8.17896894100009\\
1145	8.29071262688237\\
1146	8.69458709491361\\
1147	10.2624108708048\\
1148	9.5488739902894\\
1149	9.38059193760463\\
1150	9.42908434384794\\
1151	9.61519694583957\\
1152	9.26727699383909\\
1153	9.2632869387566\\
1154	9.51356913544689\\
1155	10.4186319975192\\
1156	9.87462986998283\\
1157	9.78356335475365\\
1158	9.97143514596914\\
1159	9.56163373794368\\
1160	9.56467675431246\\
1161	10.3900588347939\\
1162	10.0968569512706\\
1163	9.96525052838488\\
1164	9.32555141990324\\
1165	9.31697891981506\\
1166	9.22740693494248\\
1167	9.39521669670545\\
1168	9.19427178924733\\
1169	8.95438024870192\\
1170	9.16728022724167\\
1171	9.05173806661637\\
1172	8.82007082978494\\
1173	8.31485864449582\\
1174	8.15846569154344\\
1175	8.65566688566757\\
1176	8.4181675765949\\
1177	8.96881974404415\\
1178	9.62679119599635\\
1179	9.59759436575199\\
1180	9.45195217652683\\
1181	10.3055462901706\\
1182	9.83732536201799\\
1183	10.0955435335429\\
1184	9.97387493292111\\
1185	9.90459303170157\\
1186	9.90849691376802\\
1187	9.4687830100162\\
1188	9.33123157093449\\
1189	9.47183182902962\\
1190	9.67985255887477\\
1191	9.40173198452605\\
1192	10.2340905194569\\
1193	10.5779393893015\\
1194	10.6816738195519\\
1195	10.323849565568\\
1196	9.96378369569427\\
1197	11.0300902256425\\
1198	10.8124928126157\\
1199	10.1575890958308\\
1200	9.48999579474885\\
1201	9.43382147349681\\
1202	9.53411961361411\\
1203	9.70910332280114\\
1204	9.50478392719055\\
1205	9.59991257838\\
1206	9.60001289002908\\
1207	8.89839257497217\\
1208	9.09114858121629\\
1209	8.96900191159965\\
1210	8.68532732193603\\
1211	8.71288193878205\\
1212	8.71301724594974\\
1213	8.8202386329059\\
1214	9.10638642310071\\
1215	10.8633432377148\\
1216	10.3101165122929\\
1217	11.1321333675684\\
1218	10.457554797292\\
1219	10.2884143700122\\
1220	9.44463346989582\\
1221	9.49644992818134\\
1222	9.54694902059558\\
1223	9.50267213507064\\
1224	9.39605833155698\\
1225	9.37531614101956\\
1226	9.10784509808522\\
1227	9.14580795166483\\
1228	10.0889011767926\\
1229	10.0774585428775\\
1230	9.9416818249803\\
1231	9.72956623065499\\
1232	9.82315819684777\\
1233	9.37233694704834\\
1234	9.31803760392384\\
1235	9.31862604854413\\
1236	8.77706483313087\\
1237	8.62669499836636\\
1238	8.93490727434037\\
1239	9.18720546663452\\
1240	8.76289840444608\\
1241	8.61995367766156\\
1242	8.56868906129776\\
1243	9.4052689907358\\
1244	8.88352415023559\\
1245	8.64172459891697\\
1246	8.92045547330929\\
1247	8.8234033273167\\
1248	8.63794078194028\\
1249	8.67519440616489\\
1250	9.0548441018691\\
1251	8.8277706784798\\
1252	8.69861455666279\\
1253	8.52438803767544\\
1254	8.73027573279752\\
1255	8.96170931529928\\
1256	9.08349780661451\\
1257	8.97538161908466\\
1258	8.87631792477678\\
1259	8.45067754623006\\
1260	8.01521554669223\\
1261	8.10440347414543\\
1262	8.11072549116625\\
1263	8.06789245612463\\
1264	9.62339482560112\\
1265	9.42283220826246\\
1266	9.31039498487517\\
1267	8.77450263043512\\
1268	8.68462984249796\\
1269	8.72896684581582\\
1270	9.44169007213054\\
1271	9.44682632188476\\
1272	9.46655767271528\\
1273	10.5330805811794\\
1274	10.710551778633\\
1275	10.6973595549528\\
1276	9.55850333252436\\
1277	10.3084565129129\\
1278	10.0968002521791\\
1279	9.91656423572241\\
1280	9.85877389943193\\
1281	9.75878272577402\\
1282	9.60237542705042\\
1283	9.50293975874463\\
1284	9.41863735954071\\
1285	8.576867668261\\
1286	8.58291336871195\\
1287	8.54350894809301\\
1288	8.61761756293678\\
1289	8.70693065870651\\
1290	8.92507641879479\\
1291	9.18914146225609\\
1292	9.17548662250107\\
1293	8.86848263122818\\
1294	8.91777203951756\\
1295	8.95297013534744\\
1296	9.34337683251135\\
1297	9.22091638257683\\
1298	9.50092182919943\\
1299	9.66658862373923\\
1300	9.44771040571611\\
1301	9.23073537819408\\
1302	8.79869098165977\\
1303	8.47890863368866\\
1304	8.85861832593415\\
1305	8.85586546855636\\
1306	8.83574919341117\\
1307	9.04263108374789\\
1308	8.72482988085973\\
1309	8.80288497707319\\
1310	9.09391254846821\\
1311	8.96362761778032\\
1312	8.85169963183928\\
1313	8.78849044663034\\
1314	9.20733163505773\\
1315	9.34962363757705\\
1316	9.80199796294813\\
1317	9.6957472329257\\
1318	9.54754803710673\\
1319	9.29580007184817\\
1320	8.70969276198329\\
1321	8.52505919218398\\
1322	9.39208249473314\\
1323	9.46012222329491\\
1324	8.8770793447354\\
1325	8.88974180495277\\
1326	9.31362888949454\\
1327	9.69763215442462\\
1328	9.57943941759115\\
1329	9.66741207597174\\
1330	9.19068389182782\\
1331	9.54142238188643\\
1332	9.6275041749301\\
1333	9.65438923862952\\
1334	9.73685523308639\\
1335	10.7956060693401\\
1336	10.738077626933\\
1337	10.472776790633\\
1338	9.74624285609485\\
1339	9.42988966813643\\
1340	9.37325526902707\\
1341	9.3411970313527\\
1342	8.82765237345555\\
1343	9.06199972945681\\
1344	9.28172469027182\\
1345	9.32094730394931\\
1346	9.30636282440366\\
1347	9.02301906317078\\
1348	9.30836998441096\\
1349	9.4099805365713\\
1350	9.16658311787451\\
1351	10.0823844593827\\
1352	10.0458568228598\\
1353	9.57891050704952\\
1354	9.19629531313749\\
1355	9.04492980475149\\
1356	9.1617628839952\\
1357	9.02785049343522\\
1358	8.99727798494595\\
1359	9.04300429206745\\
1360	8.91526506405782\\
1361	8.93016833954164\\
1362	8.54959155967069\\
1363	8.56325907559376\\
1364	8.28966428096059\\
1365	8.20290847227113\\
1366	8.86367452844351\\
1367	8.9796436417937\\
1368	9.32818809961288\\
1369	8.94656365617695\\
1370	10.2058539084191\\
1371	9.6443775728966\\
1372	9.42080571923729\\
1373	9.27446718124919\\
1374	9.25922167171045\\
1375	9.01109676727301\\
1376	9.6461076934059\\
1377	9.39368793749295\\
1378	10.072334767401\\
1379	9.89270823308057\\
1380	9.85615385023415\\
1381	9.61482737273975\\
1382	9.52166950068188\\
1383	9.43391805488948\\
1384	9.08251023521267\\
1385	9.1624315780791\\
1386	9.09647405515108\\
1387	8.75022424623318\\
1388	9.27172676802462\\
1389	8.87906147955991\\
1390	8.76383076284015\\
1391	8.94508449814221\\
1392	9.15097158823331\\
1393	8.87371504063697\\
1394	9.0216682624396\\
1395	9.31819386840032\\
1396	9.17385068614325\\
1397	10.6933894175499\\
1398	10.5759340767534\\
1399	10.1833885609326\\
1400	9.94848595720724\\
1401	9.8494234978044\\
1402	9.64831243162009\\
1403	11.3590962713466\\
1404	11.0333417412286\\
1405	10.0105742620087\\
1406	9.78976155316252\\
1407	9.60634242581638\\
1408	9.50882859262936\\
1409	9.23895388491457\\
1410	9.37308072992492\\
1411	8.93912184784993\\
1412	8.96311424256967\\
1413	9.00848234463227\\
1414	9.08305575691778\\
1415	10.4601860359994\\
1416	10.4264856758566\\
1417	10.4161283461439\\
1418	10.2942429450284\\
1419	10.2952587626038\\
1420	10.0101624451553\\
1421	9.85579069188787\\
1422	9.79030767655856\\
1423	9.69857532947123\\
1424	9.3302728180044\\
1425	9.08841819988922\\
1426	8.99094170456375\\
1427	8.86135343007693\\
1428	9.0446761388774\\
1429	9.79475312972951\\
1430	9.86276065163751\\
1431	9.23334093751986\\
1432	9.12910214993626\\
1433	8.9479670641216\\
1434	8.86209806058559\\
1435	8.88003476350939\\
1436	11.391101773568\\
1437	11.297655797993\\
1438	10.7059987039499\\
1439	10.84986491954\\
1440	10.5016685796844\\
1441	10.1309808130635\\
1442	9.97160947748196\\
1443	10.5508135916124\\
1444	10.4256189721619\\
1445	9.74314667455155\\
1446	9.45409337065684\\
1447	9.42222915538146\\
1448	9.2690122519169\\
1449	9.46277130204224\\
1450	8.78402245471387\\
1451	8.68878198134906\\
1452	8.53386279680532\\
1453	8.78023515159295\\
1454	8.66116931808087\\
1455	8.70803235354649\\
1456	8.77413951051625\\
1457	8.53125413665325\\
1458	9.15616046601475\\
1459	9.08814176029428\\
1460	9.142132436374\\
1461	10.7068299763747\\
1462	10.3777580159043\\
1463	10.1409151979859\\
1464	10.4402237997232\\
1465	10.1244234981394\\
1466	9.9636951657567\\
1467	9.65976975645093\\
1468	9.5222821168315\\
1469	9.72848719951679\\
1470	10.5336803435977\\
1471	10.3582434047248\\
1472	10.1631775751113\\
1473	9.88194772124296\\
1474	9.11389850163437\\
1475	9.11331561792206\\
1476	9.10461527753703\\
1477	8.6757473258486\\
1478	8.72331015901035\\
1479	8.78421426477652\\
1480	8.91069206886049\\
1481	9.07609730661608\\
1482	9.23522812331687\\
1483	8.90595333151474\\
1484	9.2248539659381\\
1485	9.24900917828456\\
1486	9.29161755278821\\
1487	9.46855355498144\\
1488	8.85772924504972\\
1489	9.0652373329716\\
1490	8.92343272987319\\
1491	9.46406162686584\\
1492	9.3571256078773\\
1493	9.91928387950407\\
1494	9.80535322943937\\
1495	9.49971025024576\\
1496	9.65639927643703\\
1497	9.67384836104936\\
1498	9.68031143686574\\
1499	9.46980743262362\\
1500	9.48658118482108\\
};
\addlegendentry{$J_k$};

\addplot [color=mycolor2, only marks, mark=*, mark options={solid, mycolor2}, mark size = 1.0pt]
  table[row sep=crcr]{%
25	23.7969232587093\\
50	17.727511691114\\
75	13.4009785308729\\
100	12.4619729882598\\
125	11.5116793211982\\
150	11.1876275402757\\
175	10.731555471729\\
200	10.6879495003076\\
225	10.4322666353406\\
250	10.3313192459455\\
275	10.2437473933816\\
300	10.0673784164867\\
325	10.1057666996893\\
350	10.0367747254471\\
375	9.84954268249329\\
400	9.80138610455851\\
425	9.78841159327842\\
450	9.70692039741848\\
475	9.72805924887106\\
500	9.66633083743237\\
525	9.56176043862394\\
550	9.53440006256118\\
575	9.47051635126598\\
600	9.47993515121673\\
625	9.42854183020266\\
650	9.41788907472593\\
675	9.38635914243634\\
700	9.37906582711595\\
725	9.34690672175051\\
750	9.32869490511601\\
775	9.35163278687008\\
800	9.37654706906265\\
825	9.37067420824104\\
850	9.41845167426949\\
875	9.46150194475827\\
900	9.45784026565705\\
925	9.41964656979337\\
950	9.36018871368358\\
975	9.36077244593374\\
1000	9.38174824110681\\
1025	9.3500045099586\\
1050	9.35507322577734\\
1075	9.36768639355814\\
1100	9.38878235332638\\
1125	9.37472814181454\\
1150	9.34188780906363\\
1175	9.38053789665482\\
1200	9.32657071075654\\
1225	9.36398515520923\\
1250	9.3673125929624\\
1275	9.36666786428053\\
1300	9.39616353835622\\
1325	9.33034920706125\\
1350	9.33535580544807\\
1375	9.34319362872605\\
1400	9.36246813303317\\
1425	9.37775467701758\\
1450	9.38648758497832\\
1475	9.32977758069189\\
1500	9.33593662591971\\
};
\addlegendentry{$\mathbb{E}_\xi\big[c(\hat{\mx},\xi)\big]$};
\end{axis}
\end{tikzpicture}%

%% file: Graphics/Clamp/Clamp_quants.tex
%
%
\definecolor{mycolor1}{rgb}{0.00000,0.44700,0.74100}%
%
\begin{tikzpicture}

\begin{axis}[%
width=.9\linewidth,
height=.9\linewidth,
at={(0.758in,0.481in)},
xmin=0,
xmax=1500,
ymin=8,
ymax=15,
axis background/.style={fill=white},
xmajorgrids,
ymajorgrids,
xlabel = {Iteration}
]
\addplot [draw=none,color=mycolor1, forget plot,  name path=A]
  table[row sep=crcr]{%
1	51.06269673515\\
10	30.6528203164993\\
20	25.2127163519017\\
30	22.1366578404791\\
40	19.701398091614\\
50	16.9277829200518\\
60	14.325211894312\\
70	13.4022778756156\\
80	12.755471649848\\
90	12.1933252815881\\
100	11.8323517842594\\
110	11.3006999416125\\
120	11.1847703501922\\
130	11.0030635412099\\
140	10.9042559418028\\
150	10.8103777528979\\
160	10.5381594926851\\
170	10.5016250495042\\
180	10.4332658399735\\
190	10.4053986641275\\
200	10.3617664664881\\
210	10.1563097925106\\
220	10.165069474231\\
230	10.1206875309215\\
240	10.1298471240528\\
250	10.0759224027368\\
260	9.97039749397851\\
270	9.94439875509422\\
280	9.9275769761202\\
290	9.93070928119543\\
300	9.89175868460669\\
310	9.81144381970616\\
320	9.80816730462667\\
330	9.81510803646847\\
340	9.78683410317873\\
350	9.80738501630353\\
360	9.69405200574068\\
370	9.68136834707459\\
380	9.68094168909795\\
390	9.65762105126244\\
400	9.65246647048228\\
410	9.57260810533723\\
420	9.57429373522085\\
430	9.5766722255648\\
440	9.57434237131047\\
450	9.57985276299397\\
460	9.51992439869037\\
470	9.52399195567397\\
480	9.51189657818862\\
490	9.52271460870332\\
500	9.5163178734316\\
510	9.46245141177439\\
520	9.45659504134194\\
530	9.46382899533102\\
540	9.46040236772844\\
550	9.44645635227692\\
560	9.40400376849375\\
570	9.41019606705247\\
580	9.41380372702164\\
590	9.40748145125871\\
600	9.40810013388214\\
610	9.38707570292173\\
620	9.37555522696867\\
630	9.36844530663896\\
640	9.37753493363811\\
650	9.3717629381769\\
660	9.35220949346306\\
670	9.34828424808385\\
680	9.34672048849021\\
690	9.35416719558976\\
700	9.35721098279347\\
710	9.32153177069287\\
720	9.32385261125029\\
730	9.32373479465648\\
740	9.32249681379893\\
750	9.32265524015044\\
760	9.32193325141\\
770	9.32241048097216\\
780	9.31405223503953\\
790	9.31217632185479\\
800	9.31108758372969\\
810	9.31816550016354\\
820	9.31409095957642\\
830	9.31041666182811\\
840	9.29855112262787\\
850	9.30964953373906\\
860	9.30840922674601\\
870	9.30923455593731\\
880	9.30280310793034\\
890	9.30988589215587\\
900	9.30887419539935\\
910	9.30104380007678\\
920	9.31731237681674\\
930	9.29712931531877\\
940	9.30806636493364\\
950	9.3127660978569\\
960	9.31459741766912\\
970	9.30786132657781\\
980	9.3116806744567\\
990	9.2973337403147\\
1000	9.29387443864992\\
1010	9.30225451522333\\
1020	9.29263792105195\\
1030	9.29727192568632\\
1040	9.29363735989577\\
1050	9.29210771849803\\
1060	9.28601370917646\\
1070	9.28984751999865\\
1080	9.28314158069079\\
1090	9.28515172683321\\
1100	9.28122341431909\\
1110	9.28669967315663\\
1120	9.28615858608719\\
1130	9.29672450834784\\
1140	9.28351537371089\\
1150	9.28463888340355\\
1160	9.29188649681979\\
1170	9.29053984786669\\
1180	9.28270231774562\\
1190	9.28251571046624\\
1200	9.28784201184465\\
1210	9.28609621752688\\
1220	9.29632965260056\\
1230	9.28048340845218\\
1240	9.28779473641943\\
1250	9.2756871644171\\
1260	9.27580377868971\\
1270	9.28631903091619\\
1280	9.28684821715664\\
1290	9.28204026086358\\
1300	9.28292438272165\\
1310	9.28920302633645\\
1320	9.29067030569365\\
1330	9.27216031796159\\
1340	9.28212376824699\\
1350	9.27297443673749\\
1360	9.27560899523574\\
1370	9.27597319966009\\
1380	9.27662149063126\\
1390	9.2843050218096\\
1400	9.28203085467473\\
1410	9.28904910105431\\
1420	9.27664471390024\\
1430	9.27897332829867\\
1440	9.28609964785877\\
1450	9.28586997113569\\
1460	9.28693365861289\\
1470	9.28175181630806\\
1480	9.277307445104\\
1490	9.28185864240053\\
1500	9.28959944542206\\
};
\addplot [draw=none,color=mycolor1, forget plot,  name path=B]
  table[row sep=crcr]{%
1	51.06269673515\\
10	30.7015191440484\\
20	25.3252525707104\\
30	22.4673668694614\\
40	19.9929955473844\\
50	17.4217541773402\\
60	14.8693620397674\\
70	13.8168798990172\\
80	13.0232156322386\\
90	12.4141795241281\\
100	12.0640052973815\\
110	11.5599541788989\\
120	11.3275916377061\\
130	11.1718958500032\\
140	11.0333657942076\\
150	10.9464907610522\\
160	10.6714332361214\\
170	10.609614837703\\
180	10.56785830761\\
190	10.5181527283072\\
200	10.4703279785942\\
210	10.2958496447663\\
220	10.2666749688938\\
230	10.2562781456833\\
240	10.2228349560439\\
250	10.2051454444874\\
260	10.0680652719397\\
270	10.0395813514779\\
280	10.0176266000014\\
290	10.0014914089981\\
300	9.98532550387869\\
310	9.89618926894473\\
320	9.88500697100158\\
330	9.88359704622281\\
340	9.86842256548103\\
350	9.86632574549083\\
360	9.77308890630517\\
370	9.77500162878609\\
380	9.76594234911707\\
390	9.74484121426832\\
400	9.73845924158046\\
410	9.66486869381661\\
420	9.66300054409534\\
430	9.66624236553728\\
440	9.66509960211631\\
450	9.66523322375027\\
460	9.59674280686597\\
470	9.59578598176501\\
480	9.58756258829599\\
490	9.58103325541996\\
500	9.59159061825475\\
510	9.53834325087749\\
520	9.52632654543211\\
530	9.5202929802058\\
540	9.52163173972658\\
550	9.51435204795556\\
560	9.47625789267941\\
570	9.47788892257479\\
580	9.47339928174842\\
590	9.47019938835409\\
600	9.47760058231217\\
610	9.43988694317477\\
620	9.43490958460811\\
630	9.43870047687163\\
640	9.44391064910522\\
650	9.44238439377939\\
660	9.4031237276779\\
670	9.40586721726912\\
680	9.40240619611273\\
690	9.40000135623103\\
700	9.40622791173912\\
710	9.38037645362114\\
720	9.38345179341598\\
730	9.37532638294013\\
740	9.37027820829628\\
750	9.36970538757904\\
760	9.36984111787984\\
770	9.36983260174662\\
780	9.37466950869973\\
790	9.37855283370816\\
800	9.36674516728647\\
810	9.36799705101678\\
820	9.3576415009225\\
830	9.35462566665083\\
840	9.35855442742662\\
850	9.35318124183623\\
860	9.35893962564864\\
870	9.35753253756979\\
880	9.35688845237672\\
890	9.35556473770353\\
900	9.369442231429\\
910	9.36339679415254\\
920	9.35807554193864\\
930	9.36095360744856\\
940	9.36204473176495\\
950	9.37194770343374\\
960	9.37120359487519\\
970	9.35263101187423\\
980	9.34874458001223\\
990	9.36200538169483\\
1000	9.35146279307175\\
1010	9.36136644071463\\
1020	9.35881353451827\\
1030	9.36201097554084\\
1040	9.35744708001353\\
1050	9.35245336302748\\
1060	9.3564556482891\\
1070	9.35629540958419\\
1080	9.35296712133375\\
1090	9.35875072727581\\
1100	9.34892400147105\\
1110	9.34918941123713\\
1120	9.35056788837655\\
1130	9.34780156727742\\
1140	9.34271907819805\\
1150	9.33890953680757\\
1160	9.33781339595842\\
1170	9.33988245178438\\
1180	9.33991926474102\\
1190	9.34917433214777\\
1200	9.35445675075992\\
1210	9.34123054148668\\
1220	9.33553333365749\\
1230	9.32669879580044\\
1240	9.33139466562113\\
1250	9.33094716897035\\
1260	9.33397906672184\\
1270	9.33622951951299\\
1280	9.33689670385554\\
1290	9.33204273200729\\
1300	9.33289994608256\\
1310	9.33871436806212\\
1320	9.33790055081057\\
1330	9.33112672748286\\
1340	9.32993790714486\\
1350	9.3263755559952\\
1360	9.33055929549787\\
1370	9.32048415397337\\
1380	9.33114693895452\\
1390	9.33607684353213\\
1400	9.34098244263161\\
1410	9.33482151550655\\
1420	9.33570973748765\\
1430	9.34037966070291\\
1440	9.33938829203098\\
1450	9.33022415314659\\
1460	9.32771754244324\\
1470	9.32111255599677\\
1480	9.3157936308614\\
1490	9.31592430320623\\
1500	9.32491207829915\\
};
\addplot [color=black, line width = 1.0pt, line join = round]
  table[row sep=crcr]{%
1	51.06269673515\\
10	30.8390415531635\\
20	25.4812512440388\\
30	22.677167299771\\
40	20.4476863224377\\
50	18.1733900805137\\
60	15.6798527048531\\
70	14.3403480778025\\
80	13.3275587762036\\
90	12.7355428397222\\
100	12.280945422949\\
110	11.7963253843743\\
120	11.5423240590915\\
130	11.3645197905095\\
140	11.2702350313871\\
150	11.1518136024146\\
160	10.8732687841834\\
170	10.7576735508574\\
180	10.7543214925404\\
190	10.7074390459893\\
200	10.630254853992\\
210	10.4463515778404\\
220	10.4163111911169\\
230	10.3973477225888\\
240	10.3635561069135\\
250	10.329134439483\\
260	10.181048456066\\
270	10.1591669853367\\
280	10.1290440814337\\
290	10.1194527127041\\
300	10.1098112826341\\
310	9.99136015622351\\
320	9.98092932583493\\
330	10.0054155862308\\
340	9.96815306991626\\
350	9.96664789678054\\
360	9.86995813218723\\
370	9.85544708888042\\
380	9.85469448418017\\
390	9.82793800306579\\
400	9.82740036346912\\
410	9.7766097151159\\
420	9.75342412764795\\
430	9.74934118563526\\
440	9.74476590193148\\
450	9.75660859443406\\
460	9.70279041668725\\
470	9.70259065813201\\
480	9.68256230775587\\
490	9.67777812802726\\
500	9.68316764154756\\
510	9.63826948046408\\
520	9.6349739902845\\
530	9.62609456315974\\
540	9.61602610223794\\
550	9.6080528295762\\
560	9.56627130915232\\
570	9.55757978007048\\
580	9.5650010744754\\
590	9.56581126268617\\
600	9.56334614133861\\
610	9.52556217696111\\
620	9.5307931243535\\
630	9.52409403195767\\
640	9.52440847632303\\
650	9.52301590407619\\
660	9.49702669003911\\
670	9.48344497357717\\
680	9.47895080890557\\
690	9.4736825184023\\
700	9.46875128942916\\
710	9.44372749643354\\
720	9.45213753805837\\
730	9.45885846454595\\
740	9.45079883688633\\
750	9.45117649351594\\
760	9.44689906970825\\
770	9.44398249845712\\
780	9.44452327440962\\
790	9.44532382535395\\
800	9.44676233155824\\
810	9.43849101905112\\
820	9.43390773598756\\
830	9.44099374002973\\
840	9.44770766818878\\
850	9.42741846438613\\
860	9.43031597118018\\
870	9.4390188047989\\
880	9.44087794667289\\
890	9.43117491860419\\
900	9.43510877138012\\
910	9.43347213058267\\
920	9.44973396523025\\
930	9.43547409853133\\
940	9.44338691113269\\
950	9.43572925513176\\
960	9.43587131562979\\
970	9.42633499226613\\
980	9.42798290316488\\
990	9.42628377705769\\
1000	9.41528746475145\\
1010	9.41787093819867\\
1020	9.42862643633778\\
1030	9.42545727122473\\
1040	9.417977248239\\
1050	9.42535188174013\\
1060	9.43415952410505\\
1070	9.42813613345338\\
1080	9.42455828684135\\
1090	9.42467858446553\\
1100	9.42831766769989\\
1110	9.42961360951768\\
1120	9.43381596033185\\
1130	9.42613934442269\\
1140	9.42123031730585\\
1150	9.42722624440586\\
1160	9.42391952871065\\
1170	9.41749265598099\\
1180	9.41968501875785\\
1190	9.41046289096269\\
1200	9.4163561505639\\
1210	9.42087037131381\\
1220	9.42165858532588\\
1230	9.42525250073485\\
1240	9.42494791251518\\
1250	9.42567872502827\\
1260	9.41112129675812\\
1270	9.41260843484064\\
1280	9.41874144596209\\
1290	9.42283131035008\\
1300	9.40994509401916\\
1310	9.41340166928717\\
1320	9.40651411137357\\
1330	9.41410305140837\\
1340	9.40993285408555\\
1350	9.41728731690725\\
1360	9.42212160288723\\
1370	9.40440753363663\\
1380	9.40439294371849\\
1390	9.40397862280571\\
1400	9.40254680829596\\
1410	9.40160640247337\\
1420	9.41222051646693\\
1430	9.40987595784183\\
1440	9.40514454967049\\
1450	9.40869753158951\\
1460	9.40201011699483\\
1470	9.40036860173513\\
1480	9.39686379979047\\
1490	9.39702757866969\\
1500	9.40116896614103\\
};
\addlegendentry{median}
\addplot [draw=none,color=mycolor1, forget plot,  name path=C]
  table[row sep=crcr]{%
1	51.06269673515\\
10	30.9317994550253\\
20	25.6370893585143\\
30	22.8948980202423\\
40	20.779637707316\\
50	18.617367934205\\
60	16.3992068961503\\
70	14.9574177877988\\
80	13.8627277727765\\
90	13.0756091385292\\
100	12.5832763370464\\
110	12.0506067811179\\
120	11.8242758586974\\
130	11.6631484141673\\
140	11.4251317767272\\
150	11.3372055045991\\
160	11.0507488108032\\
170	10.9752877860821\\
180	10.8806174652429\\
190	10.8623813887949\\
200	10.8161674355339\\
210	10.6353743404749\\
220	10.5829388114844\\
230	10.602551307991\\
240	10.5592877223413\\
250	10.4958550339635\\
260	10.3465790220907\\
270	10.3306689177434\\
280	10.2922257081849\\
290	10.2941075767515\\
300	10.3034847800157\\
310	10.1810021362967\\
320	10.1702129146223\\
330	10.1372337113158\\
340	10.1104640238593\\
350	10.1037634874175\\
360	10.0556138045281\\
370	10.0114674940389\\
380	9.99228609460747\\
390	9.99075729601444\\
400	9.9986521521911\\
410	9.91562695768148\\
420	9.89919582344387\\
430	9.89196304010304\\
440	9.89415070857447\\
450	9.89032382115604\\
460	9.81864874492118\\
470	9.81914408653082\\
480	9.82455906775456\\
490	9.82222405641344\\
500	9.7916005344846\\
510	9.77373207556619\\
520	9.75958621234442\\
530	9.7226216053768\\
540	9.72342070990786\\
550	9.72621665521295\\
560	9.668264744356\\
570	9.67418716425327\\
580	9.67447356456039\\
590	9.67334675225956\\
600	9.67740417045557\\
610	9.636756376645\\
620	9.61696645311705\\
630	9.62441748321641\\
640	9.62344338842505\\
650	9.64277088903743\\
660	9.58653048260941\\
670	9.57308050881245\\
680	9.60514560681571\\
690	9.59965644336929\\
700	9.58389458931935\\
710	9.55763997743656\\
720	9.5488060484198\\
730	9.53508849525679\\
740	9.5354294671897\\
750	9.54761423987865\\
760	9.56286179519088\\
770	9.55722907228896\\
780	9.52565385003966\\
790	9.54079964488446\\
800	9.54284604136363\\
810	9.5292020540444\\
820	9.51786577214098\\
830	9.53217198478574\\
840	9.52478321692281\\
850	9.51589323945852\\
860	9.522964491156\\
870	9.52656239614964\\
880	9.52216972448775\\
890	9.50976443998279\\
900	9.51071305319872\\
910	9.53346481437728\\
920	9.51568447026121\\
930	9.51123484590715\\
940	9.5108478741497\\
950	9.50412647218738\\
960	9.49745903283242\\
970	9.51410015533294\\
980	9.51103257047026\\
990	9.52795904270099\\
1000	9.51943515646433\\
1010	9.51198553133923\\
1020	9.50750810086505\\
1030	9.50582416613727\\
1040	9.50462091004397\\
1050	9.51310094701378\\
1060	9.50071413753433\\
1070	9.5017751830961\\
1080	9.50791251556968\\
1090	9.50458406317345\\
1100	9.511845460296\\
1110	9.49460019576956\\
1120	9.50125290533617\\
1130	9.48836612509581\\
1140	9.49659331410285\\
1150	9.49549946171085\\
1160	9.50187510514821\\
1170	9.50368212826099\\
1180	9.48628141792691\\
1190	9.48373831774359\\
1200	9.49202238922988\\
1210	9.50159390528164\\
1220	9.51499378647674\\
1230	9.49487423697492\\
1240	9.4888055238025\\
1250	9.49803139638214\\
1260	9.48710036066065\\
1270	9.48545398735125\\
1280	9.49266140765823\\
1290	9.49555957415297\\
1300	9.50180773225597\\
1310	9.48397793296929\\
1320	9.50095562197923\\
1330	9.50134311751813\\
1340	9.50083078231761\\
1350	9.51108325488388\\
1360	9.4915542417176\\
1370	9.50289942752352\\
1380	9.49901219771062\\
1390	9.48823745938401\\
1400	9.50820544379257\\
1410	9.50492736591765\\
1420	9.49103170997871\\
1430	9.49547124925738\\
1440	9.48174910663836\\
1450	9.48684775382147\\
1460	9.49251788421108\\
1470	9.49521483433413\\
1480	9.49929328243426\\
1490	9.50301847902833\\
1500	9.48357831727158\\
};
\addplot [draw=none,color=mycolor1, forget plot,  name path=D]
  table[row sep=crcr]{%
1	51.06269673515\\
10	31.0388972435996\\
20	25.7562045656951\\
30	23.1128652741256\\
40	20.9876147932761\\
50	18.8082354445495\\
60	17.0186488598665\\
70	15.6821602093114\\
80	14.3146388652129\\
90	13.5476573319515\\
100	12.9584482750158\\
110	12.2699940333701\\
120	12.0679084980873\\
130	11.8481210620335\\
140	11.7664693849057\\
150	11.565431709941\\
160	11.2963677883904\\
170	11.2956283229696\\
180	11.1990720323224\\
190	11.193390718738\\
200	11.0702867493404\\
210	10.8813455969135\\
220	10.8156117644357\\
230	10.8075958356257\\
240	10.7709112850448\\
250	10.699318423144\\
260	10.5565412103431\\
270	10.5255420333954\\
280	10.5634339812548\\
290	10.537996277596\\
300	10.4633505421463\\
310	10.4090962146701\\
320	10.3804867170837\\
330	10.338000926475\\
340	10.2821948054906\\
350	10.2801053992287\\
360	10.1917043007965\\
370	10.1613272921745\\
380	10.1163402426694\\
390	10.1201620297272\\
400	10.1067622741194\\
410	10.037395349369\\
420	10.0735239771336\\
430	10.0323446126702\\
440	10.0094504762317\\
450	9.97612216068796\\
460	10.0049509331659\\
470	10.0107775517102\\
480	9.94575927620494\\
490	9.89613512069196\\
500	9.88712277423384\\
510	9.85254587662499\\
520	9.87135950848154\\
530	9.87208279400926\\
540	9.84246015309278\\
550	9.83015880309609\\
560	9.79263690052551\\
570	9.81754241564687\\
580	9.79864003510372\\
590	9.78633259808013\\
600	9.81544650482716\\
610	9.76286799517253\\
620	9.75995638829515\\
630	9.75669543491489\\
640	9.75515111236749\\
650	9.76365221164415\\
660	9.71505787284097\\
670	9.74000931752465\\
680	9.70389753526609\\
690	9.71209923039083\\
700	9.72937283903726\\
710	9.69108850526105\\
720	9.70204298564705\\
730	9.68136079471354\\
740	9.65680013932045\\
750	9.66460703107118\\
760	9.68388128594281\\
770	9.69449676114857\\
780	9.71343720027427\\
790	9.65203205332605\\
800	9.64512303482807\\
810	9.66182004696725\\
820	9.63991446959874\\
830	9.64088137794283\\
840	9.64624823781797\\
850	9.63478147375755\\
860	9.64312488426276\\
870	9.65658852863998\\
880	9.67511173422202\\
890	9.67426474865048\\
900	9.66080000353332\\
910	9.6623338486336\\
920	9.65755983929843\\
930	9.63986438526242\\
940	9.67266062194874\\
950	9.65750415022092\\
960	9.65390875687888\\
970	9.64662612033792\\
980	9.64856340494729\\
990	9.6598156201276\\
1000	9.65862932370016\\
1010	9.64927289741325\\
1020	9.63351730134055\\
1030	9.64300761168699\\
1040	9.63217457806336\\
1050	9.64561081771077\\
1060	9.64256502381656\\
1070	9.64184301838915\\
1080	9.63812696134233\\
1090	9.64426732302887\\
1100	9.65375257147926\\
1110	9.64051809350355\\
1120	9.63313023959793\\
1130	9.6207430056766\\
1140	9.6218100652798\\
1150	9.63919131122048\\
1160	9.63085927137926\\
1170	9.62411621785455\\
1180	9.61410329981925\\
1190	9.6233696428163\\
1200	9.61308566193801\\
1210	9.62197880929645\\
1220	9.61766049137031\\
1230	9.63019981864876\\
1240	9.64340967862391\\
1250	9.65557953300655\\
1260	9.64681987925491\\
1270	9.64128647634491\\
1280	9.65733978853039\\
1290	9.64505336352416\\
1300	9.65071697310213\\
1310	9.63515581813248\\
1320	9.63422144308767\\
1330	9.63663830024985\\
1340	9.63467422108673\\
1350	9.63096786906126\\
1360	9.63228327323333\\
1370	9.63448066937773\\
1380	9.64709746997684\\
1390	9.65023104962351\\
1400	9.63394543581724\\
1410	9.64132974073482\\
1420	9.62703109273594\\
1430	9.62114404306644\\
1440	9.63138711631926\\
1450	9.6175570329505\\
1460	9.62730317133195\\
1470	9.63090232473531\\
1480	9.63340883899113\\
1490	9.63264857204631\\
1500	9.6305729726543\\
};
\addplot[fill = mycolor1, fill opacity = 0.2, draw = none, area legend] fill between [of=A and B];
\addplot[fill = mycolor1, fill opacity = 0.55, draw = none, area legend] fill between [of=B and C];
\addplot[fill = mycolor1, fill opacity = 0.2, draw = none, area legend] fill between [of=C and D];
\addlegendentry{$\mathcal{Q}_{0.1,0.9}$}
\addlegendentry{$\mathcal{Q}_{0.25,0.75}$}
\end{axis}
\end{tikzpicture}%

%% file: Graphics/Beam/beam_setup.tex
\begin{tikzpicture}[x=0.9\linewidth/200,y=0.9\linewidth/200]
\draw[very thick] (10,5) -- (190,5) -- (190,65) -- (10,65) -- cycle;
\fill[pattern=north west lines] (0,5) rectangle (10,65);
\fill[pattern=north east lines,pattern color = blue] (35,10) rectangle (57,32);
\draw[blue,thick] (35,10) rectangle (57,32);
\draw[Latex-Latex] (10,0)--node[midway,below]{$3\ell$}(190,0);
\draw[Latex-Latex] (-5,5)--node[left,midway]{$\ell$}(-5,65);
\draw[Latex-,red,ultra thick] (192,10) -- (192,35);
\draw[thick, dashed, teal] (180,65) -- (180,5);
\end{tikzpicture}

%% file: Graphics/Load/load_setup.tex
\begin{tikzpicture}[x=0.9\linewidth/200,y=0.9\linewidth/200]
\draw[very thick] (10,5) -- (190,5) -- (190,50) -- (10,50) -- cycle;
\fill[pattern=north west lines] (0,5) rectangle (10,27.5);
\fill[pattern=north west lines] (190,5) rectangle (200,27.5);
\draw[Latex-Latex] (10,0)--node[midway,below]{$4\ell$}(190,0);
\draw[Latex-Latex] (-5,5)--node[left,midway]{$\ell$}(-5,50);
\draw[Latex-,red] (40,52) -- (40,65);
\draw[very thick,red] (10,50) -- (190,50);
\end{tikzpicture}

%% file: Graphics/Load/Distribution_histogram.tex
%
%
\definecolor{mycolor1}{rgb}{0.00000,0.44700,0.74100}%
\definecolor{mycolor2}{rgb}{0.85000,0.32500,0.09800}%
\begin{tikzpicture}

\begin{axis}[%
width=.9\linewidth,
height=.9\linewidth,
at={(0.758in,0.481in)},
xmin=-7,
xmax=370,
ymin=0,
ymax=0.009,
axis background/.style={fill=white},
xlabel={load case},
ylabel={probability}
]
\addplot[ybar interval, fill=mycolor1, fill opacity=.8, draw=none, area legend,line width=.1pt] table[row sep=crcr] {%
x	y\\
0.9	0.00016\\
1.898	0.00039\\
2.896	0.0004725\\
3.894	0.0004725\\
4.892	0.0004725\\
5.89	0.0005225\\
6.888	0.000575\\
7.886	0.0006525\\
8.884	0.00063\\
9.882	0.000705\\
10.88	0.0007525\\
11.878	0.0007525\\
12.876	0.0008325\\
13.874	0.0009475\\
14.872	0.00097\\
15.87	0.0010275\\
16.868	0.000985\\
17.866	0.0010425\\
18.864	0.0011775\\
19.862	0.0013325\\
20.86	0.001245\\
21.858	0.001295\\
22.856	0.001365\\
23.854	0.001655\\
24.852	0.0015425\\
25.85	0.0017375\\
26.848	0.0017425\\
27.846	0.001825\\
28.844	0.0019025\\
29.842	0.002\\
30.84	0.00208\\
31.838	0.002215\\
32.836	0.002285\\
33.834	0.0024325\\
34.832	0.00262\\
35.83	0.002485\\
36.828	0.0026775\\
37.826	0.002955\\
38.824	0.0031075\\
39.822	0.0030875\\
40.82	0.00329\\
41.818	0.0034425\\
42.816	0.003565\\
43.814	0.0037525\\
44.812	0.0038\\
45.81	0.0037125\\
46.808	0.003955\\
47.806	0.0041925\\
48.804	0.0042225\\
49.802	0.004455\\
50.8	0.0044075\\
51.798	0.0046475\\
52.796	0.0049925\\
53.794	0.00513\\
54.792	0.0054525\\
55.79	0.0054225\\
56.788	0.0055175\\
57.786	0.00549\\
58.784	0.0058\\
59.782	0.0059825\\
60.78	0.00592\\
61.778	0.00619\\
62.776	0.006135\\
63.774	0.0065775\\
64.772	0.0067025\\
65.77	0.00681\\
66.768	0.0068925\\
67.766	0.0072225\\
68.764	0.0072775\\
69.762	0.0070775\\
70.76	0.0072575\\
71.758	0.0073425\\
72.756	0.00741\\
73.754	0.0077525\\
74.752	0.0077675\\
75.75	0.00795\\
76.748	0.00788\\
77.746	0.0082575\\
78.744	0.008385\\
79.742	0.00833\\
80.74	0.0081725\\
81.738	0.0084225\\
82.736	0.008285\\
83.734	0.008455\\
84.732	0.0086325\\
85.73	0.008795\\
86.728	0.0085425\\
87.726	0.0086475\\
88.724	0.0083975\\
89.722	0.008455\\
90.72	0.0086625\\
91.718	0.0085625\\
92.716	0.00884\\
93.714	0.0088425\\
94.712	0.0086275\\
95.71	0.008725\\
96.708	0.0085275\\
97.706	0.0085225\\
98.704	0.0088175\\
99.702	0.008275\\
100.7	0.00843\\
101.698	0.00846\\
102.696	0.0082775\\
103.694	0.0082875\\
104.692	0.008395\\
105.69	0.0078975\\
106.688	0.0078475\\
107.686	0.0080825\\
108.684	0.007775\\
109.682	0.00754\\
110.68	0.0074575\\
111.678	0.0076875\\
112.676	0.00745\\
113.674	0.0073125\\
114.672	0.007235\\
115.67	0.007135\\
116.668	0.0070575\\
117.666	0.0070475\\
118.664	0.006825\\
119.662	0.00653\\
120.66	0.0064175\\
121.658	0.0064475\\
122.656	0.0061425\\
123.654	0.0058575\\
124.652	0.0056175\\
125.65	0.005695\\
126.648	0.005625\\
127.646	0.00561\\
128.644	0.0054825\\
129.642	0.0054725\\
130.64	0.00505\\
131.638	0.005195\\
132.636	0.00485\\
133.634	0.00499\\
134.632	0.0047325\\
135.63	0.0046775\\
136.628	0.0043075\\
137.626	0.004265\\
138.624	0.0040975\\
139.622	0.00403\\
140.62	0.003915\\
141.618	0.00406\\
142.616	0.00397\\
143.614	0.003565\\
144.612	0.003605\\
145.61	0.0036025\\
146.608	0.00331\\
147.606	0.003195\\
148.604	0.0032225\\
149.602	0.0032125\\
150.6	0.00303\\
151.598	0.0028925\\
152.596	0.00282\\
153.594	0.002725\\
154.592	0.002765\\
155.59	0.00264\\
156.588	0.002585\\
157.586	0.0024325\\
158.584	0.002585\\
159.582	0.002395\\
160.58	0.00227\\
161.578	0.002225\\
162.576	0.0021825\\
163.574	0.0020475\\
164.572	0.0021\\
165.57	0.0021475\\
166.568	0.0020425\\
167.566	0.001945\\
168.564	0.002005\\
169.562	0.0019725\\
170.56	0.001805\\
171.558	0.0018675\\
172.556	0.00181\\
173.554	0.0017\\
174.552	0.0017875\\
175.55	0.0017025\\
176.548	0.00177\\
177.546	0.0017425\\
178.544	0.0015925\\
179.542	0.0016725\\
180.54	0.0016875\\
181.538	0.00173\\
182.536	0.001645\\
183.534	0.001545\\
184.532	0.00161\\
185.53	0.001575\\
186.528	0.001605\\
187.526	0.00159\\
188.524	0.001595\\
189.522	0.001565\\
190.52	0.0015875\\
191.518	0.0015925\\
192.516	0.00153\\
193.514	0.001515\\
194.512	0.001475\\
195.51	0.00156\\
196.508	0.0015775\\
197.506	0.00149\\
198.504	0.0014725\\
199.502	0.00144\\
200.5	0.0014625\\
201.498	0.0015625\\
202.496	0.0013725\\
203.494	0.0014725\\
204.492	0.0014175\\
205.49	0.00147\\
206.488	0.0014525\\
207.486	0.00153\\
208.484	0.0014525\\
209.482	0.0014\\
210.48	0.00147\\
211.478	0.0015075\\
212.476	0.00144\\
213.474	0.001455\\
214.472	0.0013925\\
215.47	0.0014925\\
216.468	0.0013475\\
217.466	0.00144\\
218.464	0.0014925\\
219.462	0.00143\\
220.46	0.001395\\
221.458	0.001375\\
222.456	0.0014075\\
223.454	0.0013825\\
224.452	0.0013775\\
225.45	0.0013875\\
226.448	0.0013625\\
227.446	0.0013425\\
228.444	0.00144\\
229.442	0.0013125\\
230.44	0.0014775\\
231.438	0.00139\\
232.436	0.001325\\
233.434	0.0014325\\
234.432	0.0014975\\
235.43	0.0013325\\
236.428	0.0014\\
237.426	0.0013825\\
238.424	0.0013025\\
239.422	0.0013975\\
240.42	0.0013625\\
241.418	0.001395\\
242.416	0.0012525\\
243.414	0.001245\\
244.412	0.0012825\\
245.41	0.001325\\
246.408	0.00128\\
247.406	0.001345\\
248.404	0.0013275\\
249.402	0.001275\\
250.4	0.0012575\\
251.398	0.0012625\\
252.396	0.00126\\
253.394	0.0012275\\
254.392	0.001135\\
255.39	0.0012125\\
256.388	0.0011775\\
257.386	0.001235\\
258.384	0.00118\\
259.382	0.0011125\\
260.38	0.00116\\
261.378	0.0011075\\
262.376	0.00121\\
263.374	0.0010525\\
264.372	0.001195\\
265.37	0.0011375\\
266.368	0.0011325\\
267.366	0.0010475\\
268.364	0.0010975\\
269.362	0.0011275\\
270.36	0.0010225\\
271.358	0.0011075\\
272.356	0.0010075\\
273.354	0.0010025\\
274.352	0.0009575\\
275.35	0.0009975\\
276.348	0.001005\\
277.346	0.0010025\\
278.344	0.0009775\\
279.342	0.0010025\\
280.34	0.0009325\\
281.338	0.0009275\\
282.336	0.000995\\
283.334	0.0009125\\
284.332	0.00091\\
285.33	0.0008675\\
286.328	0.0008875\\
287.326	0.0008525\\
288.324	0.00085\\
289.322	0.0007975\\
290.32	0.0008275\\
291.318	0.000805\\
292.316	0.000765\\
293.314	0.000715\\
294.312	0.0007625\\
295.31	0.00078\\
296.308	0.0007975\\
297.306	0.000845\\
298.304	0.00083\\
299.302	0.000745\\
300.3	0.000765\\
301.298	0.0007225\\
302.296	0.0006325\\
303.294	0.0006575\\
304.292	0.0006675\\
305.29	0.000715\\
306.288	0.00064\\
307.286	0.00066\\
308.284	0.000595\\
309.282	0.0006075\\
310.28	0.0005375\\
311.278	0.00057\\
312.276	0.00059\\
313.274	0.00051\\
314.272	0.0005075\\
315.27	0.0005675\\
316.268	0.0005375\\
317.266	0.0005325\\
318.264	0.00048\\
319.262	0.0005375\\
320.26	0.0005675\\
321.258	0.00048\\
322.256	0.000485\\
323.254	0.0004675\\
324.252	0.000435\\
325.25	0.000415\\
326.248	0.0004625\\
327.246	0.00037\\
328.244	0.000455\\
329.242	0.00044\\
330.24	0.0003425\\
331.238	0.0004575\\
332.236	0.0003575\\
333.234	0.0003775\\
334.232	0.00042\\
335.23	0.0003075\\
336.228	0.000365\\
337.226	0.0003875\\
338.224	0.0003625\\
339.222	0.000345\\
340.22	0.0002925\\
341.218	0.0003425\\
342.216	0.0002875\\
343.214	0.0003375\\
344.212	0.0003\\
345.21	0.00028\\
346.208	0.0003075\\
347.206	0.0002875\\
348.204	0.00029\\
349.202	0.0002375\\
350.2	0.0002125\\
351.198	0.00023\\
352.196	0.000245\\
353.194	0.0002425\\
354.192	0.00024\\
355.19	0.0002075\\
356.188	0.000225\\
357.186	0.000205\\
358.184	0.0002175\\
359.182	0.000215\\
360.18	8.25e-05\\
361.178	8.25e-05\\
};
\end{axis}
\end{tikzpicture}%

%% file: Graphics/Dynamic/dyna_setup.tex
\begin{tikzpicture}[x=0.9\linewidth/200,y=0.9\linewidth/200]
\draw[very thick] (10,5) -- (190,5) -- (190,65) -- (10,65) -- cycle;
\draw (2.07,10) -- (9,14) -- (2.07,18) -- cycle;
\draw (2.07,60) -- (9,56) -- (2.07,52) -- cycle;
\draw (189,4) -- (193,-2.928) -- (185,-2.928) -- cycle;
\draw (0,12) circle (2);
\draw (0,16) circle (2);
\draw (0,58) circle (2);
\draw (0,54) circle (2);
\draw (191,-4.928) circle (2);
\draw (187,-4.928) circle (2);
\filldraw[fill=black!55] (-3,9) rectangle (-2,19);
\filldraw[fill=black!55] (-3,51) rectangle (-2,61);
\filldraw[fill=black!55] (184,-7.928) rectangle (194,-6.928);
\draw[Latex-Latex] (10,-12)--node[midway,below]{$3\ell$}(190,-12);
\draw[Latex-Latex] (-8,5)--node[left,midway]{$\ell$}(-8,65);
\draw[Latex-,red] (10,67) -- (10,82);
\draw[Latex-,red] (16,67) -- (16,82);
\draw[thin,red] (10,67) -- (16,67);
\end{tikzpicture}

%% file: Graphics/Dynamic/Dynamic_compliances.tex
%
%
\definecolor{mycolor1}{rgb}{0.00000,0.44700,0.74100}%
\definecolor{mycolor2}{rgb}{0.85000,0.32500,0.09800}%
\definecolor{mycolor3}{rgb}{0.92900,0.69400,0.12500}%
\begin{tikzpicture}

\begin{axis}[%
width=.9\linewidth,
height=.9\linewidth,
at={(0.758in,0.481in)},
xmin=0,
xmax=0.2,
ymode=log,
ymin=50,
ymax=2000,
yminorticks=true,
axis background/.style={fill=white},
xmajorgrids,
ymajorgrids,
yminorgrids,
xlabel = {$\omega$},
legend pos = south west,
x tick label style={
        /pgf/number format/.cd,
            fixed,
            precision=2,
        /tikz/.cd
    }
]
\addplot [color=mycolor1, line width = 1.5pt, line join = round]
  table[row sep=crcr]{%
0	207.834348778931\\
0.000803212851405623	207.840036652995\\
0.00160642570281125	207.857102455165\\
0.00240963855421687	207.885552751939\\
0.00321285140562249	207.925398477646\\
0.00401606425702811	207.976654971998\\
0.00481927710843374	208.039341996542\\
0.00562248995983936	208.113483725364\\
0.00642570281124498	208.199108754774\\
0.0072289156626506	208.29625018752\\
0.00803212851405622	208.40494565739\\
0.00883534136546185	208.52523728814\\
0.00963855421686747	208.657171855665\\
0.0104417670682731	208.800800703987\\
0.0112449799196787	208.956179903557\\
0.0120481927710843	209.123370268832\\
0.01285140562249	209.302437355532\\
0.0136546184738956	209.49345165191\\
0.0144578313253012	209.696488546541\\
0.0152610441767068	209.911628419101\\
0.0160642570281124	210.13895678103\\
0.0168674698795181	210.378564277136\\
0.0176706827309237	210.630546845719\\
0.0184738955823293	210.895005783585\\
0.0192771084337349	211.172047851109\\
0.0200803212851406	211.461785382535\\
0.0208835341365462	211.764336417834\\
0.0216867469879518	212.079824829581\\
0.0224899598393574	212.40838039886\\
0.0232931726907631	212.75013901696\\
0.0240963855421687	213.105242801143\\
0.0248995983935743	213.47384025696\\
0.0257028112449799	213.856086410149\\
0.0265060240963855	214.252143043861\\
0.0273092369477912	214.662178773047\\
0.0281124497991968	215.086369348539\\
0.0289156626506024	215.524897724211\\
0.029718875502008	215.977954423096\\
0.0305220883534137	216.445737576534\\
0.0313253012048193	216.928453283693\\
0.0321285140562249	217.426315819232\\
0.0329317269076305	217.939547836501\\
0.0337349397590361	218.468380694343\\
0.0345381526104418	219.013054709237\\
0.0353413654618474	219.573819398496\\
0.036144578313253	220.150933859365\\
0.0369477911646586	220.744667088111\\
0.0377510040160643	221.355298241575\\
0.0385542168674699	221.983117040349\\
0.0393574297188755	222.628424158918\\
0.0401606425702811	223.291531555383\\
0.0409638554216867	223.972762930878\\
0.0417670682730924	224.672454157409\\
0.042570281124498	225.390953670326\\
0.0433734939759036	226.128623043819\\
0.0441767068273092	226.885837380347\\
0.0449799196787149	227.662985880936\\
0.0457831325301205	228.460472430077\\
0.0465863453815261	229.278716139173\\
0.0473895582329317	230.118151932602\\
0.0481927710843374	230.979231242061\\
0.048995983935743	231.862422651005\\
0.0497991967871486	232.768212586227\\
0.0506024096385542	233.697106136224\\
0.0514056224899598	234.649627754077\\
0.0522088353413655	235.626322162941\\
0.0530120481927711	236.627755191162\\
0.0538152610441767	237.654514712157\\
0.0546184738955823	238.707211621858\\
0.055421686746988	239.786480863831\\
0.0562248995983936	240.892982533544\\
0.0570281124497992	242.027403045396\\
0.0578313253012048	243.190456312876\\
0.0586345381526104	244.382885036033\\
0.0594377510040161	245.605462095441\\
0.0602409638554217	246.858992050962\\
0.0610441767068273	248.144312494015\\
0.0618473895582329	249.462295898818\\
0.0626506024096386	250.813851134162\\
0.0634538152610442	252.199925457233\\
0.0642570281124498	253.621506325447\\
0.0650602409638554	255.079623530185\\
0.0658634538152611	256.575351363683\\
0.0666666666666667	258.109810900595\\
0.0674698795180723	259.684172531956\\
0.0682730923694779	261.299658587827\\
0.0690763052208835	262.9575461194\\
0.0698795180722892	264.65916989555\\
0.0706827309236948	266.405925605424\\
0.0714859437751004	268.199273243863\\
0.072289156626506	270.040740757317\\
0.0730923694779116	271.931927932045\\
0.0738955823293173	273.874510532465\\
0.0746987951807229	275.870244772348\\
0.0755020080321285	277.92097208915\\
0.0763052208835341	280.028624197878\\
0.0771084337349398	282.195228592513\\
0.0779116465863454	284.422914464248\\
0.078714859437751	286.713918919411\\
0.0795180722891566	289.070593937015\\
0.0803212851405622	291.495413452435\\
0.0811244979919679	293.990981448869\\
0.0819277108433735	296.56004028255\\
0.0827309236947791	299.205479934386\\
0.0835341365461847	301.930347859771\\
0.0843373493975904	304.737859633585\\
0.085140562248996	307.631410578932\\
0.0859437751004016	310.614588292715\\
0.0867469879518072	313.691186166558\\
0.0875502008032129	316.865218141034\\
0.0883534136546185	320.140934720983\\
0.0891566265060241	323.522840354568\\
0.0899598393574297	327.015712394745\\
0.0907630522088353	330.624621724559\\
0.091566265060241	334.354955356248\\
0.0923694779116466	338.212440884517\\
0.0931726907630522	342.203173622764\\
0.0939759036144578	346.333645865212\\
0.0947791164658635	350.610779272547\\
0.0955823293172691	355.04196049104\\
0.0963855421686747	359.635079756671\\
0.0971887550200803	364.398574265\\
0.097991967871486	369.341475120326\\
0.0987951807228916	374.473459759973\\
0.0995983935742972	379.804909358463\\
0.100401606425703	385.346972949077\\
0.101204819277108	391.111638147315\\
0.102008032128514	397.111809542619\\
0.10281124497992	403.361396067707\\
0.103614457831325	409.875408434932\\
0.104417670682731	416.670067281682\\
0.105220883534137	423.762924481893\\
0.106024096385542	431.172998536325\\
0.106827309236948	438.920926316937\\
0.107630522088353	447.029133563091\\
0.108433734939759	455.52202598719\\
0.109236947791165	464.426205219209\\
0.11004016064257	473.770711346622\\
0.110843373493976	483.587297211599\\
0.111646586345382	493.910738532497\\
0.112449799196787	504.77918509897\\
0.113253012048193	516.234557551699\\
0.114056224899598	528.322999479061\\
0.114859437751004	541.095388204434\\
0.11566265060241	554.607916609624\\
0.116465863453815	568.922752084344\\
0.117269076305221	584.108784973223\\
0.118072289156627	600.242475706522\\
0.118875502008032	617.40881161568\\
0.119678714859438	635.702384173297\\
0.120481927710843	655.228593198353\\
0.121285140562249	676.104981776524\\
0.122088353413655	698.462696215294\\
0.12289156626506	722.44805266342\\
0.123694779116466	748.224165028896\\
0.124497991967871	775.97255336208\\
0.125301204819277	805.894589347088\\
0.126104417670683	838.212537340623\\
0.126907630522088	873.169803208589\\
0.127710843373494	911.029771192165\\
0.1285140562249	952.072260513396\\
0.129317269076305	996.586116554428\\
0.130120481927711	1044.8556995905\\
0.130923694779116	1097.13798327033\\
0.131726907630522	1153.62559022076\\
0.132530120481928	1214.38951498041\\
0.133333333333333	1279.29387794347\\
0.134136546184739	1347.87507798522\\
0.134939759036145	1419.18149402742\\
0.13574297188755	1491.58171945394\\
0.136546184738956	1562.57458346015\\
0.137349397590361	1628.67544350532\\
0.138152610441767	1685.49987834603\\
0.138955823293173	1728.18028156494\\
0.139759036144578	1752.17134961188\\
0.140562248995984	1754.29553288519\\
0.14136546184739	1733.64139734666\\
0.142168674698795	1691.87410587952\\
0.142971887550201	1632.78126819129\\
0.143775100401606	1561.28962123974\\
0.144578313253012	1482.41290371562\\
0.145381526104418	1400.49608566371\\
0.146184738955823	1318.86815100129\\
0.146987951807229	1239.82016023973\\
0.147791164658635	1164.76206228891\\
0.14859437751004	1094.43744910952\\
0.149397590361446	1029.12666478234\\
0.150200803212851	968.810194131898\\
0.151004016064257	913.288262520648\\
0.151807228915663	862.262930281448\\
0.152610441767068	815.391710057664\\
0.153413654618474	772.321004801993\\
0.15421686746988	732.705920176097\\
0.155020080321285	696.221185215111\\
0.155823293172691	662.566445139282\\
0.156626506024096	631.468108118295\\
0.157429718875502	602.679155255499\\
0.158232931726908	575.977820804303\\
0.159036144578313	551.165701468786\\
0.159839357429719	528.065639039188\\
0.160642570281124	506.519577461529\\
0.16144578313253	486.386510451288\\
0.162248995983936	467.540574944429\\
0.163052208835341	449.869319005765\\
0.163855421686747	433.272147031735\\
0.164658634538153	417.658937749379\\
0.165461847389558	402.948823979121\\
0.166265060240964	389.069116665577\\
0.167068273092369	375.954362103958\\
0.167871485943775	363.545513731914\\
0.168674698795181	351.789207351356\\
0.169477911646586	340.637125332648\\
0.170281124497992	330.045439452154\\
0.171084337349398	319.974321899375\\
0.171887550200803	310.387515642987\\
0.172690763052209	301.251956500436\\
0.173493975903614	292.537440571431\\
0.17429718875502	284.216330378171\\
0.175100401606426	276.26329562684\\
0.175903614457831	268.655083588552\\
0.176706827309237	261.370315594419\\
0.177510040160643	254.389306355008\\
0.178313253012048	247.693903123817\\
0.179116465863454	241.267342613217\\
0.179919678714859	235.094123131638\\
0.180722891566265	229.159890255999\\
0.181526104417671	223.451334535131\\
0.182329317269076	217.956099805977\\
0.183132530120482	212.662700562327\\
0.183935742971888	207.560447971319\\
0.184738955823293	202.639382918024\\
0.185542168674699	197.890215630015\\
0.186345381526104	193.304271176519\\
0.18714859437751	188.873440204055\\
0.187951807228916	184.59013396836\\
0.188755020080321	180.447243886776\\
0.189558232931727	176.438104566812\\
0.190361445783133	172.556460162908\\
0.191164658634538	168.796433746153\\
0.191967871485944	165.152499374663\\
0.192771084337349	161.619456535626\\
0.193574297188755	158.192406788671\\
0.194377510040161	154.866732354899\\
0.195180722891566	151.638076495715\\
0.195983935742972	148.502325514483\\
0.196787148594378	145.455592175068\\
0.197590361445783	142.494200477035\\
0.198393574297189	139.614671654383\\
0.199196787148594	136.813711267899\\
0.2	134.088197175672\\
};
\addlegendentry{static};
\addplot [color=mycolor2, line width = 1.5pt, line join = round]
  table[row sep=crcr]{%
0	239.423482718685\\
0.000803212851405623	239.429582697922\\
0.00160642570281125	239.447888156154\\
0.00240963855421687	239.478415681609\\
0.00321285140562249	239.521192961324\\
0.00401606425702811	239.576258855542\\
0.00481927710843374	239.643663498221\\
0.00562248995983936	239.723468384177\\
0.00642570281124498	239.815746518162\\
0.0072289156626506	239.920582622918\\
0.00803212851405622	240.038073337836\\
0.00883534136546185	240.168327392219\\
0.00963855421686747	240.311465964759\\
0.0104417670682731	240.467622871491\\
0.0112449799196787	240.636944990891\\
0.0120481927710843	240.819592582677\\
0.01285140562249	241.015739650628\\
0.0136546184738956	241.225574473062\\
0.0144578313253012	241.449300006786\\
0.0152610441767068	241.687134416756\\
0.0160642570281124	241.939311693118\\
0.0168674698795181	242.206082196601\\
0.0176706827309237	242.487713384987\\
0.0184738955823293	242.784490501346\\
0.0192771084337349	243.096717359271\\
0.0200803212851406	243.424717183211\\
0.0208835341365462	243.768833532317\\
0.0216867469879518	244.129431278309\\
0.0224899598393574	244.506897627414\\
0.0232931726907631	244.901643312352\\
0.0240963855421687	245.31410379417\\
0.0248995983935743	245.744740606917\\
0.0257028112449799	246.194042764581\\
0.0265060240963855	246.662528363771\\
0.0273092369477912	247.150746168248\\
0.0281124497991968	247.659277512456\\
0.0289156626506024	248.188738153529\\
0.029718875502008	248.739780491772\\
0.0305220883534137	249.31309572883\\
0.0313253012048193	249.909416433502\\
0.0321285140562249	250.529519198338\\
0.0329317269076305	251.174227511557\\
0.0337349397590361	251.844414945259\\
0.0345381526104418	252.541008573454\\
0.0353413654618474	253.264992651058\\
0.036144578313253	254.01741272433\\
0.0369477911646586	254.799380048638\\
0.0377510040160643	255.612076339265\\
0.0385542168674699	256.456759069717\\
0.0393574297188755	257.334767231126\\
0.0401606425702811	258.24752757249\\
0.0409638554216867	259.196561536242\\
0.0417670682730924	260.183492807653\\
0.042570281124498	261.210055568899\\
0.0433734939759036	262.278103734488\\
0.0441767068273092	263.389620896596\\
0.0449799196787149	264.54673148412\\
0.0457831325301205	265.751713054264\\
0.0465863453815261	267.007009831628\\
0.0473895582329317	268.315247734526\\
0.0481927710843374	269.679251128928\\
0.048995983935743	271.10206137458\\
0.0497991967871486	272.586957516489\\
0.0506024096385542	274.137479449071\\
0.0514056224899598	275.757453666573\\
0.0522088353413655	277.451022313801\\
0.0530120481927711	279.222675638309\\
0.0538152610441767	281.077288600105\\
0.0546184738955823	283.020162125334\\
0.055421686746988	285.057069704895\\
0.0562248995983936	287.194310269073\\
0.0570281124497992	289.438768132488\\
0.0578313253012048	291.797981266498\\
0.0586345381526104	294.280219157962\\
0.0594377510040161	296.894572047908\\
0.0602409638554217	299.651053262234\\
0.0610441767068273	302.560716667595\\
0.0618473895582329	305.635792789553\\
0.0626506024096386	308.889845678036\\
0.0634538152610442	312.337955589012\\
0.0642570281124498	315.996931221354\\
0.0650602409638554	319.885557967922\\
0.0658634538152611	324.024888691064\\
0.0666666666666667	328.438585465079\\
0.0674698795180723	333.153322889232\\
0.0682730923694779	338.199265066291\\
0.0690763052208835	343.610631751724\\
0.0698795180722892	349.426372144529\\
0.0706827309236948	355.69096900921\\
0.0714859437751004	362.455400184907\\
0.072289156626506	369.778290114702\\
0.0730923694779116	377.727289372048\\
0.0738955823293173	386.380725464743\\
0.0746987951807229	395.829571497599\\
0.0755020080321285	406.179776800511\\
0.0763052208835341	417.554990304799\\
0.0771084337349398	430.099670423879\\
0.0779116465863454	443.982492244184\\
0.078714859437751	459.399792774717\\
0.0795180722891566	476.578470352043\\
0.0803212851405622	495.777126815934\\
0.0811244979919679	517.283035765224\\
0.0819277108433735	541.400058404238\\
0.0827309236947791	568.417842436846\\
0.0835341365461847	598.544591388343\\
0.0843373493975904	631.773917682391\\
0.085140562248996	667.640737020711\\
0.0859437751004016	704.811841729692\\
0.0867469879518072	740.497612949614\\
0.0875502008032129	769.862620055718\\
0.0883534136546185	786.04920273274\\
0.0891566265060241	781.821618930343\\
0.0899598393574297	753.055021260985\\
0.0907630522088353	701.79688806913\\
0.091566265060241	635.650497820509\\
0.0923694779116466	563.830818010818\\
0.0931726907630522	493.575083571183\\
0.0939759036144578	429.00833907749\\
0.0947791164658635	371.74462750113\\
0.0955823293172691	321.904242442675\\
0.0963855421686747	278.900970083828\\
0.0971887550200803	241.904646923441\\
0.097991967871486	210.071138121856\\
0.0987951807228916	182.639932350291\\
0.0995983935742972	158.965508470042\\
0.100401606425703	138.518994584605\\
0.101204819277108	120.878414248098\\
0.102008032128514	105.715851312295\\
0.10281124497992	92.784605243733\\
0.103614457831325	81.9064772345469\\
0.104417670682731	72.9574427007595\\
0.105220883534137	65.8488932800056\\
0.106024096385542	60.5023569726466\\
0.106827309236948	56.8197404538476\\
0.107630522088353	54.6581022162611\\
0.108433734939759	53.8215839919718\\
0.109236947791165	54.0760839936732\\
0.11004016064257	55.1783036130471\\
0.110843373493976	56.9037334557803\\
0.111646586345382	59.063391433465\\
0.112449799196787	61.50844033813\\
0.113253012048193	64.1269522422292\\
0.114056224899598	66.8373863700262\\
0.114859437751004	69.5816758221469\\
0.11566265060241	72.3192248704913\\
0.116465863453815	75.0221640789741\\
0.117269076305221	77.6717812769845\\
0.118072289156627	80.2559073462205\\
0.118875502008032	82.7670276401895\\
0.119678714859438	85.2009262523849\\
0.120481927710843	87.5557146910235\\
0.121285140562249	89.8311357572469\\
0.122088353413655	92.0280639234369\\
0.12289156626506	94.1481461951153\\
0.123694779116466	96.1935435536518\\
0.124497991967871	98.1667447765368\\
0.125301204819277	100.070432538481\\
0.126104417670683	101.90738732846\\
0.126907630522088	103.680418960543\\
0.127710843373494	105.392318185719\\
0.1285140562249	107.045822979321\\
0.129317269076305	108.643595597415\\
0.130120481927711	110.18820751884\\
0.130923694779116	111.682130179015\\
0.131726907630522	113.127729805607\\
0.132530120481928	114.527265435622\\
0.133333333333333	115.882888989451\\
0.134136546184739	117.196646970113\\
0.134939759036145	118.470483203762\\
0.13574297188755	119.706242287484\\
0.136546184738956	120.905673560633\\
0.137349397590361	122.070435280127\\
0.138152610441767	123.202098942812\\
0.138955823293173	124.302153645907\\
0.139759036144578	125.37201040385\\
0.140562248995984	126.413006343686\\
0.14136546184739	127.426408779055\\
0.142168674698795	128.413419121583\\
0.142971887550201	129.375176642585\\
0.143775100401606	130.312762015864\\
0.144578313253012	131.227200728389\\
0.145381526104418	132.119466264145\\
0.146184738955823	132.990483133855\\
0.146987951807229	133.841129703787\\
0.147791164658635	134.672240818049\\
0.14859437751004	135.484610318079\\
0.149397590361446	136.27899334885\\
0.150200803212851	137.056108592308\\
0.151004016064257	137.816640326721\\
0.151807228915663	138.561240415799\\
0.152610441767068	139.290530113105\\
0.153413654618474	140.005101809458\\
0.15421686746988	140.705520637754\\
0.155020080321285	141.392325970458\\
0.155823293172691	142.06603285151\\
0.156626506024096	142.727133302356\\
0.157429718875502	143.376097584851\\
0.158232931726908	144.013375358415\\
0.159036144578313	144.639396779709\\
0.159839357429719	145.254573533157\\
0.160642570281124	145.859299808733\\
0.16144578313253	146.453953201144\\
0.162248995983936	147.038895578644\\
0.163052208835341	147.614473880633\\
0.163855421686747	148.181020880675\\
0.164658634538153	148.738855899457\\
0.165461847389558	149.288285458675\\
0.166265060240964	149.829603933759\\
0.167068273092369	150.363094130067\\
0.167871485943775	150.889027862351\\
0.168674698795181	151.407666477779\\
0.169477911646586	151.919261365285\\
0.170281124497992	152.424054429532\\
0.171084337349398	152.922278537996\\
0.171887550200803	153.414157947633\\
0.172690763052209	153.89990870736\\
0.173493975903614	154.379739029786\\
0.17429718875502	154.853849657969\\
0.175100401606426	155.322434198753\\
0.175903614457831	155.785679446657\\
0.176706827309237	156.2437656875\\
0.177510040160643	156.696866986266\\
0.178313253012048	157.145151463831\\
0.179116465863454	157.588781551262\\
0.179919678714859	158.027914238819\\
0.180722891566265	158.462701310308\\
0.181526104417671	158.893289565039\\
0.182329317269076	159.319821025575\\
0.183132530120482	159.742433144812\\
0.183935742971888	160.161258992199\\
0.184738955823293	160.57642743929\\
0.185542168674699	160.988063339846\\
0.186345381526104	161.396287698211\\
0.18714859437751	161.801217829057\\
0.187951807228916	162.202967527274\\
0.188755020080321	162.601647212795\\
0.189558232931727	162.997364075609\\
0.190361445783133	163.390222210506\\
0.191164658634538	163.780322733566\\
0.191967871485944	164.167763885185\\
0.192771084337349	164.552641121126\\
0.193574297188755	164.935047200015\\
0.194377510040161	165.315072268944\\
0.195180722891566	165.692803953674\\
0.195983935742972	166.068327452737\\
0.196787148594378	166.441725633298\\
0.197590361445783	166.813079125938\\
0.198393574297189	167.182466411622\\
0.199196787148594	167.549963905135\\
0.2	167.915646038059\\
};
\addlegendentry{$P=1$};
\addplot [color=mycolor3, line width = 1.5pt, line join = round]
  table[row sep=crcr]{%
0	282.020146247426\\
0.000803212851405623	282.022373515107\\
0.00160642570281125	282.029056394511\\
0.00240963855421687	282.040198123721\\
0.00321285140562249	282.055804103052\\
0.00401606425702811	282.075881909756\\
0.00481927710843374	282.100441314214\\
0.00562248995983936	282.12949428805\\
0.00642570281124498	282.163055024122\\
0.0072289156626506	282.201139979326\\
0.00803212851405622	282.243767895751\\
0.00883534136546185	282.29095982467\\
0.00963855421686747	282.342739189573\\
0.0104417670682731	282.399131799743\\
0.0112449799196787	282.460165925913\\
0.0120481927710843	282.525872345388\\
0.01285140562249	282.596284380408\\
0.0136546184738956	282.671438000367\\
0.0144578313253012	282.75137186076\\
0.0152610441767068	282.836127385963\\
0.0160642570281124	282.925748861243\\
0.0168674698795181	283.020283505233\\
0.0176706827309237	283.11978157985\\
0.0184738955823293	283.224296481179\\
0.0192771084337349	283.333884854958\\
0.0200803212851406	283.448606712031\\
0.0208835341365462	283.56852555588\\
0.0216867469879518	283.693708522741\\
0.0224899598393574	283.824226513378\\
0.0232931726907631	283.960154361542\\
0.0240963855421687	284.10157099461\\
0.0248995983935743	284.248559615598\\
0.0257028112449799	284.401207883011\\
0.0265060240963855	284.559608135098\\
0.0273092369477912	284.723857580001\\
0.0281124497991968	284.894058559543\\
0.0289156626506024	285.070318768606\\
0.029718875502008	285.252751558148\\
0.0305220883534137	285.441476183863\\
0.0313253012048193	285.636618139621\\
0.0321285140562249	285.838309484325\\
0.0329317269076305	286.046689183866\\
0.0337349397590361	286.261903500115\\
0.0345381526104418	286.484106398447\\
0.0353413654618474	286.713459970364\\
0.036144578313253	286.950134922486\\
0.0369477911646586	287.194311080197\\
0.0377510040160643	287.446177916736\\
0.0385542168674699	287.705935151257\\
0.0393574297188755	287.973793381268\\
0.0401606425702811	288.249974754813\\
0.0409638554216867	288.534713715759\\
0.0417670682730924	288.828257800447\\
0.042570281124498	289.130868477494\\
0.0433734939759036	289.442822111156\\
0.0441767068273092	289.764410934851\\
0.0449799196787149	290.095944161699\\
0.0457831325301205	290.437749172875\\
0.0465863453815261	290.790172801548\\
0.0473895582329317	291.153582728451\\
0.0481927710843374	291.528369018194\\
0.048995983935743	291.914945777081\\
0.0497991967871486	292.313752963496\\
0.0506024096385542	292.725258385814\\
0.0514056224899598	293.149959851106\\
0.0522088353413655	293.588387556361\\
0.0530120481927711	294.041106682626\\
0.0538152610441767	294.508720258513\\
0.0546184738955823	294.991872301099\\
0.055421686746988	295.491251272393\\
0.0562248995983936	296.007593902146\\
0.0570281124497992	296.541689391392\\
0.0578313253012048	297.094384068431\\
0.0586345381526104	297.666586532688\\
0.0594377510040161	298.259273380673\\
0.0602409638554217	298.873495563576\\
0.0610441767068273	299.510385404973\\
0.0618473895582329	300.171164519196\\
0.0626506024096386	300.857152538439\\
0.0634538152610442	301.569776973301\\
0.0642570281124498	302.310584160671\\
0.0650602409638554	303.081251619559\\
0.0658634538152611	303.883601894971\\
0.0666666666666667	304.719618145051\\
0.0674698795180723	305.591461741942\\
0.0682730923694779	306.501492142984\\
0.0690763052208835	307.452289379555\\
0.0698795180722892	308.446679550215\\
0.0706827309236948	309.487763786021\\
0.0714859437751004	310.578951193788\\
0.072289156626506	311.723996379333\\
0.0730923694779116	312.927042262791\\
0.0738955823293173	314.192668974711\\
0.0746987951807229	315.525949753734\\
0.0755020080321285	316.932514881591\\
0.0763052208835341	318.418624771327\\
0.0771084337349398	319.991253535466\\
0.0779116465863454	321.658184331358\\
0.078714859437751	323.428117760941\\
0.0795180722891566	325.310794698123\\
0.0803212851405622	327.31713406621\\
0.0811244979919679	329.459386008858\\
0.0819277108433735	331.75129861848\\
0.0827309236947791	334.208294495816\\
0.0835341365461847	336.847648306957\\
0.0843373493975904	339.688649157274\\
0.085140562248996	342.752718530876\\
0.0859437751004016	346.063432883941\\
0.0867469879518072	349.646363985851\\
0.0875502008032129	353.528590255326\\
0.0883534136546185	357.73763265079\\
0.0891566265060241	362.299403469602\\
0.0899598393574297	367.234483995497\\
0.0907630522088353	372.551602965049\\
0.091566265060241	378.236482554184\\
0.0923694779116466	384.233151461686\\
0.0931726907630522	390.413373651588\\
0.0939759036144578	396.528340737002\\
0.0947791164658635	402.136701941204\\
0.0955823293172691	406.508644872265\\
0.0963855421686747	408.526590949644\\
0.0971887550200803	406.652563051038\\
0.097991967871486	399.10718746195\\
0.0987951807228916	384.427180802632\\
0.0995983935742972	362.354674435077\\
0.100401606425703	334.528135802768\\
0.101204819277108	304.24512642374\\
0.102008032128514	275.261519035775\\
0.10281124497992	250.488524557549\\
0.103614457831325	231.374004458439\\
0.104417670682731	218.021657987809\\
0.105220883534137	209.677563237965\\
0.106024096385542	205.239690058408\\
0.106827309236948	203.615086499295\\
0.107630522088353	203.894885163333\\
0.108433734939759	205.396914292149\\
0.109236947791165	207.640991502848\\
0.11004016064257	210.302982821975\\
0.110843373493976	213.170430764391\\
0.111646586345382	216.107388474572\\
0.112449799196787	219.029004987468\\
0.113253012048193	221.884023707697\\
0.114056224899598	224.643040227456\\
0.114859437751004	227.290731799696\\
0.11566265060241	229.820751892163\\
0.116465863453815	232.232387636769\\
0.117269076305221	234.528376786506\\
0.118072289156627	236.713485620057\\
0.118875502008032	238.793586587728\\
0.119678714859438	240.775064610033\\
0.120481927710843	242.664439727339\\
0.121285140562249	244.468132292053\\
0.122088353413655	246.192321936517\\
0.12289156626506	247.842868142672\\
0.123694779116466	249.425270814734\\
0.124497991967871	250.944656702997\\
0.125301204819277	252.405782104497\\
0.126104417670683	253.813045523424\\
0.126907630522088	255.170506116253\\
0.127710843373494	256.481905141401\\
0.1285140562249	257.750688603178\\
0.129317269076305	258.980029945137\\
0.130120481927711	260.172852076747\\
0.130923694779116	261.331848354168\\
0.131726907630522	262.459502146151\\
0.132530120481928	263.558105112427\\
0.133333333333333	264.629773927982\\
0.134136546184739	265.676465680795\\
0.134939759036145	266.699991894052\\
0.13574297188755	267.702031291795\\
0.136546184738956	268.684141449495\\
0.137349397590361	269.647769315951\\
0.138152610441767	270.594260777719\\
0.138955823293173	271.524869342064\\
0.139759036144578	272.440764016214\\
0.140562248995984	273.343036433871\\
0.14136546184739	274.23270733823\\
0.142168674698795	275.110732459414\\
0.142971887550201	275.978007872124\\
0.143775100401606	276.835374833468\\
0.144578313253012	277.683624225352\\
0.145381526104418	278.523500563446\\
0.146184738955823	279.355705683434\\
0.146987951807229	280.180902098434\\
0.147791164658635	280.999716025474\\
0.14859437751004	281.812740234618\\
0.149397590361446	282.620536565749\\
0.150200803212851	283.42363831671\\
0.151004016064257	284.222552360271\\
0.151807228915663	285.017761168857\\
0.152610441767068	285.809724586802\\
0.153413654618474	286.598881542436\\
0.15421686746988	287.385651582807\\
0.155020080321285	288.170436289883\\
0.155823293172691	288.953620615157\\
0.156626506024096	289.735574075667\\
0.157429718875502	290.516651904245\\
0.158232931726908	291.297196085006\\
0.159036144578313	292.077536325281\\
0.159839357429719	292.857990956372\\
0.160642570281124	293.638867787721\\
0.16144578313253	294.420464866719\\
0.162248995983936	295.203071231742\\
0.163052208835341	295.986967580951\\
0.163855421686747	296.772426915858\\
0.164658634538153	297.559715151104\\
0.165461847389558	298.349091649561\\
0.166265060240964	299.140809783655\\
0.167068273092369	299.935117402073\\
0.167871485943775	300.732257322901\\
0.168674698795181	301.532467759215\\
0.169477911646586	302.33598274648\\
0.170281124497992	303.14303253734\\
0.171084337349398	303.953843972149\\
0.171887550200803	304.768640843111\\
0.172690763052209	305.587644237502\\
0.173493975903614	306.411072847869\\
0.17429718875502	307.239143300432\\
0.175100401606426	308.072070439253\\
0.175903614457831	308.910067615005\\
0.176706827309237	309.753346961407\\
0.177510040160643	310.602119649192\\
0.178313253012048	311.456596154567\\
0.179116465863454	312.316986486463\\
0.179919678714859	313.183500431992\\
0.180722891566265	314.0563477828\\
0.181526104417671	314.935738558611\\
0.182329317269076	315.82188320716\\
0.183132530120482	316.714992828679\\
0.183935742971888	317.615279357989\\
0.184738955823293	318.522955761175\\
0.185542168674699	319.438236233286\\
0.186345381526104	320.361336380204\\
0.18714859437751	321.292473390818\\
0.187951807228916	322.231866254567\\
0.188755020080321	323.179735952972\\
0.189558232931727	324.136305694559\\
0.190361445783133	325.101801193786\\
0.191164658634538	326.076450993233\\
0.191967871485944	327.060486842398\\
0.192771084337349	328.054144132073\\
0.193574297188755	329.057662340743\\
0.194377510040161	330.071285460358\\
0.195180722891566	331.095262318711\\
0.195983935742972	332.129846751046\\
0.196787148594378	333.175297598916\\
0.197590361445783	334.231878567114\\
0.198393574297189	335.299858011087\\
0.199196787148594	336.379508788471\\
0.2	337.471108256489\\
};
\addlegendentry{$P=10$}
\end{axis}
\end{tikzpicture}%